\documentclass[12pt,reqno]{amsart}

\usepackage[all]{xy}
\usepackage{amsmath}
\usepackage{amsfonts}
\usepackage{amssymb}
\usepackage[all]{xy}
\usepackage{mathrsfs}
\usepackage{bbding}
\usepackage{txfonts}
\usepackage{amscd}
\usepackage{amsthm}

\usepackage{relsize}

\usepackage{dutchcal}

\usepackage{dsfont}

\usepackage[
    colorlinks=true,
    linkcolor=blue,
    citecolor=blue,
    urlcolor=blue,
    filecolor=blue,
    final,
    backref=page,
    hyperindex
]{hyperref}

\usepackage[shortlabels]{enumitem}
\usepackage{ifpdf}
\ifpdf
\usepackage{hyperref}
\else
\usepackage{hyperref}
\fi
\usepackage{tikz}
\usepackage[active]{srcltx}

\def\deg{{\rm deg}}

\def\M{M}
\def\Z{\mathbb{Z}}
\def\N{\mathbb{N}}

\def\C{\mathbb{C}}

\numberwithin{equation}{section}
\newtheorem{theo}{Theorem}[section]
\newtheorem{defi}[theo]{Definition}
\newtheorem{coro}[theo]{Corollary}
\newtheorem{lemm}[theo]{Lemma}
\newtheorem{prop}[theo]{Proposition}

\newtheorem{case}{Case}

\newtheorem{rema}[theo]{Remark}

\begin{document}

\title[Representations of the Loop mirror BMS-Kac-Moody algebra]{Representations of the loop mirror BMS-Kac-Moody algebra}

\author{Haibo Chen}
\address{School of Science, Jimei University, Xiamen, Fujian 361021, China}
\email{hypo1025@jmu.edu.cn}

\author{Yongliang Wang}
\address{School of Science, Hefei University of Technology, Hefei, Anhui 230601, China}
\email{2026010019@mail.hfut.edu.cn}

\author{Ying Xu$^*$}
\address{School of Science, Hefei University of Technology, Hefei, Anhui 230601, China}
\email{xuying@hfut.edu.cn}

\date{\today}

\vspace{-5mm}

\begin{abstract}
In this paper, we  construct a new family of non-weight modules
$
\Omega(\lambda,\zeta,\alpha,\beta,\gamma,\sigma,\linebreak h)
$ over the loop mirror BMS-Kac-Moody algebra $\mathcal{G}$, which is a
Takiff-type extension of the loop mirror Heisenberg--Virasoro algebra. We establish criteria for  irreducibility, determine their isomorphism classes, and classify all $\mathcal{G}$-modules that are free of rank one. Moreover, we determine all self-extensions of these modules: for every such module $M$, $\operatorname{Ext}^{1}_{U(\mathcal{G})}(M,M)\cong \mathbb{C}^{6}$, and we construct explicit $1$-cocycles realizing every cohomology class. Every nonsplit self-extension of a simple module is indecomposable of composition length two. Restricting these modules to the loop mirror Heisenberg--Virasoro subalgebra $\mathcal{H}$ yields new $\mathcal{H}$-modules. Among them the $\sigma$-type modules $\Psi(\lambda,\zeta,\sigma,\alpha,\beta,\gamma,h)$ constitute a new family with no counterpart in the non-loop or non-Takiff case. For the degenerate case $\sigma=0$, we establish criteria for irreducibility and isomorphism, and describe certain maximal submodules and their simple quotients. Finally, we construct and classify a new class of tensor product $\mathcal{H}$-modules.

\end{abstract}


\footnote {
$^*$Corresponding author: Y. Xu (xuying@hfut.edu.cn).
}
\maketitle
\vspace{-10mm}

\tableofcontents
\section{Introduction}
The representation theory of infinite-dimensional Lie algebras---notably the Virasoro, Heisenberg--Virasoro, and BMS algebras, together with their super and loop analogues---has long occupied a central position in algebra, geometry, and mathematical physics, owing largely to its profound connections with conformal field theory, vertex operator algebras, and string theory.

Recently, attention has increasingly shifted toward non-weight modules,
 which lack a Cartan weight-space decomposition and exhibit richer, more intricate structures whose classification remains largely open. 
 Such modules were investigated for \(\mathfrak{sl}_{n+1}\) by Nilsson \cite{N1} and independently by Tan and Zhao \cite{TZ1}. Concurrently, substantial progress has been made in the context of loop algebras: Chen and Yao \cite{CY} investigated non-weight modules over loop-Virasoro and Block-type Lie algebras, deriving simple loop-Virasoro modules by tensoring irreducible modules with highest-weight or Whittaker modules; Chen and Guo \cite{CG1} constructed and classified  free \(U(\mathcal{h})\)-modules of rank one over the BMS--Kac--Moody algebra modulo the center.
 In a related direction, Takiff-type extensions have attracted considerable attention since Takiff \cite{T} introduced them in invariant theory. For instance, Qiao and Zhu \cite{QZ} constructed non-weight modules over the Takiff Lie algebra for \(\mathcal{sl}_2\) via tensor products of rank-one free \(U(\mathcal{h})\)-modules with irreducible highest-weight modules, completely characterizing their irreducibility and isomorphism conditions. Collectively, these developments demonstrate that loop and Takiff directions constitute a rapidly expanding frontier in the non-weight representation theory of infinite-dimensional Lie algebras. Moreover, these non-weight modules were used to study tensor product modules involving restricted modules over  Virasoro algebras, further enriching the representation theory of infinite-dimensional Lie algebras (see, e.g., \cite{MZ,TZ2,TZ3}). Motivated by these developments, Chen and Bai \cite{CB} constructed non-weight modules by tensoring finitely many irreducible polynomial   modules with irreducible restricted modules.

However, the non-weight representation theory of the loop mirror BMS--Kac--Moody algebra \(\mathcal{G}\) has remained completely unexplored. Unlike its classical cases, \(\mathcal{G}\) is endowed with a natural \(\frac{1}{2}\mathbb{Z}\)-grading in which  each infinite-dimensional homogeneous component is labelled by two loop indexes. By virtue of being a Takiff-type extension of the loop mirror Heisenberg--Virasoro algebra, it exhibits phenomena with no analogue in the one-parameter or non-loop cases. Moreover, such as extensions of non-weight modules of these algebras, have rarely been addressed in the literatures.

This paper investigates the non-weight representations of the loop mirror BMS--Kac--Moody algebra \(\mathcal{G}\). We 
 construct a family of non-weight \(\mathcal{G}\)-modules \(\Omega(\lambda,\zeta,\alpha,\beta,\gamma,\sigma,h)\), classify all rank-one free \(U(\mathcal{h})\)-modules over \(\mathcal{G}\), and determine their irreducibility and isomorphism classes. Restricting these modules to the loop mirror Heisenberg--Virasoro subalgebra \(\mathcal{H}\) yields \(\sigma\)-type and degenerate \(\mathcal{H}\)-modules whose irreducibility, maximal submodules, quotients, and isomorphism classes are studied. We further construct and classify tensor product modules obtained by tensoring finitely many irreducible degenerate \(\mathcal{H}\)-modules with simple restricted \(\mathcal{H}\)-modules. These results fill a gap in the literature, reveal phenomena absent in the non-loop case, and extend earlier work on the BMS algebra, the mirror Heisenberg--Virasoro algebra, the loop Virasoro algebra, and the BMS--Kac--Moody algebra \cite{CB,CG,CY,GMZ}.

The main conclusions of this paper are summarized as follows:

\begin{itemize}
    \item \textbf{ Non-weight \(\mathcal{G}\)-modules (Section \ref{section 3}).} \\
    The family \(\Omega(\lambda,\zeta,\alpha,\beta,\gamma,\sigma,h)\) is constructed on \(\mathbb{C}[x,y]\) (Proposition \ref{prop 3.1}); it is simple if and only if \((\alpha,\gamma)\neq(0,0)\) (Theorem \ref{theo 3.3}). Their isomorphism classes are determined (Theorem \ref{theo3.2}), and they yield a complete classification of all rank-one free \(U(\mathfrak{h})\)-modules over \(\mathcal{G}\), with \(\mathfrak{h}=\mathbb{C}L_0\oplus\mathbb{C}W_0\) (Theorem \ref{th4.1}).
    
    \item \textbf{Self-extensions of \(\mathcal{G}\)-modules (Section \ref{section 4}).} \\
    A family of explicit cocycles \(\tau_{\boldsymbol{v}}\) is constructed, and computing the first cohomology group shows that every \(1\)-cocycle is cohomologous to some \(\tau_{\boldsymbol{v}}\), so that \(H^1(\mathcal{G},\operatorname{End}_{\mathbb{C}}(M))\cong \operatorname{Ext}^1_{U(\mathcal{G})}(M,M)\cong \mathbb{C}^6\) (Theorem \ref{theo 4.1}); if $(\alpha,\gamma)\neq (0,0)$,  every nonzero class gives a nonsplit self-extension that is indecomposable of composition length two (Corollary \ref{coro 4.5}).

    \item \textbf{New simple \(\mathcal{H}\)-modules (Section \ref{section 5}).} \\
    Restricting the \(\mathcal{G}\)-modules to the loop mirror Heisenberg--Virasoro subalgebra \(\mathcal{H}\) yields the \(\sigma\)-type modules \(\Psi(\lambda,\zeta,\sigma,\alpha,\beta,\gamma,h)\) \((\sigma\neq0)\), simple if and only if \((\alpha,\gamma)\neq(0,0)\) (Proposition \ref{prop 5.1}), and the degenerate modules \(\Omega_{\mathcal{H}}(\lambda,\zeta,\alpha,\beta,\gamma,h)\) \((\sigma=0)\), simple if and only if \(\deg_y(h(y))=1\) and \((\alpha,\gamma)\neq(0,0)\) (Theorem \ref{prop 6.1}); certain maximal submodules, simple quotients, and isomorphism classes are determined (Propositions \ref{prop 5.2} and \ref{prop 5.3}, Theorem \ref{theo 5.2} and Corollary \ref{coro 6.5}).

    \item \textbf{Tensor product modules (Section \ref{section 6}).} \\
    For \(\mathrm{T} = \left(\bigotimes_{i=1}^{m} \Omega_{\mathcal{H}}(\lambda_i, \zeta_i, \alpha_i, \beta_i, \gamma_i)\right) \otimes V\), where \(V\) is a simple restricted \(\mathcal{H}\)-module, \(\mathrm{T}\) is irreducible if and only if the pairs \((\lambda_i,\zeta_i)\) are pairwise distinct (Theorem \ref{th5.1} and Corollary \ref{coro6.6}); the isomorphism classes of such modules are completely determined (Theorem \ref{theo5.3}).
\end{itemize}


Throughout this paper, $\C$, $\C^*$, $\Z$, $\Z^*$, $\N$, and $\Z_{+}$ denote the sets of complex numbers, nonzero complex numbers, integers, nonzero  integers, nonnegative integers, and positive integers, respectively. All vector spaces (respectively, algebras and modules) are assumed to be over $\C$. The zero operator and the identity operator on any vector space $V$ are denoted by $0_V$ and $\operatorname{Id}_{V}$, respectively.

\section{Preliminaries}  
This section reviews the definitions of the loop mirror BMS-Kac-Moody algebra and summarizes known results concerning the BMS algebra and the mirror Heisenberg-Virasoro algebra.

\subsection{\texorpdfstring{Loop mirror BMS-Kac-Moody algebra $\mathcal{G}$}{1}}

For a complex Lie algebra $\mathcal{L}$ and $n\in\mathbb Z_+$, define the associated Takiff algebra by $\mathcal L_n=\mathcal L\otimes(\mathbb C[t]/ t^{n+1}\mathbb C[t])$. Following the pioneering work \cite{T}, in which these algebras were studied for $n=1$, the Lie algebra $\mathcal{L}_{n}$ is also known as the generalized or $n$-th Takiff algebra. When $n=1$, we refer to $\mathcal L_1$ simply as the Takiff algebra.

Let $\mathcal{G}$ denote the loop mirror BMS-Kac-Moody algebra. It has a $\mathbb{C}$-basis
\[
\left\{ L_{m,i},\; W_{m,i},\; I_{r,i},\; B_{r,i},\; C_{k,i} \;\middle|\; m,i\in \mathbb{Z},\; r\in\mathbb{Z}+\tfrac{1}{2},\; k=1,2,3,4 \right\},
\]
with the following Lie brackets:
\begin{align*}
&[L_{m,i},L_{n,j}]= (m-n)L_{m+n,i+j}+\frac{m^3-m}{12}\delta_{m+n,0}\,C_{1,i+j}, \notag \\
&[L_{m,i},W_{n,j}]=(m-n)W_{m+n,i+j}+\frac{m^3-m}{12}\delta_{m+n,0}\,C_{2,i+j}, \notag \\
&[L_{m,i},I_{r,j}]=-r\,I_{m+r,i+j}, \notag \ 
[L_{m,i},B_{r,j}]=-r\,B_{m+r,i+j}, \notag \ 
[W_{m,i},I_{r,j}]=-r\,B_{m+r,i+j}, \notag \\
&[I_{r,i},I_{s,j}]=r\,\delta_{r+s,0}\,C_{3,i+j}, \notag \ 
[I_{r,i},B_{s,j}]=r\,\delta_{r+s,0}\,C_{4,i+j}, \notag \\
&[W_{m,i},W_{n,j}]=[W_{m,i},B_{r,j}]=[B_{r,i},B_{s,j}]=[C_{k,i},\mathcal{G}]=0
\end{align*}
for all $m,n,i,j\in \mathbb{Z}$, $r,s\in\mathbb{Z}+\frac{1}{2}$, and $k\in\{1,2,3,4\}$. 
The algebra $\mathcal{G}$ admits a $\frac{1}{2}\Z$-grading 
\begin{equation}\label{grading 1.1}
\aligned
\mathcal{G}=\bigoplus_{n\in\frac{1}{2}\Z}\mathcal{G}_n=\bigoplus_{n\in\Z}\Big\{L_{n,i},W_{n,i},\delta_{n,0}C_{k,i}\mid i\in \Z,k=1,2,3,4\Big\}\oplus\bigoplus_{n\in\frac{1}{2}+\Z}\Big\{I_{n,i},B_{n,i}\mid i\in \Z\Big\}.
\endaligned
\end{equation}

The subalgebras $\mathcal{V}=\{L_{m,i}, C_{1,i}\mid m,i\in \mathbb{Z}\}$, $\mathcal{W}=\{L_{m,i},W_{m,i}, C_{1,i},C_{2,i}\mid m,i\in \mathbb{Z}\}$, and $\mathcal{H}=\{L_{m,i}, I_{r,i},C_{1,i},C_{3,i}\mid m,i\in \mathbb{Z}, r\in\mathbb{Z}+\frac{1}{2}\}$ are isomorphic to the loop Virasoro algebra, the loop BMS algebra, and the loop mirror Heisenberg-Virasoro algebra, respectively. Moreover, $\mathcal{G}$ is a Takiff-type extension of $\mathcal{H}$. For simplicity, we denote $X_p:=X_{p,0}$, where $X\in\{L,W,I,B,C\}$ and $p\in\frac{1}{2}\mathbb{Z}$.
The algebra $\mathcal{G}$ also contains several classical Lie algebras: the Virasoro algebra $\hat{\mathcal{V}}=\{L_m, C_1\mid m\in \mathbb{Z}\}$, the BMS algebra $\hat{\mathcal{W}}=\{L_{m},W_m, C_{1},C_2\mid m\in \mathbb{Z}\}$, and the mirror Heisenberg-Virasoro algebra $\hat{\mathcal{H}}=\{L_{m}, I_r,C_{1},C_3\mid m\in \mathbb{Z}, r\in\mathbb{Z}+\frac{1}{2}\}$. The so-called mirror BMS-Kac-Moody algebra $\hat{\mathcal{G}}=\{L_m,W_m,I_r,B_r,C_i\mid m\in\mathbb{Z},r\in\Z+\frac{1}{2},i=1,2,3,4\}$ is also a subalgebra of $\mathcal{G}$.

\subsection{Some known results}
Fix any $\lambda\in\C^*, \alpha, \beta\in\C$ and $h(y)=\sum_{i=0}^{\infty}\xi_iy^i\in\C[y]$. For any $r\in\frac{1}{2}\Z,m\in\Z$,  define
\begin{equation}\label{formula 2.1}
\aligned
g_{r}(y)=r\frac{h(y)-h(\alpha)}{y-\alpha},
\endaligned
\end{equation}
and
\begin{equation}\label{formula 2.2}
\aligned
h_{m}(y)&=mh(y)-\alpha(m-1) g_m(y).
\endaligned
\end{equation}
Then we have
\begin{equation}\label{formula 2.3}
\aligned
mh_m(y)-nh_n(y)
+m(y-m\alpha)\frac{d}{dy}h_n(y)
-n(y-n\alpha)\frac{d}{dy}h_m(y)
=(m-n)h_{m+n}(y),
\endaligned
\end{equation}
and
\begin{equation}\label{formula 2.4}
\aligned
m(y-m\alpha)\frac{d}{dy}g_r(y)-rg_r(y)
-r\frac{d}{dy}h_m(y)
=-rg_{m+r}(y).
\endaligned
\end{equation}
 
Let $\Psi_{\hat{\mathcal{W}}}(\lambda,\alpha,h)$
denote the vector space $\mathbb C[x,y]$, endowed with the following $\hat{\mathcal W}$-action:
\begin{eqnarray*}
&&L_m f(x,y)=\lambda^m(x+h_{m}(y))f(x+m,y)+m\lambda^m(y-m\alpha)\frac{\partial }{\partial y}f(x+m,y), 
\\&& 
W_m f(x,y)=\lambda^m(y-m\alpha)f(x+m,y),
\ 
C_1 f(x,y) = C_2 f(x,y) =0,\label{formular 2.1.3.3}
\end{eqnarray*}
where $f(x,y)\in\C[x,y]$ and $m\in\Z$.

We recall several established properties of $\Psi_{\hat{\mathcal{W}}}(\lambda,\alpha,h)$.
\begin{lemm}[\rm see \cite{CG}]\label{le2.1}
Let $\lambda,\lambda^{\prime}\in\mathbb C^*$, $\alpha,\alpha^{\prime}\in\mathbb C$, and $h(y),h^{\prime}(y)\in\mathbb C[y]$, and let $h_m(y),h^{\prime}_m(y)$ be defined by (\ref{formula 2.2}). Then the following statements hold.
\begin{itemize}
\item[\rm(1)] $\Psi_{\hat{\mathcal{W}}}(\lambda,\alpha,h)$ is a $\hat{\mathcal{W}}$-module.
\item[\rm(2)] $\Psi_{\hat{\mathcal{W}}}(\lambda,\alpha,h)$ is simple if and only if $\alpha\neq0$.
\item[\rm(3)] $\Psi_{\hat{\mathcal{W}}}(\lambda,\alpha,h)\cong \Psi_{\hat{\mathcal{W}}}(\lambda^{\prime},\alpha^{\prime},h^{\prime})$ if and only if $\lambda=\lambda^{\prime},\alpha=\alpha^{\prime},h(y)=h^{\prime}(y)$.
\item[\rm(4)] 
Every $\hat{\mathcal{W}}$-module that is free of rank one over $U(\mathbb C L_0\oplus\mathbb C W_0)$ is isomorphic to $\Psi_{\hat{\mathcal{W}}}(\lambda,\alpha,h)$ for some $\lambda\in\mathbb C^*$, $\alpha\in\mathbb C$, and $h(y)\in\mathbb C[y]$.
\end{itemize}
\end{lemm}
Note that the underlying Virasoro module of the $\hat{\mathcal{W}}$-module $\Psi_{\hat{\mathcal{W}}}(\lambda,\alpha,h)$ was studied in \cite{CG2}.
Let $\Psi_{\hat{\mathcal{H}}}(\lambda,\alpha,\beta)=\C[x]$
be a vector space and define an action of $\hat{\mathcal{H}}$ by
\begin{equation*}
\aligned
L_mf(x)&=\lambda^m(x+m\alpha)f(x+m),\ 
I_rf(x)=\beta\lambda^{r-\frac{1}{2}}f(x+r),\
C_1 f(x)= C_3 f(x) =0
\endaligned
\end{equation*}
for all $m\in\Z$, $r\in\Z+\frac{1}{2}$, and $f(x)\in\C[x]$. The following results concerning $\Psi_{\hat{\mathcal{H}}}(\lambda,\alpha,\beta)$ have been established.

\begin{lemm}[\rm see \cite{GMZ}]\label{le2.1.1}
For $\lambda\in\C^*$ and $\alpha,\beta\in\C$,
\begin{itemize}
\item[\rm(1)] $\Psi_{\hat{\mathcal{H}}}(\lambda,\alpha,\beta)$ is an $\hat{\mathcal{H}}$-module;
\item[\rm(2)] $\Psi_{\hat{\mathcal{H}}}(\lambda,\alpha,\beta)$ is irreducible if and only if $(\alpha, \beta)\neq(0,0)$.
\end{itemize}
\end{lemm}

We conclude this subsection with the following technical result.
\begin{lemm}[\rm see \cite{LZ0}]\label{le2.3}
Let $P$ be a complex vector space, and let $P_1\subseteq P$ be a subspace. Suppose that $\lambda_1,\lambda_2,\dots,\lambda_s\in\C^*$ are pairwise distinct and that $v_{i,j}\in P$ and $f_{i,j}(t)\in\C[t]$ satisfy ${\rm deg}(f_{i,j}(t))=j$ for $i=1,2,\dots,s;j=0,1,2,\dots,k$. If $\sum_{i=1}^{s}\sum_{j=0}^{k}\lambda_i^mf_{i,j}(m)v_{i,j}\in P_1$
for every $m\in\mathbb Z$ 
,  then $v_{i,j}\in P_1$ for all $i, j$.

\end{lemm}

\section{\texorpdfstring{Non-weight $\mathcal{G}$-modules $\Omega(\lambda, \zeta, \alpha,\beta,\gamma,\sigma,h)$}{Non-weight modules over G}}\label{section 3}

In this section, we construct a family of non-weight $\mathcal G$-modules, study their isomorphism classes and simplicity, and use them to classify $U(\mathcal{h})$-free modules over $\mathcal{G}$.

\subsection{Construction and properties}
\begin{prop}\label{prop 3.1}
  Let $\lambda,\zeta\in\mathbb C^*$, $\alpha,\beta,\gamma,\sigma\in\mathbb C$, and $h(y)\in\mathbb C[y]$. Set $\Omega(\lambda,\zeta,\alpha,\beta,\gamma,\sigma,h)\cong\mathbb C[x,y]$ as a vector space and define a 
 ${\mathcal{G}}$-action on $\Omega(\lambda, \zeta, \alpha,\beta,\gamma,\sigma,h)$ by
\begin{equation*}
\aligned
&L_{m,i} f(x,y)=\lambda^m\zeta^i(x+h_{m}(y)+i\sigma(y-m\alpha))f(x+m,y)+m\lambda^m\zeta^i(y-m\alpha)\frac{\partial }{\partial {y}}{f(x+m,y)},\\
&W_{m,i} f(x,y)=\lambda^m\zeta^i(y-m\alpha)f(x+m,y),\\
&I_{r,i} f(x,y)=\lambda^{r-\frac{1}{2}}\zeta^i(\beta+\gamma g_{r}(y)+i\sigma \gamma)f(x+r,y)+r\gamma\lambda^{r-\frac{1}{2}}\zeta^i\frac{\partial }{\partial {y}}{f(x+r,y)},\\ 
&B_{r,i} f(x,y)=\gamma\lambda^{r-\frac{1}{2}}\zeta^if(x+r,y),\ C_{k,i} f(x,y)=0,
\endaligned
\end{equation*}
where $m,i\in\Z$, $r\in\Z+\frac{1}{2}$, and $k\in\{1,2,3,4\}$. These formulas define a $\mathcal G$-module structure on $\Omega(\lambda,\zeta,\alpha,\beta,\gamma,\sigma,h)$.
\end{prop}
\begin{proof}
For any $m,n,i,j\in\mathbb{Z}$ and $f(x,y)\in\mathbb{C}[x,y]$, $L_{m,i}L_{n,j}f(x,y)=\lambda^{m+n}\zeta^{i+j}(Af(x+m+n,y)+B\frac{\partial }{\partial y}f(x+m+n,y)+C\frac{\partial^2}{\partial y^2}f(x+m+n,y))$, where $A=x+h_m(y)+i\sigma(y-m\alpha)\bigl(x+m+h_n(y)+j\sigma(y-n\alpha)+m(y-m\alpha)
\left(\frac{d}{dy}h_n(y)+j\sigma\right)\bigr)$, $B=
n(y-n\alpha)(x+h_m(y)+i\sigma(y-m\alpha)\bigr)+m(y-m\alpha)
\bigl(x+m+h_n(y)+j\sigma(y-n\alpha)\bigr)+mn(y-m\alpha)$, $C=mn(y-m\alpha)(y-n\alpha)$.

By (\ref{formula 2.3}), we have
\begin{equation*}
\aligned
\relax[L_{m,i},L_{n,j}]f(x,y)
=&(m-n)\lambda^{m+n}\zeta^{i+j}
\bigl(x+h_{m+n}(y)+(i+j)\sigma(y-(m+n)\alpha)\bigr)
f(x+m+n,y)\\
&
+(m-n)(m+n)\lambda^{m+n}\zeta^{i+j}
(y-(m+n)\alpha)
\frac{\partial }{\partial y}f(x+m+n,y)\\
=&(m-n)L_{m+n,i+j}f(x,y)+\frac{m^3-m}{12}\delta_{m+n,0}C_{1,i+j}f(x,y).
\endaligned
\end{equation*}

For other similar verification calculations, please refer to Appendix \ref{ Proposition 3.1} for details. 
\end{proof}

\begin{rema}\label{rema3.2}
Let $\boldsymbol{\omega}=(\lambda, \zeta, \alpha,\beta,\gamma,\sigma,h)$ and $\rho_{\boldsymbol{\omega}}:\mathcal G\to\operatorname{End}_{\mathbb C}(\C[x,y])$ be the corresponding representation of $\Omega(\lambda, \zeta, \alpha,\beta,\gamma,\sigma,h)$. Taylor’s formula gives
\begin{equation}\label{Taylor's formula}
\aligned
e^{t\frac{\partial}{\partial{x}}}f(x,y)
=
\sum_{j\ge0}
\frac{t^j}{j!}\Big(\frac{\partial}{\partial{x}}\Big)^jf(x,y)
=f(x+t,y)
\endaligned
\end{equation}
for every $t\in\mathbb C$ and $f\in \C[x,y]$. 
Consequently, the action in Proposition \ref{prop 3.1} admits the equivalent realization
\begin{equation*}\label{act. def 3.1.2}
\aligned
\rho_{\boldsymbol{\omega}}(L_{m,i})
&=\lambda^m\zeta^i
\left[
x+h_m(y)+i\sigma(y-m\alpha)
+m(y-m\alpha)\frac{\partial}{\partial{y}}
\right]e^{m\frac{\partial}{\partial{x}}},
\rho_{\boldsymbol{\omega}}(W_{m,i})
=\lambda^m\zeta^i(y-m\alpha)e^{m\frac{\partial}{\partial{x}}},\\
\rho_{\boldsymbol{\omega}}(I_{r,i})
&=\lambda^{r-\frac12}\zeta^i
\left[
\beta+\gamma g_r(y)+i\sigma\gamma
+r\gamma\frac{\partial}{\partial{y}}
\right]e^{r\frac{\partial}{\partial{x}}},
\rho_{\boldsymbol{\omega}}(B_{r,i})
=\gamma\lambda^{r-\frac12}\zeta^i e^{r\frac{\partial}{\partial{x}}},\
\rho_{\boldsymbol{\omega}}(C_{k,i})=0_{\C[x,y]},
\endaligned
\end{equation*}
where $m,i\in\mathbb Z, r\in\mathbb Z+\frac12, k=1,2,3,4$.
Polynomial coefficients denote multiplication operators, and operator products act from right to left. Expanding the exponentials yields a realization by formal series of differential operators, each of which acts as a finite sum on any given polynomial.


\end{rema}

\begin{rema}\label{rema3.3}
We can restrict the ${\mathcal{G}}$-module \(\Omega(\lambda,\zeta,\alpha,\beta,\gamma,\sigma,h)\) to obtain a \(\hat{\mathcal{W}}\)-module, a \(\mathcal{W}\)-module, or a \(\hat{\mathcal{G}}\)-module.
\end{rema}

\begin{theo}\label{theo 3.3}
    The $\mathcal{G}$-module $\Omega(\lambda,\zeta, \alpha,\beta,\gamma,\sigma,h)$ is simple if and only if $(\alpha, \gamma)\neq(0,0)$.
\end{theo}
\begin{proof} Restricting the module to $\hat{\mathcal W}$, we obtain $\Psi_{\hat{\mathcal{W}}}(\lambda,\alpha,h)$. Thus, 
if $\alpha\ne0$, simplicity follows from Lemma \ref{le2.1}(2).
    For the remainder of the proof, we assume $\alpha=0$. 

\begin{case}
    $\gamma\neq0$.
\end{case}
Let $M$ be a nonzero submodule of $\Omega(\lambda,\zeta, 0,\beta,\gamma,\sigma,h)$. Choose a nonzero element $f(x,y)=\sum_{i=0}^kx^if_i(y)\in M$ of minimal degree in $x$, where $f_{k}(y)\ne0$. Suppose that ${{\rm deg}_x}f(x,y)>0$. Then $\sum_{i=0}^k(x+r)^if_i(y)=\gamma^{-1}\lambda^{\frac{1}{2}-r}B_{r,0}\Big{(}\sum_{i=0}^kx^if_i(y)\Big{)}\in M$
for any $r\in\Z+\frac{1}{2}$. By Lemma \ref{le2.3}, we obtain $f_k(y)\in M$, which contradicts the minimality of ${\rm deg}_x(f(x,y))$. Hence, $f(x,y)=f(y)$. The actions of $L_{0}$ and $W_{0}$ imply that $f(y)\mathbb C[x,y]\subseteq M$. Choose a nonzero $f(y)\in M\cap\mathbb C[y]$ of minimal degree, and suppose that ${\rm deg}_y(f(y))>0$. Clearly, $\bigl(\beta +\gamma g_{r}(y)+i\sigma\gamma\bigl)f(y)+r\gamma\frac{d}{dy}f(y)=\lambda^{\frac{1}{2}-r}\zeta^{-i}I_{r,i}f(y)\in M$
for all $r\in\Z+\frac{1}{2},i\in\Z$, which implies that $\frac{d}{dy}f(y)\in M$. This contradicts the minimality assumption on $f(y)$. Therefore, $M=\C[x,y]$.

\begin{case}
    $\gamma=0$.
\end{case}
Recall that we are still in the standing assumption $\alpha=0$ fixed at the beginning of the proof; indeed the following claim fails when $\alpha\ne0$.
Then $y^i\C[x,y]$ is a submodule of $\Omega(\lambda,\zeta, 0,\beta,0,\sigma,h)$ for every $i\in\N$. In particular, $y\mathbb C[x,y]$ is a nonzero proper submodule.
\end{proof}


\begin{rema}
    For every $i\in\N$, consider the quotient module
 \[
     \Phi_i=y^i\Omega(\lambda,\zeta, 0,\beta,0,\sigma,h)/ y^{i+1}\Omega(\lambda,\zeta, 0,\beta,0,\sigma,h).
 \] Supposing $0\neq y^{i}f(x)\in \Phi_i$, 
 we have $C_{k,j}(y^{i}f(x))=0=B_{r,j}(y^{i}f(x))$,
\begin{equation*}
\aligned
L_{m,j}(y^{i}f(x)) &= \lambda^{m}\zeta^j(x+mh(y)+j\sigma y)y^{i}f(x+m)+\lambda^{m}\zeta^jmiy^{i}f(x+m) \\
&\equiv \lambda^{m}\zeta^j(x+m(h(0)+i))y^{i}f(x+m),\  {\rm mod} 
\ y^{i+1}\mathbb{C}[x, y],
\endaligned
\end{equation*}
$W_{m,j}(y^{i}f(x)) \equiv 0,\  {\rm mod} 
\ y^{i+1}\mathbb{C}[x, y]$ and $I_{r,j}(y^{i}f(x)) \equiv \lambda^{r-\frac{1}{2}}\zeta^j\beta y^{i}f(x+r),\  {\rm mod} 
\ y^{i+1}\mathbb{C}[x, y]$
for all $m,j\in\Z, r\in\Z+\frac{1}{2}$ and $k\in\{1,2,3,4\}$. 
Accordingly, let $\Psi_{\mathcal{H}}(\lambda,\zeta,\alpha,\beta)$ denote $\C[x]$ endowed with the following ${\mathcal{H}}$-action:
\begin{equation*}
\aligned
L_{m,i}f(x)&=\lambda^m\zeta^i(x+m\alpha)f(x+m),\
I_{r,i}f(x)=\beta\lambda^{r-\frac{1}{2}}\zeta^if(x+r),\
C_{1,i} f(x)= C_{3,i} f(x) =0
\endaligned
\end{equation*}
for all $m,i\in\Z,r\in\Z+\frac{1}{2}$, where $\lambda,\zeta\in\C^*,\alpha,\beta\in\C$. Then, we have $\operatorname{Res}_{\hat{\mathcal{H}}}^{\mathcal{H}}\Psi_{\mathcal{H}}(\lambda,\zeta,\alpha,\beta)\cong \Psi_{\hat{\mathcal{H}}}(\lambda,\alpha,\beta)$ and
 $\operatorname{Res}_{\mathcal{H}}^{\mathcal{G}}\Phi_i
     \cong \Psi_{\mathcal{H}}(\lambda,\zeta,h(0)+i,\beta)$.
By Lemma \ref{le2.1.1}, $\Phi_i$ is irreducible if and only if $(h(0)+i,\beta)\neq (0,0)$.
\end{rema}

\begin{theo}\label{theo3.2}
   Let $\lambda, \lambda^\prime,\zeta,\zeta^\prime\in\mathbb{C}^*, \alpha,\alpha^\prime,\beta,\beta^\prime,\gamma,\gamma^\prime,\sigma,\sigma'\in\C,h(y),h^\prime(y)\in\C[y]$. We conclude that $\Omega(\lambda, \zeta, \alpha,\beta,\gamma,\sigma,h)\cong \Omega(\lambda^\prime, \zeta^\prime,\alpha^\prime,\beta^\prime,\gamma^\prime,\sigma',h^\prime)$ if and only if 
   \[
       (\lambda, \zeta, \alpha,\beta,\gamma,\sigma,h)=(\lambda^\prime, \zeta^\prime,\alpha^\prime,\beta^\prime,\gamma^\prime,\sigma',h^\prime).
   \]
\end{theo}
\begin{proof}
    Let $\phi:\Omega(\lambda,\zeta, \alpha,\beta,\gamma,\sigma,h)\to \Omega(\lambda^\prime, \zeta^\prime,\alpha^\prime,\beta^\prime,\gamma^\prime,\sigma',h^\prime)$ be an isomorphism of ${\mathcal{G}}$-modules. For all $i,j\in\N$, the  identity $\phi(L_{0}^{i}W_{0}^j1)=L_{0}^{i}W_{0}^j\phi(1)$ implies that $\phi(x^{i}y^j)=x^{i}y^j \phi(1)$. Since $\phi$ is an isomorphism, 
    $\phi(1)\in\C^*$ 
    and $\phi(f(x,y))=\phi(1)f(x,y)$. Then, for any $m,i\in\mathbb{Z},r\in\Z+\frac{1}{2}$ and $f(x,y)\in\mathbb{C}[x,y]$, the identities 
    $\phi(X1)=X\phi(1),X\in\{L_{m,i}, W_{m,i},I_{r,i},B_{r,i} \}$ imply,
    respectively, that
\begin{equation*}
\aligned
&\phi(1)\lambda^m\zeta^i\left[x+h_{m}(y)+i\sigma(y-m\alpha)\right]=\phi(1)(\lambda^\prime)^m(\zeta^\prime)^i\left[x+h^\prime_{m}(y)+i\sigma'(y-m\alpha')\right],
\\&\phi(1)\lambda^{r-\frac{1}{2}}\zeta^i\left[(\beta+\gamma g_r(y)+i\sigma \gamma)\right] 
= \phi(1)(\lambda^\prime)^{r-\frac{1}{2}}(\zeta^\prime)^i\left[(\beta^\prime+\gamma^\prime g^\prime_r(y)+i\sigma' \gamma')\right],
\\&\phi(1)\lambda^m\zeta^i(y-m\alpha)
=\phi(1)(\lambda^\prime)^m(\zeta^\prime)^i(y-m\alpha^\prime),\ \  \phi(1)\gamma\lambda^{r-\frac{1}{2}}\zeta^i
=\phi(1)\gamma^\prime(\lambda^\prime)^{r-\frac{1}{2}}(\zeta^\prime)^i.
\endaligned
\end{equation*}
Comparison of the two module actions yields $(\lambda, \zeta, \alpha,\beta,\gamma,\sigma,h)=(\lambda^\prime, \zeta^\prime,\alpha^\prime,\beta^\prime,\gamma^\prime,\sigma',h^\prime)$. The converse follows immediately.\end{proof}

\subsection{\texorpdfstring{Classification of $U(\mathcal{h})$-free modules of rank one}{Classification of U(CL0 CW0)-free modules of rank 1}}
Denote $\mathcal{h}=\C L_0\oplus \C W_0
$.
Let $V$ be a $\mathcal{G}$-module which is a free $U(\mathcal{h})$-module of rank one. From the definition of $\mathcal{G}$, we have $L_{0}W_{0}=W_{0}L_{0}$. Let $\mathbf{1}$ be the $\C[L_{0},W_{0}]$-basis of $V$ and assume that $L_{0}\mathbf{1}=x\mathbf{1}$ and $W_{0}\mathbf{1}=y\mathbf{1}$. Then $V=\C[L_{0},W_{0}]\mathbf{1}=\C[x,y]\mathbf{1}$.

For any $m,i\in\Z, r\in\Z+\frac{1}{2}, k\in \{1,2,3,4\}$, let polynomials $l_{m,i},w_{m,i},i_{r,i},b_{r,i}, c_{k,i}\in \C[x,y]$ be such that $L_{m,i}\mathbf{1}=l_{m,i}(x,y)\mathbf{1},W_{m,i}\mathbf{1}=w_{m,i}(x,y)\mathbf{1},I_{r,i}\mathbf{1}=i_{r,i}(x,y)\mathbf{1},
B_{r,i}\mathbf{1}=b_{r,i}(x,y)\mathbf{1},
\linebreak C_{k,i}\mathbf{1}=c_{k,i}(x,y)\mathbf{1}.$

By Lemma \ref{le2.1}, we obtain $l_{m,0}(x,y)=\lambda^m(x+h_{m}(y)),\ w_{m,0}(x,y)=\lambda^m(y-m\alpha),\ c_{1,0}(x,y)=c_{2,0}(x,y)=0.$

The following is an immediate consequence of the relations of $\mathcal{G}$, by induction on $p$.

\begin{lemm}\label{le4.1}
    For all $m,i\in\mathbb{Z}$, $r\in\mathbb{Z}+\frac{1}{2}$, $p\in\mathbb{N}$, and $k\in\{1,2,3,4\}$, the following formulas hold:
    \begin{align*}
        L_{m,i}L_{0}^{p} &= (L_{0}+m)^{p}L_{m,i}, &
        L_{m,i}W_{0}^{p} &= W_{0}^{p}L_{m,i}+pmW_{0}^{p-1}W_{m,i},\\
        W_{m,i}L_{0}^{p} &= (L_{0}+m)^{p}W_{m,i}, &
        W_{m,i}W_{0}^{p} &= W_{0}^{p}W_{m,i},\\
        I_{r,i}L_{0}^{p} &= (L_{0}+r)^{p}I_{r,i}, &
        I_{r,i}W_{0}^{p} &= W_{0}^{p}I_{r,i}+prW_{0}^{p-1}B_{r,i},\\
        B_{r,i}L_{0}^{p} &= (L_{0}+r)^{p}B_{r,i}, &
        B_{r,i}W_{0}^{p} &= W_{0}^{p}B_{r,i},\\
        C_{k,i}L_{0}^{p} &= L_{0}^{p}C_{k,i}, &
        C_{k,i}W_{0}^{p} &= W_{0}^{p}C_{k,i}.
    \end{align*}
\end{lemm}

Using Lemma \ref{le4.1}, we easily obtain the following formulas.
\begin{lemm}\label{Lemma 4.2}
    For any $m,i\in \Z, r\in \Z+\frac{1}{2}, k\in \{1,2,3,4\}$, $f(x,y)\in \C[x,y]$, we have
\begin{equation*}\label{formula4.2.1}
\aligned
L_{m,i} f(x,y)\mathbf{1}=&f(x+m,y)l_{m,i}(x,y)\mathbf{1}+m\frac{\partial {f(x+m,y)}}{\partial {y}}w_{m,i}(x,y)\mathbf{1},\\
W_{m,i} f(x,y)\mathbf{1}=&f(x+m,y)w_{m,i}(x,y)\mathbf{1},\\
I_{r,i} f(x,y)\mathbf{1}=&f(x+r,y)i_{r,i}(x,y)\mathbf{1}+r\frac{\partial {f(x+r,y)}}{\partial {y}}b_{r,i}(x,y)\mathbf{1},\\
B_{r,i} f(x,y)\mathbf{1}=&f(x+r,y)b_{r,i}(x,y)\mathbf{1},\\
C_{k,i} f(x,y)\mathbf{1}=&f(x,y)c_{k,i}(x,y)\mathbf{1}.
\endaligned
\end{equation*}

\end{lemm}

\begin{lemm}\label{le4.2}
For all $m,i\in\mathbb{Z}$, $r\in\Z+\frac{1}{2}$, and $k\in\{1,2,3,4\}$, the   polynomials $w_{m,i}(x,y)=w_{m,i}(y)$, $b_{r,i}(x,y)=b_{r,i}(y)$, and $c_{k,i}(x,y)=c_{k,i}(y)$ belong to $\C[y]$.
\end{lemm}
\begin{proof}
Suppose $w_{m,i}(x,y)=\sum_{p=0}^{r_{m,i}}q_{m,i,p}(y)x^p$, where $r_{m,i}\in\N$, $q_{m,i,p}(y)\in\C[y]$ and $q_{m,i,r_{m,i}}(y)\neq0$.
It follows that
\begin{equation*}
\begin{aligned}
0 &=[W_{m,i}, W_{n,j}]\mathbf{1} 
  = W_{m,i}W_{n,j}\mathbf{1}-W_{n,j}W_{m,i}\mathbf{1} \\
  &= \left(\left(\sum_{p=0}^{r_{n,j}} q_{n,j,p}(y)(x+m)^p\right)\left(\sum_{p=0}^{r_{m,i}} q_{m,i,p}(y)x^p\right)\right.  - \left.\left(\sum_{p=0}^{r_{m,i}} q_{m,i,p}(y)(x+n)^p\right)\left(\sum_{p=0}^{r_{n,j}} q_{n,j,p}(y)x^p\right)\right)\mathbf{1} \\
  &= \left(q_{m,i,r_{m,i}}q_{n,j,r_{n,j}}(m r_{n,j}-n r_{m,i})x^{r_{m,i}+r_{n,j}-1} + \text{lower-degree terms in } x\right)\mathbf{1}.
\end{aligned}
\end{equation*}
Hence $r_{m,i}=0$ for all $m,i\in\mathbb{Z}$. Thus, $w_{m,i}(x,y)=w_{m,i}(y)\in\mathbb{C}[y]$. Similarly, $b_{r,i}(x,y)=b_{r,i}(y)\in\mathbb{C}[y]$ for any $r\in\Z+\frac{1}{2},i\in\Z$. Furthermore, $c_{k,i}(x,y)=c_{k,i}(y)\in\C[y]$ for any $k\in \{1,2,3,4\}$ follows from
$[W_{m,i}, C_{k,i}]\mathbf{1}=c_{k,i}(x+m,y)w_{m,i}(y)\mathbf{1}-w_{m,i}(y)c_{k,i}(x,y)\mathbf{1}=0.$

Suppose $w_{m,i}(y)\mathbf{1}=W_{m,i}\mathbf{1}=0$ for some $m,i\in\Z$. Clearly, for any $n\in\Z\backslash\{0, 2m\},j\in\Z$, we have $W_{n,i+j}\mathbf{1}=
    (n-2m)^{-1}[L_{n-m,j},W_{m,i}]\mathbf{1}=0.$
Then we obtain that $2mW_{0}\mathbf{1}+\frac{m^3-m}{12}C_{2}\mathbf{1}=0$
for any $m\in\Z^*$, which implies  that $W_{0}\mathbf{1}=0$. 
This yields a contradiction.
\end{proof}

\begin{lemm}\label{le4.4}
For all $m,i\in\mathbb{Z}$, we have the following statements.
\begin{itemize}
\item[{\rm (1)}] 
$l_{m,i}(x,y)=\lambda^m\zeta^ix + l_{m,i}^{(0)}(y)$  for some $l_{m,i}^{(0)}(y)\in \mathbb{C}[y]$.
\item[{\rm (2)}] $w_{m,i}(y)=\lambda^m\zeta^i(y-m\alpha)\in\C[y]$ for some $\lambda,\zeta\in \C^*, \alpha\in\C$.
\item[{\rm (3)}] $c_{2,i}(y)=0$.

\end{itemize}
\end{lemm}
\begin{proof}
    According to the relation  $[L_{m,i},W_{n,j}]\mathbf{1}=(m-n)W_{m+n,i+j}\mathbf{1}+\frac{m^3-m}{12}\delta_{m+n,0}C_{2,i+j}\mathbf{1},$
 we have
\begin{equation}\label{formula 3.2}
\aligned
&(m-n)w_{m+n,i+j}(y)\mathbf{1}+\frac{m^3-m}{12}\delta_{m+n,0}c_{2,i+j}(y)\mathbf{1}\\
=&\big(l_{m,i}(x,y)-l_{m,i}(x+n,y)\big)w_{n,j}(y)\mathbf{1}+m\frac{d}{d{y}}w_{n,j}(y)w_{m,i}(y)\mathbf{1}
\endaligned
\end{equation}
for any $m,n,i,j\in\Z$. So, we can deduce ${\rm deg}_x(l_{m,i}(x,y))\leq1$ for all $m,i\in\Z$. Then, for any $m,i\in\Z$, let us set $l_{m,i}(x,y)=l_{m,i}^{(1)}(y)x+l_{m,i}^{(0)}(y)$, where $l_{m,i}^{(1)}(y),l_{m,i}^{(0)}(y)\in\C[y]$ and $l_{0,0}^{(1)}(y)=1,l_{0,0}^{(0)}(y)=0$. Hence, we have
\begin{equation}\label{formula 3.3}
\aligned
 (m-n)w_{m+n,i+j}(y)+\frac{m^3-m}{12}\delta_{m+n,0}c_{2,i+j}(y)
=-nw_{n,j}(y)l_{m,i}^{(1)}(y)+m\frac{d}{d{y}}{w_{n,j}(y)}w_{m,i}(y)
\endaligned
\end{equation}
for any $m,n,i,j\in\Z$. Letting $m=n\neq 0$ and $i=j$ in (\ref{formula 3.3}), we obtain $l_{m,i}^{(1)}(y)=\frac{d}{dy}{w_{m,i}(y)}$. 

Taking $m=\pm1$ in (\ref{formula 3.3}) yields
\begin{eqnarray}
    &&nw_{n,j}(y)\frac{d}{d{y}}{w_{1,i}(y)}-\frac{d}{d{y}}{w_{n,j}(y)}w_{1,i}(y)=(n-1)w_{n+1,i+j}(y),\label{formula 3.4} \\
    &&nw_{n,j}(y)\frac{d}{d{y}}{w_{-1,i}(y)}+\frac{d}{d{y}}{w_{n,j}(y)}w_{-1,i}(y)=(n+1)w_{n-1,i+j}(y) \label{formula 3.5}
\end{eqnarray}
for all $n\in\Z$. Setting $i=0$ and substituting $n= 1$ and $n=-1$ into equations (\ref{formula 3.4}) and (\ref{formula 3.5}) respectively, we have $(y-\alpha)\frac{d}{dy}w_{1,j}(y)-w_{1,j}(y)=0,\ (y+\alpha)\frac{d}{dy}w_{-1,j}(y)-w_{-1,j}(y)=0.$
It follows that $w_{1,j}(y)=\lambda\zeta_j(y-\alpha)$ and $w_{-1,j}(y)=\lambda^{-1}u_{j}(y+\alpha)$, where $\zeta_j,u_{j}\in\C^*$. Taking $n=-1$ in (\ref{formula 3.4}), we obtain that $w_{0,i+j}(y)=\zeta_i u_jy.$
  It is clear that  $w_{0,i}(y)=\zeta_i y=u_i y$, $\zeta_0=u_0=1$ and $u_j=\zeta_{-j}^{-1}\in\C^*$. Taking $n=0 $ in (\ref{formula 3.4}), we have $\zeta_j\zeta_i=\zeta_{i+j}$ and $w_{0,i}(y)=\zeta_1^iy:=\zeta^iy$ for some $\zeta\in\C^*$. 
  Setting $n=i=0$ in (\ref{formula 3.2}), we can obtain $mw_{m,j}(y)+\frac{m^3-m}{12}\delta_{m,0}c_{2,j}(y)
=m\frac{d}{d{y}}w_{0,j}(y)w_{m,0}(y)$
for any $m,j\in\Z$,  which means that $w_{m,i}(y)=\frac{d}{dy}{w_{0,i}(y)}w_{m,0}(y)$
   for all $m,i\in\Z$.
  Hence, we can verify that $w_{m,i}(y)=\lambda^m\zeta^i(y-m\alpha)$ .
  Letting $m=j=0$ and $n\neq 0$ in (\ref{formula 3.2}), we can obtain
  $-n\lambda^n\zeta^i(y-n\alpha)=-nl_{0,i}^{(1)}(y)\lambda^n(y-n\alpha).$
Therefore, $l_{0,i}^{(1)}(y)=\zeta^i=\frac{d}{dy}w_{0,i}(y)$ and $l_{m,i}^{(1)}(y)=\frac{d}{dy}w_{m,i}(y)=\lambda^m\zeta^i$ for all $m,i\in\Z$. Taking $n=-m$ and $j=0$ in (\ref{formula 3.2}), we have $2m\zeta^iy+\frac{m^3-m}{12}c_{2,i}(y)=m\lambda^m\zeta^i\lambda^{-m}(y+m\alpha)+m\lambda^{-m}\lambda^{m}\zeta^i(y-m\alpha)
=2m\zeta^iy,$
which implies that $c_{2,i}(y)=0$.
\end{proof}

\begin{lemm}\label{le4.5}
For all $m,i\in\mathbb{Z}$, we have $l_{m,i}(x,y)=\lambda^m\zeta^i\big(x+h_m(y)+i\sigma(y-m\alpha)\big), c_{1,i}(y)=0$ for some $\sigma\in\C$.
\end{lemm}
\begin{proof}
  By $(m-n)L_{m+n,i+j}+\frac{m^3-m}{12}\delta_{m+n,0}C_{1,i+j} =[L_{m,i},L_{n,j}]$, we have
\begin{equation}\label{formula4.5.1}
\aligned
&(m-n)\lambda^{m+n}\zeta^{i+j}x+(m-n)l_{m+n,i+j}^{(0)}(y)+\frac{m^3-m}{12}\delta_{m+n,0}c_{1,i+j}(y)-(m-n)\lambda^{m+n}\zeta^{i+j}x
\\
=&m\lambda^n\zeta^{j}l_{m,i}^{(0)}(y)-n\lambda^m\zeta^{i}l_{n,j}^{(0)}(y)+m\lambda^m\zeta^{i}\frac{d}{d{y}}{l_{n,j}^{(0)}(y)}(y-m\alpha)-n\lambda^n\zeta^{j}\frac{d}{d{y}}{l_{m,i}^{(0)}(y)}(y-n\alpha)
\endaligned
\end{equation}
for all $m,n, i,j\in\Z$. Taking $m=j=0$ and $n\ne 0$ in (\ref{formula4.5.1}), we have
\begin{equation}\label{formula4.5.2}
\aligned
l_{n,i}^{(0)}(y)=\lambda^n\zeta^ih_n(y)+\lambda^n(y-n\alpha)   \frac{d}{d{y}}{l_{0,i}^{(0)}(y)}. 
\endaligned
\end{equation}
Letting $m=n=1,j=0$ in (\ref{formula4.5.1}) and combining with (\ref{formula4.5.2}), we obtain
\begin{equation*}\label{formula4.5.3}
\aligned
0=
  l_{1,i}^{(0)}(y)-\lambda \zeta^ih(y)
+\lambda \zeta^i(y-\alpha)\frac{d}{d{y}} h(y)- (y-\alpha)
\frac{d }{d y}l_{1,i}^{(0)}(y),
\endaligned
\end{equation*}
and hence $(y-\alpha)^2\frac{d^2}{d{y^2}}l_{0,i}^{(0)}(y)   =0$.
Since $l_{0,i}^{(0)}(y)\in\C[y]$, then
\begin{equation}\label{formula4.5.4}
\aligned
\frac{d^2}{d{y^2}}l_{0,i}^{(0)}(y)   =0.
\endaligned
\end{equation}
Taking $m=-n=1,j=0$ in (\ref{formula4.5.1}) and combining with (\ref{formula4.5.2}), we get
\begin{equation*}\label{formula4.5.5}
\aligned
2l_{0,i}^{(0)}(y)=\lambda^{-1}l_{1,i}^{(0)}(y)+\zeta^ih_{-1}(y)+\zeta^i(y-\alpha)\frac{d}{d{y}} h_{-1}(y)+\lambda^{-1}(y+\alpha)\frac{d}{d{y}} l_{1,i}^{(0)}(y),
\endaligned
\end{equation*}
and hence
\begin{equation}\label{formula4.5.6}
\aligned
(y^2-\alpha^2)\frac{d^2}{d{y^2}}l_{0,i}^{(0)}(y)+2y\frac{d}{d{y}}l_{0,i}^{(0)}(y)-2l_{0,i}^{(0)}(y)=0.
\endaligned
\end{equation}
Through (\ref{formula4.5.4}), we obtain $y\frac{d}{d{y}}l_{0,i}^{(0)}(y)=l_{0,i}^{(0)}(y)$, which means that $l_{0,i}^{(0)}(y)$ is either a proportional function of $y$ or a zero function, i.e. 
\begin{equation}\label{formula4.5.7}
\aligned
l_{0,i}^{(0)}(y)=k_iy
\endaligned
\end{equation}
 for some $k_i\in \C$. Letting $m=1,n=0$ in (\ref{formula4.5.1}), combining (\ref{formula4.5.2}) and (\ref{formula4.5.7}) ,  we have $l_{1,i+j}^{(0)}(y)=\zeta^jl_{1,i}^{(0)}(y)+\lambda\zeta^i(y-\alpha)\frac{d}{d{y}}l_{0,j}^{(0)}(y),$
and hence $k_{i+j}=\zeta^jk_i+\zeta^ik_j,\ {\rm\ i.e.}\ \frac{k_{i+j}}{\zeta^{i+j}}=\frac{k_i}{\zeta^i}+\frac{k_j}{\zeta^j}.$
Then we can verify $k_i=i\sigma \zeta^i$ for some $\sigma=\frac{k_1}{\zeta}\in \C$. Therefore, by (\ref{formula4.5.2}), we have $l_{m,i}(x,y)=\lambda^m\zeta^i\big(x+h_m(y)+i\sigma(y-m\alpha)\big), c_{1,i}(y)=0$. Indeed, $c_{1,i}(y)=0$ can be seen directly by applying (\ref{formula4.5.1}) with $n=-m$ and $j=0$, which gives $[L_{m,i},L_{-m,0}]\mathbf{1}=2mL_{0,i}\mathbf{1}+\frac{m^3-m}{12}C_{1,i}\mathbf{1}$; substituting (\ref{formula4.5.2}) and $l_{0,i}^{(0)}(y)=k_iy=i\sigma\zeta^iy$ shows that all the $L$-terms cancel, whence $\frac{m^3-m}{12}c_{1,i}(y)=0$ and therefore $c_{1,i}(y)=0$ (as $\frac{m^3-m}{12}\ne0$ for $m=2$).
\end{proof}

\begin{lemm}\label{le4.6}
For all $i\in\Z,r\in\mathbb{Z}+\frac{1}{2}$, we have $b_{r,i}(y)=\gamma\lambda^{r-\frac{1}{2}}\zeta^i$ for some $\gamma\in \C$.
\end{lemm}
\begin{proof}
Since
\begin{equation}\label{formula 3.6}
\aligned
-rb_{m+r,i+j}(y)=&-rB_{m+r,i+j} \mathbf{1}
=[L_{m,i},B_{r,j}] \mathbf{1}
=m\lambda^m\zeta^i(y-m\alpha)\frac{d }{d {y}}b_{r,j}(y)-r\lambda^m\zeta^ib_{r,j}(y),
\endaligned
\end{equation}
by fixing $i=j=0$, and letting $m=-1, r=\frac{1}{2}$ and $m=1, r=-\frac{1}{2}$ in (\ref{formula 3.6}) respectively, we can obtain that
\begin{equation*}\label{formula 4.6.2}
\aligned
&b_{-\frac{1}{2},0}(y)=2\lambda^{-1}(y+\alpha)\frac{d}{dy}b_{\frac{1}{2},0}(y)+\lambda^{-1}b_{\frac{1}{2},0}(y),\
b_{\frac{1}{2},0}(y)=2\lambda(y-\alpha)\frac{d}{dy}b_{-\frac{1}{2},0}(y)+\lambda b_{-\frac{1}{2},0}(y).
\endaligned
\end{equation*}
By substituting the elimination method, we have the following second-order homogeneous linear differential equation
\begin{equation}\label{formula 3.7}
\aligned
(y^2-\alpha^2)\frac{d^2}{dy^2}b_{\frac{1}{2},0}(y)+(2y-\alpha)\frac{d}{dy}b_{\frac{1}{2},0}(y)=0.
\endaligned
\end{equation}
Letting $b(y):=\frac{d}{dy}b_{\frac{1}{2},0}(y)$, 
we calculate that 
    $b(y)
    =\frac{S_1}{y^2-\alpha^2}\sqrt{\frac{y-\alpha}{y+\alpha}},$
where $S_1\in\C$.
If $\alpha\neq 0$ and letting $\theta:=\sqrt{\frac{y-\alpha}{y+\alpha}}$, then we have 
$\frac{d}{dy}\log{\theta}=\frac{1}{\theta}\frac{d\theta}{dy}
=\frac{\alpha}{y^2-\alpha^2}.$
Hence, 
$\frac{d\theta}{dy}=\frac{\alpha}{y^2-\alpha^2}\sqrt{\frac{y-\alpha}{y+\alpha}}$.
It is not difficult to obtain that $b_{\frac{1}{2},0}(y)=K_0\sqrt{\frac{y-\alpha}{y+\alpha}}+\gamma.$
 If $\alpha=0$, then $b(y)=\frac{S_1}{y^2}$ and $b_{\frac{1}{2},0}(y)=\frac{-S_1}{y}+\gamma.$
 However, since $b_{\frac{1}{2},0}(y)\in\C[y]$, we obtain $b_{\frac{1}{2},0}(y)=\gamma$ and $b_{-\frac{1}{2},0}(y)=\lambda^{-1}\gamma$. Taking $i=j=0$, $0\neq m=\frac{1}{2}-r$ and $0\neq m=-\frac{1}{2}-r$  in (\ref{formula 3.6}), respectively, we obtain this pair of first-order non-homogeneous linear differential equations with mirror symmetry: 
\begin{equation}\label{formula 3.8}
\aligned
\left(r-1/2\right)\left[y+\left(r-1/2\right)\alpha\right]\frac{d}{dy}b_{r,0}(y)+rb_{r,0}(y)=r\gamma\lambda^{r-\frac{1}{2}},
\endaligned
\end{equation}
\begin{equation}\label{formula 3.8.1}
\aligned
\left(r+1/2\right)\left[y+\left(r+1/2\right)\alpha\right]\frac{d}{dy}b_{r,0}(y)+rb_{r,0}(y)=r\gamma\lambda^{r-\frac{1}{2}}.
\endaligned
\end{equation}
By subtracting (\ref{formula 3.8}) from (\ref{formula 3.8.1}), we get a first-order homogeneous linear differential equation $(y+2r\alpha)\frac{d}{dy}b_{r,0}(y)=0.$
 It follows from  $b_{r,0}(y)\in\mathbb{C}[y]$ that we have $\frac{d}{dy}b_{r,0}(y)=0$  and
 $b_{r,0}(y)=\gamma\lambda^{r-\frac{1}{2}}$ for any $r\in\Z+\frac{1}{2}$. 
 Due to
\begin{equation*}
\aligned
-rb_{r,i}(y)\mathbf{1}
&=[L_{0,i},B_{r,0}]\mathbf{1}=\zeta^i
\bigl(x+i\sigma y\bigr)b_{r,0}(y)\mathbf{1}
-\zeta^i
\bigl(x+r+i\sigma y\bigr)b_{r,0}(y)\mathbf{1}=-r\zeta^ib_{r,0}(y)\mathbf{1},
\endaligned
\end{equation*}
we can verify $b_{r,i}(y)=\gamma\lambda^{r-\frac{1}{2}}\zeta^i$.
\end{proof}

\begin{lemm}\label{le4.7}
For all $i\in\Z, r\in\mathbb{Z}+\frac{1}{2},k\in\{3,4\}$, we have $i_{r,i}(x,y)=\lambda^{r-\frac{1}{2}}\zeta^i(\beta+\gamma g_r(y)+i\sigma \gamma)$ and $c_{k,i}(y)=0$, where $\beta\in\C$. 
\end{lemm}
\begin{proof}
According to $[W_{m,i}, I_{r,j}]=-rB_{m+r,i+j}$, we have $\lambda^m\zeta^i(y-m\alpha)(i_{r,j}(x+m,y)-i_{r,j}(x,y))=0$
for any $m,i\in\mathbb{Z},r\in\Z+\frac{1}{2}$. It is easy to obtain $i_{r,i}(x,y)=i_{r,i}(y)\in\mathbb{C}[y]$ for $r\in\mathbb{Z}+\frac{1}{2},i\in\Z$. Thus, we can deduce that $c_{4,i}(y)=0$ from $[I_{r,i},B_{-r,0}]=rC_{4,i}$ for any $r\in\Z+\frac{1}{2},i\in\Z$.


Then, by $-rI_{m+r,i+j}\mathbf{1}=[L_{m,i},I_{r,j}]\mathbf{1}$, we have
\begin{equation}\label{formula 3.9}
\aligned
-ri_{m+r,i+j}(y) 
=m\lambda^m\zeta^i(y-m\alpha)\frac{d }{d {y}}i_{r,j}(y) -r\lambda^m\zeta^ii_{r,j}(y) -r\gamma\lambda^{m+r-\frac{1}{2}}\zeta^{i+j}\Bigg(\frac{d }{d {y}}h_{m}(y)+i\sigma\Bigg).
\endaligned
\end{equation}
Fixing $i=j=0$, and letting $m=1,r=-\frac{1}{2}$ and $m=-1,r=\frac{1}{2}$ in (\ref{formula 3.9}) respectively, we have 
\begin{eqnarray}
\label{formula33.18}
&&i_{\frac{1}{2},0}(y)=2\lambda(y-\alpha)\frac{d}{dy}i_{-\frac{1}{2},0}(y)+\lambda i_{-\frac{1}{2},0}(y)+\gamma\frac{d}{dy}h_1(y),\\
&& \label{formula33.19} i_{-\frac{1}{2},0}(y)=2\lambda^{-1}(y+\alpha)\frac{d}{dy}i_{\frac{1}{2},0}(y)+\lambda^{-1} i_{\frac{1}{2},0}(y)+\gamma\lambda^{-1}\frac{d}{dy}h_{-1}(y).
\end{eqnarray}
Substituting \eqref{formula33.19} into \eqref{formula33.18}, we obtain
\begin{equation*}\label{formula 4.7.1.1}
\aligned
(y^2-\alpha^2)\frac{d^2}{dy^2}i_{\frac{1}{2},0}(y)+(2y-\alpha)\frac{d}{dy} i_{\frac{1}{2},0}(y)=-\frac{\gamma}{4}\left(2(y-\alpha)\frac{d^2}{dy^2}h_{-1}(y)+\frac{d}{dy}h_{-1}(y)+\frac{d}{dy}h_1(y)\right).
\endaligned
\end{equation*}
By integrating once, we deduce 
\begin{equation*}\label{formula 4.7.1.2}
\aligned
(y^2-\alpha^2)\frac{d}{dy}i_{\frac{1}{2},0}(y)-\alpha i_{\frac{1}{2},0}(y)
&=\frac{\gamma}{4}\left(2(y+\alpha)\frac{d}{dy}h(y)-6\alpha\frac{h(y)-h(\alpha)}{y-\alpha}-2h(y)\right)+K_2,
\endaligned
\end{equation*}
where $K_2\in\C$. Letting $g(y):=g_1(y)=\frac{h(y)-h(\alpha)}{y-\alpha}$, we obtain $h(y)=(y-\alpha)g(y)+h(\alpha),\ \frac{d}{dy}h(y)=(y-\alpha)\frac{d}{dy}g(y)+g(y), g_r(y)=rg(y)$
and
\begin{equation}\label{formula 3.10}
\aligned
(y^2-\alpha^2)\frac{d}{dy}i_{\frac{1}{2},0}(y)-\alpha i_{\frac{1}{2},0}(y)=\frac{\gamma}{2}\left((y^2-\alpha^2)\frac{d}{dy}g(y)-\alpha g(y)-h(\alpha)\right)+K_2.
\endaligned
\end{equation}
By taking the derivative of the two sides of equation (\ref{formula 3.10}) with respect to $y$, we have
\begin{equation}\label{formula 3.11}
\aligned
(y^2-\alpha^2)\frac{d^2}{dy^2}i_{\frac{1}{2},0}(y)+(2y-\alpha) \frac{d}{dy}i_{\frac{1}{2},0}(y)=\frac{\gamma}{2}\left((y^2-\alpha^2)\frac{d^2}{dy^2}g(y)+(2y-\alpha) \frac{d}{dy}g(y)\right).
\endaligned
\end{equation}
Letting $\hat{i}_{\frac{1}{2},0}(y):=i_{\frac{1}{2},0}(y)-\frac{\gamma}{2}g(y)$, then equation (\ref{formula 3.11}) can be transformed into the same homogeneous form as equation (\ref{formula 3.7}).
Since $i_{\frac{1}{2},0}(y)\in\C[y]$, we can easily obtain that $i_{\frac{1}{2},0}(y)=\beta+\frac{1}{2}\gamma g(y),
    i_{-\frac{1}{2},0}(y)=\beta\lambda^{-1}-\frac{1}{2}\gamma\lambda^{-1}g(y)$,
where $\beta\in\C$. Taking $i=j=0,r=\frac{1}{2}$ and $m=s-\frac{1}{2}$ for any $s\in\Z+\frac{1}{2}$ in (\ref{formula 3.9}), we can obtain that $i_{s,0}(y)=\lambda^{s-\frac{1}{2}}(\beta+\gamma g_s(y)).$
Taking $m=j=0$ in (\ref{formula 3.9}), we can obtain that $i_{r,i}(y)=\zeta^ii_{r,0}(y)+\gamma \lambda^{r-\frac{1}{2}}\zeta^ii\sigma$. 
Obviously, we can determine that $i_{r,i}(y)=\lambda^{r-\frac{1}{2}}\zeta^i(\beta+\gamma g_r(y)+i\sigma\gamma)$ for any $r\in\Z+\frac{1}{2},i\in\Z$.
Therefore, considering $[I_{r,i},I_{-r,0}]\mathbf{1}=rC_{3,i}\mathbf{1}$, we have $c_{3,i}(y)=0$.
\end{proof}

Lemmas \ref{Lemma 4.2}, \ref{le4.4}, \ref{le4.5}, \ref{le4.6}, and \ref{le4.7} immediately yield the main result of this paper.
\begin{theo}\label{th4.1}
 Let $V$ be a $\mathcal{G}$-module that is free of rank one over $U(\mathcal{h})$. Then $V\cong\Omega(\lambda, \zeta, \alpha, \beta,\linebreak\gamma,\sigma,h)$ for some $\lambda,\zeta\in\C^*$, $\alpha, \beta, \gamma,\sigma\in\C$, and $h(y)\in\C[y]$.
\end{theo}

By Remark \ref{rema3.3} and Theorem \ref{th4.1}, we have the following corollary.

\begin{coro}
    Let $V$ be a $\hat{\mathcal{G}}$-module that is free of rank one over $U(\mathcal{h})$. Then, for any $\zeta\in\C^*$ and $\sigma\in\C$, $V\cong {\operatorname{Res}}_{\hat{\mathcal{G}}}^{\mathcal{G}}\Omega(\lambda, \zeta,\alpha,\beta,\gamma,\sigma,h)$ for some $\lambda\in\C^*$, $\alpha, \beta, \gamma\in\C$, and $h(y)\in\C[y]$.
\end{coro}






\section{\texorpdfstring{ Self-extensions of $\mathcal{G}$-modules $\Omega(\lambda, \zeta, \alpha,\beta,\gamma,\sigma,h)$}{1}}\label{section 4}

We first recall the standard
cohomological interpretation of self-extensions; see
\cite{ChevalleyEilenberg,Fuks,Weibel}. 
    Let $\mathcal{L}$ be a Lie algebra over $\C$, and let $V$ be an $\mathcal{L}$-module. A linear map $\tau:\mathcal{L}\to V$ is called a $1$-cocycle 
    if 
    $\tau\big([X,Y]\big)=X\tau(Y)-Y\tau(X)$
    for all $X,Y\in \mathcal{L}$. A $1$-cocycle of the form $\tau_{T}(X):=XT$ for $T\in V$ is called a $1$-coboundary determined by $T$. The quotient of the space $Z^1(\mathcal{L},V)$ of $1$-cocycles by the space $B^1(\mathcal{L},V)$ of $1$-coboundaries is denoted by $H^1(\mathcal{L},V)$.
    Let $M,N$ be $\mathcal{L}$-modules. An extension of $M$ by $N$ is a short exact sequence of $\mathcal{L}$-modules
$0\longrightarrow N\xrightarrow{\iota} E
\xrightarrow{\pi}M\longrightarrow0$, where $\iota,\pi$ are  $\mathcal{L}$-module homomorphisms with $\iota$ injective, $\pi$ surjective and $\operatorname{im} (\iota)=\ker (\pi)$. The extension is called split if there exists an $\mathcal{L}$-module homomorphism $\mathbf{s}:M\to E$ such that $\pi\circ\mathbf{s}=\operatorname{Id}_{M}$, and nonsplit otherwise. When $N=M$, such an extension is called a self-extension of $M$. Two self-extensions $0\longrightarrow M\xrightarrow{\iota} E
\xrightarrow{\pi}M\longrightarrow0$, $0\longrightarrow M\xrightarrow{\iota'} E'
\xrightarrow{\pi'}M\longrightarrow0$ are called equivalent if there exists an $\mathcal{L}$-module isomorphism $\psi:E\to E'$ satisfying $\psi\circ \iota=\iota'$ and $\pi'\circ\psi=\pi$. The equivalence classes of self-extensions are classified by $\operatorname{Ext}^1_{U(\mathcal{L})}(M,M)\cong H^1(\mathcal{L},\operatorname{End}_{\mathbb C}(M))$, with the zero element corresponding to the split self-extension.

Denote $M=\Omega(\lambda,\zeta,\alpha,\beta,\gamma,\sigma,h)$
and $\boldsymbol{\omega}=(\lambda,\zeta,\alpha,\beta,\gamma,\sigma,h)
.$
Let $\rho_{\boldsymbol{\omega}}: \mathcal G\longrightarrow\operatorname{End}_{\mathbb C}(M)$
be the corresponding representation. 
The representation $\rho_{\boldsymbol{\omega}}$
induces a $\mathcal{G}$-module structure on $\operatorname{End}_{\mathbb C}(M)$ by $XT=[\rho_{\boldsymbol{\omega}}(X),T],
X\in \mathcal{G},T\in\operatorname{End}_{\mathbb C}(M).$
This interpretation is also compatible with the classical deformation
theory of Lie algebra homomorphisms; see
\cite{NijenhuisRichardson}. Indeed, let $\Lambda=\mathbb C[\varepsilon]/(\varepsilon^2)$
be the algebra of dual numbers. We use the same symbols $\rho_{\boldsymbol{\omega}}(X),\tau(X)$ and $T$ for their natural $\Lambda$-linear extensions $\rho_{\boldsymbol{\omega}}(X)\otimes \operatorname{Id}_{\Lambda},\tau(X)\otimes\operatorname{Id}_{\Lambda}$ and $T\otimes\operatorname{Id}_{\Lambda}$ to $M\otimes_{\C}\Lambda$, respectively; accordingly, $\varepsilon\tau(X)$ and $\varepsilon T$ stand for $\tau(X)\otimes(\varepsilon\operatorname{Id}_{\Lambda})$ and $T\otimes(\varepsilon\operatorname{Id}_{\Lambda})$.
A first-order deformation of
$\rho_{\boldsymbol{\omega}}$ is a 
$\Lambda$-linear representation that can be written as $\rho_\varepsilon(X)
=
\rho_{\boldsymbol{\omega}}(X)+\varepsilon\tau(X)\in\operatorname{End}_{\Lambda}(M\otimes_{\C}\Lambda).$
Comparing the coefficients of $\varepsilon$ in $\rho_\varepsilon([X,Y])
=
[\rho_\varepsilon(X),\rho_\varepsilon(Y)]$
shows that $\tau$ satisfies 
\begin{equation}\label{formula 4.1}
\aligned
\tau\big([X,Y]\big)=[\rho_{\boldsymbol{\omega}}(X), \tau(Y)]+[\tau(X),\rho_{\boldsymbol{\omega}}(Y)],
\endaligned
\end{equation}
for all $X,Y\in \mathcal{G}$.  Moreover, conjugation by
$\operatorname{Id}_{M\otimes_{\C}\Lambda}-\varepsilon T$ changes $\rho_{\boldsymbol{\omega}}(X)$ into $$(\operatorname{Id}_{M\otimes_{\C}\Lambda}-\varepsilon T)\rho_{\boldsymbol{\omega}}(X)(\operatorname{Id}_{M\otimes_{\C}\Lambda}+\varepsilon T)
=
\rho_{\boldsymbol{\omega}}(X)
+\varepsilon[\rho_{\boldsymbol{\omega}}(X),T],$$
so that $1$-coboundaries $\tau_T(X)=[\rho_{\boldsymbol{\omega}}(X),T]$ correspond to infinitesimally trivial
changes of basis. Thus $H^1\bigl(\mathcal G,\operatorname{End}_{\mathbb C}(M)\bigr)
=
Z^1\bigl(\mathcal G,\operatorname{End}_{\mathbb C}(M)\bigr)/
B^1\bigl(\mathcal G,\operatorname{End}_{\mathbb C}(M)\bigr).$

Since every short
exact sequence of complex vector spaces splits, 
the underlying vector space of any self-extension of $M$ is non-canonically isomorphic to $M\oplus M$. Under such an identification, its $\mathcal G$-action has the form
\begin{equation}\label{formula 4.3}
\aligned
X
\begin{pmatrix}
    f_1\\
    f_2
\end{pmatrix}
=\begin{pmatrix}
    \rho_{\boldsymbol{\omega}}(X)f_1+\tau(X)f_2\\
    \rho_{\boldsymbol{\omega}}(X)f_2 
\end{pmatrix}
=\begin{pmatrix}
    \rho_{\boldsymbol{\omega}}(X) & \tau(X) \\
    0_M & \rho_{\boldsymbol{\omega}}(X)
\end{pmatrix}\begin{pmatrix}
    f_1\\
    f_2
\end{pmatrix}, \ 
\begin{pmatrix}
    f_1\\
    f_2
\end{pmatrix}\in M\oplus M.
\endaligned
\end{equation} 
Letting 
\begin{equation*}
\aligned
\Theta: M\oplus M&\to  M\otimes_{\C} \Lambda,\\
(f_1,f_2)^{\operatorname{T}}&\mapsto f_2+\varepsilon f_1,
\endaligned
\end{equation*}
we have
\begin{equation*}
\aligned
\Theta\left( X
\begin{pmatrix}
    f_1\\
    f_2
\end{pmatrix}  \right)=&\rho_{\boldsymbol{\omega}}(X)f_2 +\varepsilon\left(\rho_{\boldsymbol{\omega}}(X)f_1+\tau(X)f_2 \right)
=(\rho_{\boldsymbol{\omega}}(X)+\varepsilon\tau(X))(f_2+\varepsilon f_1)=\rho_{\boldsymbol{\varepsilon}}(X)\Theta\left( 
\begin{pmatrix}
    f_1\\
    f_2
\end{pmatrix}  \right).
\endaligned
\end{equation*}
 Hence, $\Theta$ is a module isomorphism and $\tau$ determines a self-extension
\begin{equation}\label{formula 4.4}
\aligned
0\longrightarrow M\xrightarrow{\iota} E_\tau\xrightarrow{\pi} M\longrightarrow0,
\endaligned
\end{equation}
where $\iota(f)=(f,0)^{\operatorname{T}},\pi\left((f_1,f_2)^{\operatorname{T}}\right)=f_2$. It is worth noting that
$\iota\circ\tau(X)\circ\pi=\begin{pmatrix}
    0_M & \tau(X)\\
    0_M & 0_M
\end{pmatrix},$
$$\Theta\begin{pmatrix}
    \rho_{\boldsymbol{\omega}}(X) & \tau(X) \\
    0_M & \rho_{\boldsymbol{\omega}}(X)
\end{pmatrix}\Theta^{-1}=\rho_{\boldsymbol{\varepsilon}}(X), \ \Theta\begin{pmatrix}
    0_M & \tau(X) \\
    0_M & 0_M
\end{pmatrix}\Theta^{-1}=\varepsilon\tau(X),\ 
\Theta\begin{pmatrix}
    \rho_{\boldsymbol{\omega}}(X) & 0_M \\
    0_M & \rho_{\boldsymbol{\omega}}(X)
\end{pmatrix}\Theta^{-1}=\rho_{\boldsymbol{\omega}}(X).$$
A different such identification replaces  $\tau$ by
$\tau+\tau_T$ for some $T\in\operatorname{End}_{\mathbb C}(M)$.
Consequently, two cocycles determine equivalent self-extensions
inducing the identity maps on the two copies of $M$ if and only if
their difference is a $1$-coboundary. Therefore, one has the standard
natural isomorphism
\begin{equation}\label{formula 4.5}
\aligned
H^1(\mathcal{G},\operatorname{End}_{\mathbb C}(M))\cong \operatorname{Ext}^1_{U(\mathcal{G})}(M,M).
\endaligned
\end{equation}

We shall first construct an explicit family of $1$-cocycles and then prove that every $1$-cocycle is cohomologous to a member of this family, thereby obtaining a classification of the self-extensions of $M$.
Throughout
this section, denote $h'(y)=\frac{d}{d{y}}h(y)$, $h'(\alpha)=\frac{d}{d{y}}h(y)\big{|}_{y=\alpha}$ for any $h(y)\in \C[y]$ and put $\Delta_tF(x,y)=F(x+t,y)-F(x,y)$
for $F(x,y)\in\mathbb C[x,y]$ and $t\in\mathbb C$. Obviously, $\Delta_t$ is surjective on $\C[x,y]$.

\subsection{\texorpdfstring{$1$-cocycles $\tau_{\boldsymbol{v}}$}{1}}
Let $\boldsymbol{v}=(u,v,a,b,c,s,q)$, where $u,v,a,b,c,s\in \C,q\in\C[y]$. Set $\lambda_{\varepsilon}=\lambda(1+\varepsilon u), \zeta_{\varepsilon}=\zeta(1+\varepsilon v), 
\alpha_{\varepsilon}=\alpha+\varepsilon a,\beta_{\varepsilon}=\beta+\varepsilon b, 
\gamma_{\varepsilon}=\gamma+\varepsilon c,  \sigma_{\varepsilon}=\sigma+\varepsilon s,
h_{\varepsilon}=h+\varepsilon q$.
For a polynomial $F$ in the module parameters, let $F_{\varepsilon}$ be the expression obtained by the above substitutions, and write $F_{\varepsilon}=F+\varepsilon\,\delta F$, where the first-order variation $\delta F:=[\varepsilon]F_{\varepsilon}$ is the coefficient of $\varepsilon$ in $F_{\varepsilon}$.
Then, we have $h_\varepsilon(\alpha_\varepsilon)
=h(\alpha)+\varepsilon\bigl(q(\alpha)+ah'(\alpha)\bigr)$. 
Furthermore, $\frac1{y-\alpha-\varepsilon a}
=\frac1{y-\alpha}
+\varepsilon\frac{a}{(y-\alpha)^2}.$
It follows that 
\begin{equation}\label{B7}
\aligned
(g_r)_\varepsilon(y)
&=
r\frac{h_\varepsilon(y)-h_\varepsilon(\alpha_\varepsilon)}
        {y-\alpha_\varepsilon}
=g_r(y)+\varepsilon r\left(
 \frac{q(y)-q(\alpha)}{y-\alpha}
 +a\frac{h(y)-h(\alpha)-(y-\alpha)h'(\alpha)}
          {(y-\alpha)^2}
 \right).
\endaligned
\end{equation}
Since $h_m(y)=mh(y)-\alpha(m-1)g_m(y),$
its first-order variation is
\begin{equation*}\label{B8}
\aligned
\delta h_m(y)
=mq(y)-a(m-1)g_m(y)-\alpha m(m-1)\left(
 \frac{q(y)-q(\alpha)}{y-\alpha}
 +a\frac{h(y)-h(\alpha)-(y-\alpha)h'(\alpha)}
          {(y-\alpha)^2}
\right).
\endaligned
\end{equation*}
Substituting these parameters into the action formulas of $M$, write $\rho_{\varepsilon}(X)
=\rho_{\boldsymbol{\omega}}(X)+\varepsilon \tau_{\boldsymbol{v}}(X), X\in \mathcal{G}.$
Then $\tau_{\boldsymbol{v}}$ satisfies (\ref{formula 4.1}). In this way, we obtain an explicit family of $1$-cocycles. 

\begin{prop}\label{prop 4.1}
 $\tau_{\boldsymbol{v}}$ is a $1$-coboundary if and only if 
\begin{equation}\label{formula 4.8}
\aligned
v=a=s=0,b=\frac{u\beta}{2}, c=\frac{u\gamma}{2},q(\alpha)=0.
\endaligned
\end{equation}

\end{prop}
\begin{proof}
\rm{(1)} Suppose that there exists $T\in \operatorname{End}_{\C}(M)$ such that $\tau_{\boldsymbol{v}}(X)=[\rho_{\boldsymbol{\omega}}(X),T]$ for any $X\in\mathcal{G}$. Since  the actions of $L_0$ and $W_0$ are respectively multiplication by $x$ and $y$, and $\tau_{\boldsymbol{v}}(L_{0})=\tau_{\boldsymbol{v}}(W_{0})=0_{\C[x,y]}$, we can obtain that $[x,T]=[y,T]=0_{\C[x,y]}$ and $T$ is multiplication by some polynomial $p(x,y)\in\C[x,y]$.

For $f\in\C[x,y]$, the coefficient of $\varepsilon$ in the action of $W_{m,i}$ is $\tau_{\boldsymbol{v}}(W_{m,i})f
=\lambda^m\zeta^i
\big((mu+iv)(y-m\alpha)-ma\big)
e^{m\frac{\partial}{\partial{x}}}f.$
On the other hand, $[\rho_{\boldsymbol{\omega}}(W_{m,i}),T]f
=\lambda^m\zeta^i(y-m\alpha)\Delta_mp
e^{m\frac{\partial}{\partial{x}}}f. $
Consequently,
\begin{equation}\label{formula 4.10}
\aligned
(mu+iv)(y-m\alpha)-ma=(y-m\alpha)\Delta_mp.
\endaligned
\end{equation}
Taking $(m,i)=(0,1)$ and $(1,0)$ in (\ref{formula 4.10}), and then setting $y=\alpha$ in the second resulting identity, we have $v=a=0$. It follows that $\Delta_1p=u$, and hence
\begin{equation}\label{formula 4.11}
\aligned
p(x,y)=ux+p_0(y)
\endaligned
\end{equation}
for some $p_0\in\C[y]$ 
. 
Similarly, comparison on $L_{0,1},B_{r,i}$ and $I_{r,i}$, respectively, gives $s=0,c=\frac{u\gamma}{2},b=\frac{u\beta}{2}.$

For $f\in\M$, by (\ref{formula 4.11}), comparison on $L_{1}$ gives $\tau_{\boldsymbol{v}}(L_{1})f
=\lambda\Big(
u\bigl(x+h(y)\bigr)
+u(y-\alpha)\frac{\partial }{\partial y}
+q(y)
\Big)e^{\frac{\partial}{\partial{x}}}f,$
and
$\relax[\rho_{\boldsymbol{\omega}}(L_{1}),T]f
=\lambda\Big(
u\bigl(x+h(y)\bigr)
+u(y-\alpha)\frac{\partial }{\partial y}
+(y-\alpha)\frac{d}{dy}p_0(y)
\Big)e^{\frac{\partial}{\partial{x}}}f.$
Comparison gives $q(y)=(y-\alpha)\frac{d}{dy}p_0(y)$
and hence $q(\alpha)=0$.

\rm{(2)} Conversely, suppose that (\ref{formula 4.8}) holds. Then $q(y)=(y-\alpha)r_0(y)$ for some $r_0(y)\in\C[y]$. Choose $p_0(y)\in\C[y]$ such that $\frac{d}{dy}p_0(y)=r_0(y)$, and let $T$ be multiplication by $p(x,y)=ux+p_0(y)$. Then 
\begin{equation}\label{formula 4.14.0}
\aligned
\Delta_tp=ut
\endaligned
\end{equation}
for any $t\in\C$. The variation of $g_t(y)$ in the direction $q(y)$ is $t\frac{q(y)-q(\alpha)}{y-\alpha}=tr_0(y)$, and hence the corresponding variation of $h_m(y)$ is
\begin{equation}\label{formula 4.14.1}
\aligned
m q(y)-\alpha(m-1)m r_0(y)
&=m(y-\alpha)r_0(y)-\alpha m(m-1)r_0(y)
=m(y-m\alpha)r_0(y).
\endaligned
\end{equation}
Using (\ref{formula 4.14.0}) and (\ref{formula 4.14.1}), direct calculation gives, for every $f\in M,m,i\in\Z,r\in\Z+\frac{1}{2}$, $[ \rho_{\boldsymbol{\omega}}(W_{m,i}),T]f=\tau_{\boldsymbol{v}}(W_{m,i})f$, $[ \rho_{\boldsymbol{\omega}}(B_{r,i}),T]f=\tau_{\boldsymbol{v}}(B_{r,i})f$
and
\begin{equation*}
\aligned
\relax[\rho_{\boldsymbol{\omega}}(L_{m,i}),T]f
=&\lambda^m\zeta^im\left(
u\bigl(x+h_m(y)+i\sigma(y-m\alpha)\bigr)+(y-m\alpha)r_0(y)
+mu(y-m\alpha)
\frac{\partial }{\partial y}
\right)e^{m\frac{\partial}{\partial{x}}}f,
\endaligned
\end{equation*}
which is exactly $\tau_{\boldsymbol{v}}(L_{m,i})f$. Similarly, 
\begin{equation*}
\aligned
\relax[\rho_{\boldsymbol{\omega}}(I_{r,i}),T]f
=&\lambda^{r-\frac12}\zeta^ir\left(
u\bigl(\beta+\gamma g_r(y)+i\sigma\gamma\bigr)+\gamma r_0(y)
+ur\gamma\frac{\partial }{\partial y}
\right)e^{r\frac{\partial}{\partial{x}}}f=\tau_{\boldsymbol{v}}(I_{r,i})f.
\endaligned
\end{equation*}
Both sides vanish on the central generators. Therefore,
$\tau_{\boldsymbol{v}}(X)=[\rho_{\boldsymbol{\omega}}(X),T]
\text{\ for all }X\in\mathcal G$,
so $\tau_{\boldsymbol{v}}$ is a 1-coboundary.
\end{proof}


\subsection{Classification of self-extensions}

\begin{lemm}\label{lemm 4.1}
    Every $1$-cocycle $\tau$ is cohomologous to a $1$-cocycle $\tau'$ satisfying 
    \begin{equation}\label{formula 4.15}
\aligned
\tau'(L_{0})=\tau'(W_{0})=0_{\C[x,y]}.
\endaligned
\end{equation}
\end{lemm}
\begin{proof}
Let $0\to M\xrightarrow{\iota} E_{\tau}\xrightarrow{\pi} M\to 0$
be the self-extension determined by $\tau$, where $\iota(f_1)=(f_1,0),\pi(f_1,f_2)=f_2$ for any $f_1,f_2\in\C[x,y]$.
Recall that $U(\mathcal{h})\cong \C[x,y]$ and $M$ is a free $U(\mathcal{h})$-module of rank one generated by $\mathbf{1}$. In particular, $M$ is projective as a $U(\mathcal{h})$-module. Therefore, after restricting the above exact sequence to $U(\mathcal{h})$, there exists a $U(\mathcal{h})$-module homomorphism $\mathbf{s}:M\to E_{\tau}$ such that $\pi\circ \mathbf{s}=\operatorname{Id}_{\C[x,y]}$. We make this splitting explicit. Since $\pi(\mathbf{s}(f))=f$ for every $f\in M$, there exists a uniquely determined $\C$-linear operator $S\in \operatorname{End}_{\C}(M)$ such that $\mathbf{s}(f)=(S(f),f)^{\operatorname{T}}$. Because $\mathbf{s}$ is $U(\mathcal{h})$-linear, for $X\in\{L_{0}, W_{0} \}$ and $f\in M$, we have $X(\mathbf{s}(f))=\mathbf{s}(X(f))$. By (\ref{formula 4.3}), we have 
\begin{equation*}
\aligned
X\big(\mathbf{s}(f)\big)
&=X\begin{pmatrix}
    S(f)\\
    f
\end{pmatrix}
=\begin{pmatrix}
    \rho_{\boldsymbol{\omega}}(X)S(f)+\tau(X)f\\
    \rho_{\boldsymbol{\omega}}(X)f
\end{pmatrix},\ \mathbf{s}\big(X (f)\big)
=\mathbf{s}\bigl(\rho_{\boldsymbol{\omega}}(X)f\bigr)
=\begin{pmatrix}
    S\rho_{\boldsymbol{\omega}}(X)f\\
    \rho_{\boldsymbol{\omega}}(X)f
\end{pmatrix}.
\endaligned
\end{equation*}
Comparison of the first components gives
$\rho_{\boldsymbol{\omega}}(X)S(f)+\tau(X)f
=S\rho_{\boldsymbol{\omega}}(X)f$.
Therefore,
\begin{equation}\label{formula 4.16.1}
\aligned
\tau(X)+[\rho_{\boldsymbol{\omega}}(X),S]=0_{\C[x,y]}.
\endaligned
\end{equation}
Now define
$\tau'(X)
=\tau(X)+[\rho_{\boldsymbol{\omega}}(X),S],
 X\in\mathcal G.$
Since $X\mapsto[\rho_{\boldsymbol{\omega}}(X),S]$ is the $1$-coboundary determined by S, the maps $\tau$ and $\tau'$ are cohomologous. Moreover, applying (\ref{formula 4.16.1}) to $X=L_{0}$ and $X=W_{0}$, respectively, gives
$\tau'(L_{0})=0_{\C[x,y]},
\tau'(W_{0})=0_{\C[x,y]}.$
Thus every $1$-cocycle is cohomologous to a $1$-cocycle satisfying (\ref{formula 4.15}).


\end{proof}

\begin{lemm}\label{lemm4.2}
    Suppose that $\tau$ satisfies (\ref{formula 4.15}). Then there exist polynomials $\Tilde{l}_{m,i},\Tilde{w}_{m,i},\Tilde{i}_{r,i},\Tilde{b}_{r,i}\in\C[x,y]$ such that
\begin{eqnarray}
&\tau(W_{m,i})=\Tilde{w}_{m,i}e^{m\frac{\partial}{\partial{x}}},\ &\tau(L_{m,i})=\left(\Tilde{l}_{m,i}+m\Tilde{w}_{m,i}\frac{\partial}{\partial{y}}\right)e^{m\frac{\partial}{\partial{x}}},\label{formula 4.16}\\ &\tau(B_{r,i})=\Tilde{b}_{r,i}e^{r\frac{\partial}{\partial{x}}},\ &\tau(I_{r,i})=\left(\Tilde{i}_{r,i}+r\Tilde b_{r,i}\frac{\partial}{\partial{y}}\right)e^{r\frac{\partial}{\partial{x}}}\label{formula 4.16.2}
\end{eqnarray}
for any $m,i\in\Z,r\in\Z+\frac{1}{2}$. Moreover, for every central generator $C_{k,i}$, there exists $\Tilde{c}_{k,i}\in\C$ such that $\tau(C_{k,i})=\Tilde{c}_{k,i}\operatorname{Id}_{\C[x,y]}$, where $i\in\Z,k=1,2,3,4$.
    
\end{lemm}
\begin{proof}
    Let $X$ be a homogeneous generator of degree $d$, so that $[L_{0},X]=-dX$.
 Applying (\ref{formula 4.1}) and $\tau(L_0)=0_{\C[x,y]},\rho_{\boldsymbol{\omega}}(L_0)=x\in\operatorname{End}_{\C}(M)$, we have $\big[x,\tau(X)   \big]=-d\tau(X)$.
By Taylor's formula (\ref{Taylor's formula}),   $\big[x,e^{t\frac{\partial}{\partial{x}}}\big]=-te^{t\frac{\partial}{\partial{x}}}$ for every $t\in\C$. Set $T_X=\tau(X)e^{-d\frac{\partial}{\partial{x}}}$ for such a homogeneous generator $X$. Since $-d\tau(X)=\big[x,\tau(X)   \big]$, we can obtain that $\big[x,T_X   \big]=0_{\C[x,y]}.$
Then $T_X$ commutes with multiplication by $x$, and hence is a $\C[x]$-linear endomorphism of $\C[x,y]$. 
Since $\rho_{\boldsymbol{\omega}}(W_0)=y,\tau(W_0)=0_{\C[x,y]}$, the cocycle identity (\ref{formula 4.1}) applied to $[W_0,W_{m,i}]=0,[W_0,L_{m,i}]=-mW_{m,i}$ for any $m,i\in\Z$, gives $\big[y,\tau(W_{m,i})  \big]=0_{\C[x,y]},\big[y,\tau(L_{m,i})  \big]=-m\tau(W_{m,i})$. Likewise, from $[W_{0}, B_{r,i}]=0,[W_{0},I_{r,i}]=-rB_{r,i}$, we have $\big[y,\tau(B_{r,i})  \big]=0_{\C[x,y]},\big[y,\tau(I_{r,i})  \big]=-r\tau(B_{r,i})$.

Write $\tau(W_{m,i})=T_{W_{m,i}}e^{m\frac{\partial}{\partial{x}}}$. Then we obtain $[y,T_{W_{m,i}}]e^{m\frac{\partial}{\partial{x}}}=0_{\C[x,y]}$, and hence $[y,T_{W_{m,i}}]=0_{\C[x,y]}$, which means that the operator $T_{W_{m,i}}$ commutes with $y$. Therefore, $T_{W_{m,i}}(f)=fT_{W_{m,i}}(1)$ for any $f\in\C[x,y]$. It must therefore be multiplication by $\Tilde{w}_{m,i}:=T_{W_{m,i}}(1)\in\C[x,y]$ for any $m,i\in\Z$. Thus, $\tau(W_{m,i})=\Tilde{w}_{m,i}e^{m\frac{\partial}{\partial{x}}}$. Similarly, we have $\tau(B_{r,i})=\Tilde{b}_{r,i}e^{r\frac{\partial}{\partial{x}}}$ for any $r\in\Z+\frac{1}{2},i\in\Z$. Now write $\tau(L_{m,i})=T_{L_{m,i}}e^{m\frac{\partial}{\partial{x}}}$. We can obtain $[y,T_{L_{m,i}}]=-m\Tilde{w}_{m,i}$ by $\big[y,\tau(L_{m,i})  \big]=-m\tau(W_{m,i})$. Since $\Big[ y,\frac{\partial}{\partial{y}}  \Big]=-1$, we have $\left[ y,  T_{L_{m,i}}-m \Tilde{w}_{m,i}  \frac{\partial}{\partial{y}}  \right]=0_{\C[x,y]}$.
Moreover, both $T_{L_{m,i}}$ and $\Tilde{w}_{m,i}\frac{\partial}{\partial{y}}$ commute with $x$. Hence, $T_{L_{m,i}}-m\Tilde{w}_{m,i}\frac{\partial}{\partial{y}}$ commutes with both $x$ and $y$, and is therefore multiplication by a polynomial $\Tilde{l}_{m,i}\in\C[x,y]$. Consequently, $\tau(L_{m,i})=\Big( \Tilde{l}_{m,i}+m\Tilde{w}_{m,i}\frac{\partial}{\partial{y}}  \Big)e^{m\frac{\partial}{\partial{x}}}$. By the same argument, we have $\tau(I_{r,i})=\Big( \Tilde{i}_{r,i}+r\Tilde{b}_{r,i}\frac{\partial}{\partial{y}}  \Big)e^{r\frac{\partial}{\partial{x}}}$. 

Finally, let $C_{k,i}$ be a central generator. Since $\rho_{\boldsymbol{\omega}}(C_{k,i})=0_{\C[x,y]}$ and $[C_{k,i},X]=0$ for any $X\in \mathcal{G},i\in\Z,k=1,2,3,4$, (\ref{formula 4.1}) gives $\big[\tau(C_{k,i}),\rho_{\boldsymbol{\omega}}(X)  \big]=0$.
Thus, $\tau(C_{k,i})$ commutes with every operator $\rho_{\boldsymbol{\omega}}(X)$. In particular, it commutes with $\rho_{\boldsymbol{\omega}}(L_0)=x$ and $\rho_{\boldsymbol{\omega}}(W_0)=y$. Hence, $\tau(C_{k,i})$ is multiplication by some polynomial $\Tilde{c}_{k,i}\in\C[x,y]$. Since it also commutes with $\rho_{\boldsymbol{\omega}}(W_{1,0})=\lambda(y-\alpha)e^{\frac{\partial}{\partial{x}}}$, we have, for any $0\neq f(x,y)\in M$, $0=\big[  \rho_{\boldsymbol{\omega}}(W_{1,0}), \Tilde{c}_{k,i}(x,y) \big]f(x,y)=\lambda(y-\alpha)\Delta_1\Tilde{c}_{k,i}(x,y)f(x+1,y),$
which means $\Tilde{c}_{k,i}(x,y)=\Tilde{c}_{k,i}(y)\in\C[y]$. From 
    $$0=\big[  \rho_{\boldsymbol{\omega}}(L_{1,0}), \Tilde{c}_{k,i}(y) \big]f(x,y)=\lambda(y-\alpha)\frac{d}{dy}\Tilde{c}_{k,i}(y)f(x+1,y),$$
we have $\frac{d}{dy}\Tilde{c}_{k,i}(y)=0$, and hence $\Tilde{c}_{k,i}(y)=\Tilde{c}_{k,i}\in\C$. Therefore, $\tau(C_{k,i})=\Tilde{c}_{k,i}\operatorname{Id}_{\C[x,y]}$ for $i\in\Z,k=1,2,3,4$.
\end{proof}

\begin{lemm}\label{tilde w}
    For any $u\in\C$, one can choose a cohomologous normalized cocycle and $v,a\in\C$ such that 
 \begin{equation}\label{B4}
\aligned
\Tilde w_{m,i}=\delta w_{m,i}
=\lambda^m\zeta^i
 \left((mu+iv)(y-m\alpha)-ma\right).
\endaligned
\end{equation}
\end{lemm}
\begin{proof}
Applying (\ref{formula 4.1}) to $[L_{1,0},W_{1,0}]=0$, we have $$(x+h(y))\Delta_1\Tilde{w}_{1,0}(x,y)+(y-\alpha)\left(\frac{\partial}{\partial{y}}\Tilde{w}_{1,0}(x+1,y)-\Delta_1\Tilde{l}_{1,0}(x,y)    \right)=0.$$
Thus, $\Delta_1\Tilde{w}_{1,0}(x,\alpha)=0$. Define $\Delta_1p=\frac{\Tilde{w}_{1,0}(x,y)-\lambda(u(y-\alpha)-a)}{\lambda(y-\alpha)}\in\C[x,y]$, where $a=-\lambda^{-1}\Tilde{w}_{1,0}(0,\alpha)$. Replacing $\tau$ by $\tau-[\rho_{\boldsymbol{\omega}}(\cdot),p]$, we can obtain $\Tilde{w}_{1,0}=\lambda(u(y-\alpha)-a)\in\C[y]$. Furthermore, $[W_{m,i},W_{1,0}]=0$ reduces to $\lambda(y-\alpha)\Delta_1\Tilde{w}_{m,i}(x,y)=0$. Therefore, $\Tilde{w}_{m,i}\in\C[y]$.
We next apply (\ref{formula 4.1}) to $[L_{m,i},W_{n,j}]
=(m-n)W_{m+n,i+j}
+\frac{m^3-m}{12}\delta_{m+n,0}C_{2,i+j}.$
For $m+n\neq0$ or $m=0,\pm 1$, we have
\begin{equation}\label{B2}
\aligned
m(y-m\alpha)
 \frac{d}{d y}\Tilde w_{n,j}-n\Tilde w_{n,j}
-\lambda^{\,n-m}\zeta^{\,j-i}((y-n\alpha)
       \Delta_n\Tilde l_{m,i}    
-m\Tilde w_{m,i})
=(m-n)\lambda^{-m}\zeta^{-i}
 \Tilde w_{m+n,i+j}.
\endaligned
\end{equation}
Taking $(n,j)=(1,0)$ shows that $\Delta_1\Tilde l_{m,i}$ is independent of $x$, so $\frac{\partial^2}{\partial{x^2}}\Tilde l_{m,i}=0$ and $\Delta_n\Tilde l_{m,i}=n\frac{\partial}{\partial{x}}\Tilde l_{m,i}$. Taking $m=n\ne 0,j=0$ gives 
\begin{equation}\label{B2.0}
\aligned
(y-m\alpha)\frac{\partial}{\partial{x}}\Tilde l_{m,i}=\Tilde w_{m,i}+\zeta^i\Big((y-m\alpha)\frac{d}{d{y}}\Tilde w_{m,0}-\Tilde w_{m,0}\Big).
\endaligned
\end{equation}
 Taking $m=0,n=1,j=0$ in (\ref{B2}), we have $\Tilde{w}_{1,i}=\zeta^i\Tilde{w}_{1,0}+\lambda(y-\alpha)\frac{\partial}{\partial{x}}\Tilde{l}_{0,i}$. Substituting this relation into (\ref{B2}) with $(m,n,i,j)=(1,1,0,i)$ gives $(y-\alpha)^2\frac{\partial}{\partial{y}}\frac{\partial}{\partial{x}}\Tilde{l}_{0,i}=0$. Returning to (\ref{B2}) with $(m,n)=(0,1)$, we obtain $\frac{\partial}{\partial{x}}\Tilde{l}_{0,i+j}=\zeta^i\frac{\partial}{\partial{x}}\Tilde{l}_{0,j}+\zeta^j\frac{\partial}{\partial{x}}\Tilde{l}_{0,i}$. Hence, there is a constant $v\in\C$ such that $\frac{\partial}{\partial{x}}\Tilde{l}_{0,i}=iv\zeta^i$. Then letting $m=0,n\neq0,j=0$ in (\ref{B2}), we have
\begin{equation}\label{B2.1}
\aligned
\Tilde{w}_{n,i}=\zeta^i\Tilde{w}_{n,0}+iv\lambda^n\zeta^i(y-n\alpha).
\endaligned
\end{equation}
By (\ref{B2}) and (\ref{B2.0}), we can obtain $\frac{\partial}{\partial{x}}\Tilde{l}_{1,0}=\lambda u$. Taking $(m,n,i,j)=(1,-1,0,0)$ in (\ref{B2}) gives $\lambda[(y-\alpha)\frac{d}{d{y}}\Tilde{w}_{-1,0}+\Tilde{w}_{-1,0}]+2uy-a=0$. Then we have $\Tilde{w}_{-1,0}=\lambda^{-1}(-u(y+\alpha)+a)$. Taking $(m,n,i,j)=(-1,2,0,0)$ and using $\frac{\partial}{\partial{x}}\Tilde{l}_{-1,0}=\frac{d}{d{y}}\Tilde{w}_{-1,0}=-\lambda^{-1}u$, we have  $\Tilde{w}_{2,0}=\lambda^{2}(2u(y-2\alpha)-2a)$. The choice $(m,n,i,j)=(1,-2,0,0)$ similarly gives $\Tilde{w}_{-2,0}=\lambda^{-2}(-2u(y+2\alpha)+2a)$.
By induction on $m\neq0$, we obtain 
\begin{equation}\label{B2.2}
\aligned
\Tilde{w}_{m,0}=\lambda^m\left(mu(y-m\alpha)-ma\right).
\endaligned
\end{equation}
  Finally, (\ref{B2}) with $(m,n,i,j)=(1,-1,0,i)$ gives $\Tilde{w}_{0,i}=iv\zeta^iy$. Together with (\ref{B2.1}) and (\ref{B2.2}), this proves (\ref{B4}). 
\end{proof}

\begin{lemm}\label{tilde l}
    For the representative in the preceding lemma there exist $s\in\C$ and $q\in\C[y]$ such that
    \begin{equation}\label{B6}
\begin{aligned}
\Tilde l_{m,i}=\delta l_{m,i}
=\lambda^m\zeta^i\left(
 (mu+iv)\bigl(x+h_m(y)+i\sigma(y-m\alpha)\bigr)+
 \delta h_m(y)+is(y-m\alpha)-im\sigma a
\right).
\end{aligned}
\end{equation}
\end{lemm}
\begin{proof}
    By Lemma \ref{tilde w}, we can verify 
\begin{equation}\label{B6.0}
\begin{aligned}
\frac{\partial}{\partial{x}}\Tilde l_{m,i}
=\lambda^m\zeta^i(mu+iv).
\end{aligned}
\end{equation}
   Consequently, $q(y):=\lambda^{-1}\Tilde l_{1,0}(x,y)-u(x+h(y))$ 
is independent of $x$. By (\ref{B6.0}), $P_i:=\zeta^{-i}\Tilde l_{0,i}(x,y)-ivx-i^2v\sigma y$ is independent of $x$. Applying the cocycle identity
to $[L_{m,i},L_{n,j}]
=(m-n)L_{m+n,i+j}
+\frac{m^3-m}{12}\delta_{m+n,0}C_{1,i+j}$
and comparing the multiplication parts, for $m+n\neq 0$ or $m=0,\pm 1$, 
\begin{equation}\label{B5}
\begin{aligned}
&(m-n)\Tilde l_{m+n,i+j}(x,y) \\
=&A_{m,i,m,0}\Tilde l_{n,j}(x+m,y)-A_{m,i,m,n}\Tilde l_{n,j}(x,y)+A_{n,j,n,m}\Tilde l_{m,i}(x,y)-A_{n,j,n,0}\Tilde l_{m,i}(x+n,y)\\
&+B_{m,i,m}\frac{\partial}{\partial y}\Tilde l_{n,j}(x+m,y)-B_{n,j,n}\frac{\partial}{\partial y}\Tilde l_{m,i}(x+n,y)-nC_{m,i}\Tilde w_{n,j}(y)+mC_{n,j}\Tilde w_{m,i}(y),
\end{aligned}
\end{equation}
 where $A_{p,q,k,l}=\lambda^p\zeta^q(x+l+h_k+q\sigma(y-k\alpha)),B_{p,q,k}=\lambda^p\zeta^qk(y-k\alpha),C_{p,q}=\lambda^p\zeta^q\bigl(h_p'(y)+q\sigma\bigr)$.
   Taking $(m,n)=(1,1),(1,-1),(0,1)$ in (\ref{B5}), gives respectively
   \[
       \frac{d^2}{dy^2}P_i=0,\ y\frac{d}{dy}P_i=P_i, \ P_{i+j}=P_i+P_j.
   \]
   Therefore, there is a constant $s\in\C$ such that
   \begin{equation}\label{B6.1}
\begin{aligned}
\Tilde{l}_{0,i}=\zeta^i\left(ivx+i^2v\sigma y+isy\right).
\end{aligned}
\end{equation}
Define 
\begin{equation}\label{B6.2}
\begin{aligned}
d_{m,i}=\Tilde{l}_{m,i}-\lambda^m\zeta^i[(mu+iv)(x+h_m(y)+i\sigma(y-m\alpha))+\delta h_m+is(y-m\alpha)-im\sigma a].
\end{aligned}
\end{equation}
By the definitions of $q$ and $s$, (\ref{B6.1}) and (\ref{B6.2}) give $d_{1,0}=d_{0,i}=0$. Direct substitution into the cocycle identity gives the homogeneous differential equation
\begin{equation*}\label{B6.3}
\begin{aligned}
(m-n)d_{m+n,i+j}=\lambda^m\zeta^i\Big(-nd_{n,j}+m(y-m\alpha)\frac{d}{dy}d_{n,j}\Big)+\lambda^n\zeta^j\Big(md_{m,i}-n(y-n\alpha)\frac{d}{dy}d_{m,i}\Big).
\end{aligned}
\end{equation*}
Taking $(m,n)=(1,-1)$ gives $\frac{d}{dy}((y-\alpha)d_{-1,0})$, hence $d_{-1,0}=0$. The choice $(-1,2)$ and $(1,-2)$ give $(y+\alpha)\frac{d}{dy}d_{2,0}+2d_{2,0}=0$ and $(y-\alpha)\frac{d}{dy}d_{-2,0}+2d_{-2,0}=0$. Their only polynomial solutions are zero. By induction on $m\neq0$, $d_{m,0}=0$. Its $m=0$ specialization and $d_{0,i}=0$ give $d_{n,i+j}=\zeta^id_{n,j}=0$.
\end{proof}

\begin{lemm}\label{tilde b}
    For the same representative there exists $c\in\C$ such that
    \begin{equation}\label{B11}
\aligned
\Tilde b_{r,i}=\delta b_{r,i}
=\lambda^{r-\frac12}\zeta^i
 \left(
   \gamma\left(\left(r-1/2\right)u+iv\right)+c
 \right).
\endaligned
\end{equation}
\end{lemm}
\begin{proof}
    By Lemma \ref{lemm4.2}, 
the cocycle identity applied to $[W_{1,0},B_{r,j}]=0$ gives $\lambda(y-\alpha)\Delta_1\Tilde{b}_{r,j}=0$. Hence, $\Tilde{b}_{r,j}\in\C[y]$.
Next, applying the cocycle identity to $[L_{m,i},B_{r,j}]=-rB_{m+r,i+j}$ gives
{\begin{equation}\label{B10}
\aligned
&\lambda^m\zeta^i
 \Big(
   m(y-m\alpha)\frac{d}{dy}\Tilde b_{r,j}-r\Tilde b_{r,j}
 \Big)
-r\gamma\lambda^{m+r-\frac12}\zeta^{i+j}(mu+iv)
=-r\Tilde b_{m+r,i+j}.
\endaligned
\end{equation}}
Define $d_{r,i}=\Tilde b_{r,i}-\gamma\lambda^{r-\frac{1}{2}}\zeta^i[(r-1/2)u+iv].$
Substituting into (\ref{B10}) gives
\begin{equation}\label{B11.1}
\aligned
\lambda^m\zeta^i\Big(m(y-m\alpha)\frac{d}{dy}d_{r,j}-rd_{r,j}\Big)=-rd_{m+r,i+j}.
\endaligned
\end{equation}
Taking $m=0$ yields $d_{r,i}=\zeta^id_{r,0}$. Taking $i=j=0$ and $(m,r)=(1,-\frac{1}{2}),(-1,\frac{1}{2})$ give 
\[
    \lambda\Big(2(y-\alpha)\frac{d}{dy}d_{-\frac{1}{2},0}+d_{-\frac{1}{2},0}\Big)=d_{\frac{1}{2},0},\ \lambda^{-1}\Big(2(y+\alpha)\frac{d}{dy}d_{\frac{1}{2},0}+d_{\frac{1}{2},0}\Big)=d_{-\frac{1}{2},0}.
\]
Then we have $(y^2-\alpha^2)\frac{d^2}{dy^2}d_{\frac{1}{2},0}+(2y-\alpha)\frac{d}{dy}d_{\frac{1}{2},0}=0$ and hence $d_{\frac{1}{2},0}=c\in\C$. Taking $(r,j)=(\frac{1}{2},0)$ in (\ref{B11.1}) gives $d_{m+\frac{1}{2},i}=c\lambda^m\zeta^i$, proving (\ref{B11}).
\end{proof}

\begin{lemm}\label{tilde i}
    For the same representative there exists $b\in\C$ such that
    \begin{equation}\label{B13red}
\aligned
\Tilde i_{r,i}=\delta i_{r,i}
=\lambda^{r-\frac12}\zeta^i
\left(\left(\left(r-1/2\right)u+iv\right)
  \left(\beta+\gamma g_r(y)+i\sigma\gamma\right)+b+c g_r(y)
+\gamma \delta g_r(y)
        +i(s\gamma+\sigma c)
\right).
\endaligned
\end{equation}
\end{lemm}
\begin{proof}
    Define $$d_{r,i}=
\Tilde i_{r,i}
-\lambda^{r-\frac12}\zeta^i\left(
\left(\left(r-1/2\right)u+iv\right)
  \bigl(\beta+\gamma g_r(y)+i\sigma\gamma\bigr)+c g_r(y)
+\gamma \delta g_r(y)
        +i(s\gamma+\sigma c)
\right).$$
Applying (\ref{formula 4.1}) to $[W_{m,i},I_{r,j}]=-rB_{m+r,i+j}$ gives $\lambda^m\zeta^i(y-m\alpha)\Delta_md_{r,j}=0$. Taking $m=1$ shows that $d_{r,j}\in\C[y]$.
Next, applying the cocycle identity to $[L_{m,i},I_{r,j}]=-rI_{m+r,i+j}$ gives (\ref{B11.1}). Hence, $d_{\frac{1}{2},0}=b\in\C$ This proves (\ref{B13red}).
\end{proof}

\begin{theo}\label{theo 4.1}
    For every $\boldsymbol{\omega}$, every $1$-cocycle is cohomologous to a $1$-cocycle $\tau_{\boldsymbol{v}}$, i.e. 
\begin{equation}\label{formula 4.17}
\aligned
H^1(\mathcal{G},\operatorname{End}_{\C}(M))\cong \operatorname{Ext}^1_{U(\mathcal{G})}(M,M)\cong \C^6.
\endaligned
\end{equation} 
\begin{proof}
Lemmas (\ref{tilde w}), (\ref{tilde l}), (\ref{tilde b}) and (\ref{tilde i}) give the normalized operator form (\ref{formula 4.16}) and (\ref{formula 4.16.2}).
It remains to determine the central terms. Applying the
cocycle identity to  $[L_{m,i},L_{-m,j}]$ and
$[L_{m,i},W_{-m,j}]$, 
$[I_{r,i},I_{-r,j}]$ and $[I_{r,i},B_{-r,j}]$, and substituting (\ref{B4}), (\ref{B6}), (\ref{B11}) and (\ref{B13red}), we obtain 
   $\Tilde c_{1,i+j}
=\Tilde c_{2,i+j}
=\Tilde c_{3,i+j}
=\Tilde c_{4,i+j}=0.$
These calculations show that there exist $u,v,a,b,c,s\in\mathbb C, q(y)\in\mathbb C[y],$
and an operator $T\in\operatorname{End}_{\mathbb C}(M)$
such that
\begin{equation*}\label{formula 4.19.1}
\aligned
\tau(X)=\tau_{(u,v,a,b,c,s,q)}(X)
       +[\rho_{\boldsymbol{\omega}}(X),T]
\endaligned
\end{equation*} 
for every generator $X$ of $\mathcal G$. Hence $\tau$ is
cohomologous to 
$\tau_{(u,v,a,b,c,s,q)}$.
By Proposition \ref{prop 4.1}, two parameter directions $\boldsymbol{v},\boldsymbol{v}'$ determine the same
cohomology class if and only if $v=v', a=a', s=s',$
and $b-\frac{u\beta}{2}
=b'-\frac{u'\beta}{2},
c-\frac{u\gamma}{2}
=c'-\frac{u'\gamma}{2},
q(\alpha)=q'(\alpha).$
Consequently, the assignment 
$[\tau_{\boldsymbol{v}} ] \mapsto \left(
v,\ a,\ b-\frac{u\beta}{2},\
c-\frac{u\gamma}{2},\ s,\ q(\alpha)
\right)$
defines a linear isomorphism $H^1(\mathcal{G},\operatorname{End}_{\C}(M))\cong \C^6.$
\end{proof}

\end{theo}

\begin{coro}\label{coro 4.5}
    Suppose that $(\alpha,\gamma)\neq(0,0)$.  For every nonzero class $[\tau]\in \operatorname{Ext}^1_{U(\mathcal{G})}(M,M)$, the corresponding self-extension (\ref{formula 4.4}) is nonsplit, and $E_{\tau}$ has composition length two and is indecomposable.
\end{coro}
\begin{proof}
By Theorem \ref{theo 3.3}, the assumption $(\alpha,\gamma)\neq(0,0)$ implies that
$M$ is irreducible. Since $[\tau]\neq0$ in
$\operatorname{Ext}_{U(\mathcal G)}^{1}(M,M)$, the self-extension $0\longrightarrow M\xrightarrow{\iota}E_{\tau}
\xrightarrow{\pi}M\longrightarrow0$
does not split; otherwise, by (\ref{formula 4.5}), $\tau$ would be a $1$-coboundary and
hence $[\tau]=0$.
Moreover, $\iota(M)$ is a simple submodule of $E_{\tau}$, and $E_{\tau}/\iota(M)\cong M$
is also simple. Therefore, $0\subsetneq\iota(M)\subsetneq E_{\tau}$
is a composition series, so $E_{\tau}$ has composition length two.
Finally, if $E_{\tau}$ were decomposable, then its composition length two would force it to be a direct sum of two simple modules. Hence, $E_{\tau}$ would be semisimple, and the above self-extension would split, a contradiction.  Thus $E_{\tau}$ is
indecomposable.
\end{proof}

\begin{rema}\label{rema4.6}
The cocycle computation in Theorem \ref{theo 4.1} does not use the simplicity of $M$; the simplicity of $M$ (equivalently, the assumption $(\alpha,\gamma)\neq(0,0)$) is used only in Corollary \ref{coro 4.5} to conclude that the self-extensions are nonsplit and indecomposable.
\end{rema}

\section{\texorpdfstring{ A new class of simple $\mathcal{H}$-modules 
from $\Omega(\lambda, \zeta, \alpha,\beta,\gamma,\sigma,h)$}{A new class of irreducible modules over hat{mathcal{H}}from Omega(lambda,zeta,alpha,beta,gamma,h)}}\label{section 5}

\subsection{\texorpdfstring{ $\mathcal{H}$-modules $\Omega_{\mathcal{H}}(\lambda,\zeta,\alpha,\beta,\gamma,\sigma,h)$ and the $\sigma$-type}{1}}

By restricting the action in Proposition \ref{prop 3.1} to $\mathcal{H}$, we obtain an $\mathcal{H}$-module, denoted by $$\Omega_{\mathcal{H}}(\lambda,\zeta,\alpha,\beta,\gamma,\sigma,h)=\operatorname{Res}_{\mathcal{H}}^{\mathcal{G}}\Omega(\lambda, \zeta,\alpha,\beta,\gamma,\sigma,h)$$ as a vector space. And the action of $\mathcal{H}$ is given by
\begin{eqnarray}\label{act. def 6.1}
&L_{m,i} f(x,y)&=\lambda^m\zeta^i(x+h_{m}(y)+i\sigma(y-m\alpha))f(x+m,y)+m\lambda^m\zeta^i(y-m\alpha)\frac{\partial }{\partial {y}}{f(x+m,y)},\nonumber\\
&I_{r,i} f(x,y)&=\lambda^{r-\frac{1}{2}}\zeta^i(\beta+\gamma g_{r}(y)+i\sigma \gamma)f(x+r,y)+r\gamma\lambda^{r-\frac{1}{2}}\zeta^i\frac{\partial }{\partial {y}}{f(x+r,y)},\\ 
&C_{1,i}f(x,y)&=C_{3,i}f(x,y)=0\nonumber 
\end{eqnarray}
for any $m,i\in\Z,r\in\Z+\frac{1}{2}, f(x,y)\in\C[x,y]$, where $\lambda,\zeta\in\C^*,\alpha,\beta,\gamma,\sigma\in\C,h\in \C[y]$, and $h_m(y)$ and $g_r(y)$ are defined in (\ref{formula 2.2}) and (\ref{formula 2.1}), respectively.

\begin{defi}\label{def 5.1}
    An $\mathcal{H}$-module $\Omega_{\mathcal{H}}(\lambda,\zeta,\alpha,\beta,\gamma,\sigma,h)$ is said to be of $\sigma$-type if $\sigma\neq0$. In this case, we denote it by  $\Psi(\lambda,\zeta,\sigma,\alpha,\beta,\gamma,h)$.
\end{defi}

\begin{rema}\label{rema 5.1}
    For any $\Psi(\lambda,\zeta,\sigma,\alpha,\beta,\gamma,h)$,  we can recover the actions of $W_{m,i}$ and $B_{r,i}$ through
$\frac{\zeta^i}{\sigma}(\zeta^{-1}L_{m,1}-L_{m,0})f(x,y)=\zeta^i\lambda^m(y-m\alpha)f(x+m,y)=W_{m,i}f(x,y)$
and $ \frac{\zeta^i}{\sigma}(\zeta^{-1}I_{r,1}-I_{r,0})f(x,y)= \zeta^i\gamma\lambda^{r-\frac12}f(x+r,y)=B_{r,i}f(x,y),$
respectively. This means that when $\sigma\neq 0$, $\mathcal{G}$-module $\Omega(\lambda,\zeta,\alpha,\beta,\gamma,\sigma,h)$ and $\mathcal{H}$-module $\Omega_{\mathcal{H}}(\lambda,\zeta,\alpha,\beta,\gamma,\sigma,h)$ have exactly the same submodule structure.
\end{rema}

By Remark \ref{rema 5.1}, we obtain a result of the irreducible condition that is formally consistent with Theorem \ref{theo 3.3}.
\begin{prop}\label{prop 5.1}
The $\mathcal{H}$-module $\Psi(\lambda,\zeta,\sigma,\alpha,\beta,\gamma,h)$ is simple if and only if $(\alpha, \gamma)\neq(0,0)$.
\end{prop}

\begin{theo}\label{th4.4}
    Suppose $\sigma,\sigma'\in\C^*$. Let
$\lambda,\lambda',\zeta,\zeta'\in\mathbb C^*,
\alpha,\alpha',\beta,\beta',\gamma,\gamma'\in\mathbb C$,
$h,h'\in\mathbb C[y]$.
Then
$\Psi(\lambda,\zeta,\sigma,\alpha,\beta,\gamma,h)
\cong
\Psi(\lambda',\zeta',\sigma',\alpha',\beta',\gamma',h')$
as $\mathcal H$-modules if and only if
\begin{equation*}
\aligned
\lambda=\lambda',\zeta=\zeta',\beta=\beta',
\sigma\alpha=\sigma'\alpha',
\sigma\gamma=\sigma'\gamma',
h'(y)=h\!\left(\frac{\sigma'}{\sigma}y\right).
\endaligned
\end{equation*}
\begin{proof}
\rm(1)   Suppose that $\varphi:\Psi(\lambda,\zeta,\sigma,\alpha,\beta,\gamma,h)
\rightarrow
\Psi(\lambda',\zeta',\sigma',\alpha',\beta',\gamma',h')$
 is an isomorphism. 
For any $f\in\C[x,y]$, the identity $\varphi(L_{0}f)=L_{0}\varphi(f)$ implies $\varphi(xf)=x\varphi(f)$. Clearly, $L_{0,2}-2\zeta L_{0,1}+\zeta^2L_{0}$ is the annihilation operator on $\C[x,y]$.
The surjectivity of $\varphi$ implies that the same operator vanishes on $\varphi(\C[x,y])$. Consequently,
\begin{equation*}
\aligned
0
=&\varphi\Big((L_{0,2}-2\zeta L_{0,1}+\zeta^2L_{0})1\Big)
\\=&(L_{0,2}-2\zeta L_{0,1}+\zeta^2L_{0})\varphi(1)
=\Big((\zeta'-\zeta)^2x
+2\zeta'\sigma'(\zeta'-\zeta)y\Big)\varphi(1),
\endaligned
\end{equation*}
which implies $\zeta=\zeta'$.
By the identity $\varphi(L_{0,1}f)=L_{0,1}\varphi(f)$, we have $\zeta\bigl(x\varphi(f)+\sigma\varphi(yf)\bigr)
=\varphi\bigl(\zeta(x+\sigma y)f\bigr)=\zeta(x+\sigma'y)\varphi(f),$
which means $\varphi(yf)=\frac{\sigma'}{\sigma}y\varphi(f)$ and $\varphi(f(x,y))
=f\!\left(x,\frac{\sigma'}{\sigma}y\right)\varphi(1)$.
Since $\varphi$ is surjective, there exists $k\in\mathbb C[x,y]$ such that
$1=k\!\left(x,\frac{\sigma'}{\sigma}y\right)\varphi(1)$,
which means $\varphi(1)\in\mathbb C^*$.

From  $\lambda(\sigma'y-\sigma\alpha)\varphi(1)=\varphi\bigl((\zeta^{-1}L_{1,1}-L_{1,0})1\bigr)
=(\zeta^{-1}L_{1,1}-L_{1,0})\varphi(1)=\lambda'(\sigma'y-\sigma'\alpha')\varphi(1),$
we have $\lambda'=\lambda,\sigma'\alpha'=\sigma\alpha$.
Similarly, we can obtain $\sigma'\gamma'=\sigma\gamma$ by $\varphi\Big((\zeta^{-1}I_{\frac12,1}-I_{\frac12,0})1\Big)=\varphi(\sigma\gamma)$. 
Due to 
\begin{equation*}
\aligned
&\left(x+h\!\left(\frac{\sigma'}{\sigma}y\right)\right)\varphi(1)=\varphi(\lambda^{-1}L_{1,0}1)
=\lambda^{-1}L_{1,0}\varphi(1)
=(x+h'(y))\varphi(1),
\endaligned
\end{equation*}
we have 
$h'(y)=h\!\left(\frac{\sigma'}{\sigma}y\right)$.
By the relations already established,
$g_r'(y)=\frac{\sigma'}{\sigma}
g_r\!\left(\frac{\sigma'}{\sigma}y\right),
\gamma'g_r'(y)=\gamma g_r\!\left(\frac{\sigma'}{\sigma}y\right).$
It is clear that $\beta'=\beta$ from
\begin{equation*}
\aligned
&\left(\beta+\gamma g_{\frac12}
\!\left(\frac{\sigma'}{\sigma}y\right)\right)\varphi(1)=\varphi(I_{\frac12,0}1)
=I_{\frac12,0}\varphi(1)
=\Big(\beta'+\gamma'g_{\frac12}'(y)\Big)\varphi(1).
\endaligned
\end{equation*}

\rm(2)  
Conversely, define
$\varphi(f(x,y))=f\!\left(x,\frac{\sigma'}{\sigma}y\right).$
This is a linear bijection. By (\ref{formula 2.1}) and (\ref{formula 2.2}),
$g_r'(y)=\frac{\sigma'}{\sigma}
g_r\!\left(\frac{\sigma'}{\sigma}y\right),
h_m'(y)=h_m\!\left(\frac{\sigma'}{\sigma}y\right)$.
Moreover,
$\varphi\!\left(\frac{\partial }{\partial y}f\right)
=\frac{\sigma}{\sigma'}
\frac{\partial}{\partial y}\varphi(f)$.
Substitution into the action formulas gives
$\varphi(L_{m,i}f)=L_{m,i}\varphi(f)$,
$\varphi(I_{r,i}f)=I_{r,i}\varphi(f)$.
The central elements act trivially on both modules. Hence $\varphi$ is an $\mathcal H$-isomorphism.
\end{proof}

\end{theo}

\subsection{\texorpdfstring{ Degenerate $\Omega_{\mathcal{H}}(\lambda,\zeta,\alpha,\beta,\gamma,\sigma,h)$ }{1}}
In this subsection, we consider the other case, $\sigma=0$. For simplicity, we denote the $\mathcal{H}$-module by 
$\Omega_{\mathcal{H}}(\lambda,\zeta,\alpha,\beta,\gamma,h)$. We present the $\mathcal{H}$-action as follows.
\begin{equation}\label{act. def 5.1}
\aligned
&L_{m,i} f(x,y)=\lambda^m\zeta^i(x+h_{m}(y))f(x+m,y)+m\lambda^m\zeta^i(y-m\alpha)\frac{\partial }{\partial {y}}{f(x+m,y)},\\
&I_{r,i} f(x,y)=\lambda^{r-\frac{1}{2}}\zeta^i(\beta+\gamma g_{r}(y))f(x+r,y)+r\gamma\lambda^{r-\frac{1}{2}}\zeta^i\frac{\partial }{\partial {y}}{f(x+r,y)},\\ 
&C_{1,i}f(x,y)=C_{3,i}f(x,y)=0
\endaligned
\end{equation}
for any $m,i\in\Z,r\in\Z+\frac{1}{2}, f(x,y)\in\C[x,y]$, where $\lambda,\zeta\in\C^*,\alpha,\beta,\gamma\in\C,h\in \C[y]$, and $h_m(y)$ and $g_r(y)$ are defined in (\ref{formula 2.2}) and (\ref{formula 2.1}) respectively.  Let $g(y)=\frac{h(y)-h(\alpha)}{y-\alpha},\ D=\frac{\partial}{\partial{y}}+g(y),\ s=h(\alpha).$
Then $h(y)=s+(y-\alpha)g(y),g_r(y)=rg(y), h_m(y)=ms+m(y-m\alpha)g(y)$, and (\ref{act. def 5.1}) is transformed into
\begin{equation*}\label{act. def 6.2}
\aligned
&L_{m,i} f(x,y)=\lambda^m\zeta^i(x+ms+m(y-m\alpha)D)f(x+m,y),  \\
&I_{r,i} f(x,y)=\lambda^{r-\frac{1}{2}}\zeta^i(\beta+r\gamma D)f(x+r,y),
C_{1,i}f(x,y)=C_{3,i}f(x,y)=0.
\endaligned
\end{equation*}

\begin{rema}\label{re 6.1}
    For any $e(x)\in\C[x]$, the operators $De(x)-e(x)D$ and $Dy-yD$ on $\C[x,y]$ are the annihilation operator $\operatorname{0}_{\C[x,y]}$ and the identity operator $\operatorname{Id}_{\C[x,y]}$, respectively.
\end{rema}


\begin{rema}\label{rem 6.1}
    By Remark \ref{re 6.1}, it is clear that $xDf=D(xf)\in D\C[x,y]$, 
and
\[
   yD^2(f)= (Dy-{\operatorname{Id}_{\C[x,y]}})(D(f))=D\Big(yD(f)-f\Big) \in D\C[x,y]
\]
for any $f\in\C[x,y]$. 
\end{rema}

\subsubsection{Irreducibility and maximal submodules}

\begin{theo}\label{prop 6.1}
    $\Omega_{\mathcal{H}}(\lambda,\zeta,\alpha,\beta,\gamma,h)$ is simple if and only if ${\rm deg}_y\ (h(y))=1$ and $(\alpha,\gamma)\neq(0,0)$.
\end{theo}
\begin{proof}
\rm(1) 
If $h\in\C$, then $g(y)=0$, $D=\frac{\partial}{\partial{y}}$ and $\C[x]$ is a nonzero proper $\mathcal{H}$-submodule. 

If ${\rm deg}_y\ (h(y))\ge 2$, then ${\rm deg}_y\ (g(y))\ge 1 $ and $ {\rm deg}_y\ (Df(x,y))={\rm deg}_y\ (f(x,y))+{\rm deg}_y\ (g(y)) \ge 1$ for any $0\neq f(x,y)\in\C[x,y]$. Hence, we have $1\notin D\C[x,y]$ and $D(1)=g(y)\neq 0$, which means that $D\C[x,y]$ is a nonzero proper subspace. By Remark \ref{rem 6.1}, $\mathcal{H}D(f(x,y))\in D\C[x,y]$. Hence,
$D\C[x,y]$ is a nonzero proper submodule.

If $(\alpha,\gamma)=(0,0)$, 
obviously, $y\C[x,y]$ is a nonzero proper submodule.
Therefore, simplicity implies  ${\rm deg}_y(h(y))=1$ and $(\alpha,\gamma)\neq (0,0)$.

\rm(2) Suppose that $h(y)=cy+d$ and $(\alpha,\gamma)\neq (0,0)$, where $c\in\C^*$. Then $g(y)=c$ and $D=\frac{\partial}{\partial{y}}+c$. It is clear that $\ker(D)= {0}$. Let $M$ be a nonzero $\mathcal{H}$-submodule of $\Omega_{\mathcal{H}}(\lambda,\zeta,\alpha,\beta,\gamma,h)$. Assume $f(x,y)=\sum_{j=0}^nx^jf_j(y)\in M$, where $f_j(y)\in \C[y]$ and $f_n(y)\neq 0$.
If $\gamma\neq 0$, then the highest-degree coefficient of $\lambda^{\frac{1}{2}-r}I_{r}f(x,y)=(\beta+r\gamma D)f(x+r,y)$
with respect to $r$ is $\gamma D(f_n(y))$. From Lemma \ref{le2.3}, we obtain that
 \begin{equation}\label{formula 5.2}
\aligned
0\neq  D(f_n(y))\in M\cap \C[y].
\endaligned
\end{equation} 
If $\gamma=0$, then $\alpha\neq 0$ and the highest-degree coefficient of $\lambda^{-m}L_mf(x,y)=(x+ms+myD-m^2\alpha D)f(x+m,y)$
with respect to $m$ is $-\alpha D(f_n(y))$. Thus (\ref{formula 5.2}) also holds in this case.

Letting $q(y)\in M\cap \C[y]$ and assuming $\gamma\neq 0$, we have $D(q(y))\in M\cap \C[y]$ 
since $\lambda^{\frac{1}{2}-r}I_{r}q(y)=\beta q(y)+r\gamma D(q(y))\in M\cap \C[y].$
Moreover, independently of whether $\gamma$ vanishes,
We observe the quadratic term in $m$ for 
\begin{equation}\label{formula 5.3}
\aligned
(\lambda^{-m}L_{m}-L_0-ms)q(y)=myD(q(y))-m^2\alpha D(q(y))\in M\cap \C[y],
\endaligned
\end{equation} 
and still obtain $D(q(y))\in M\cap \C[y]$. Comparing the coefficients of $m$ in (\ref{formula 5.3}), we have $yD(q(y))\in M\cap \C[y]$. Then we can obtain that $\frac{\partial}{\partial{y}}q(y)\in M\cap \C[y]$ from $D(q(y))=\frac{\partial}{\partial{y}}q(y)+cq(y)\in M\cap \C[y]$ for any $q(y)\in M\cap \C[y]$. Therefore, we verify $1\in M\cap \C[y]$ and $cy=yD(1)\in M\cap \C[y]$. By induction on $k$ in $yD(y^k)=ky^k+cy^{k+1}\in M\cap \C[y]$, we  have $y^k\in M\cap \C[y]$ for any $k\in\N$ and $\C[y]\subset M\cap \C[y]$. For all $i,k\in\N$, the action of $L_{0}^i$ on $y^k$ implies that $x^iy^k\in M$ and $\C[x,y]= M$. 
\end{proof}

The proof of Theorem \ref{prop 6.1} provides nonzero proper submodules in all reducible cases. Next, we study the maximality of some of these submodules and construct the corresponding simple quotients. The following results do not give a classification of all maximal submodules.

\begin{prop}\label{prop 5.2}
Suppose that $\alpha=\gamma=0$.
\begin{itemize}
    \item[(1)]  If $(\beta,h(0))\ne(0,0)$, then $y\mathbb C[x,y]$ is a maximal submodule.
    \item[(2)] If $\beta=h(0)=0$, then $y\mathbb C[x,y]+x\mathbb C[x]$ is a maximal submodule, and the corresponding quotient is the one-dimensional trivial module.
\end{itemize}
\end{prop}
\begin{proof}
The space $y\mathbb C[x,y]$ is a submodule, and $\C[x,y]/y\C[x,y]\cong \C[x]$ 
with the action $L_{m,i}f(x)
=\lambda^m\zeta^i(x+mh(0))f(x+m),
I_{r,i}f(x)
=\beta\lambda^{r-\frac12}\zeta^if(x+r)$ for every $f(x)\in\mathbb C[x]$.

\rm(1) Suppose that $(\beta,h(0))\ne(0,0)$, and let
$f(x)=\sum_{j=0}^{n}a_jx^j, a_n\ne0$,
belong to a nonzero submodule of this quotient. If $\beta\ne0$, comparison of the highest-degree coefficients in $r$ in
$\lambda^{\frac12-r}I_{r,0}f(x)=\beta f(x+r)$
shows that the submodule contains the nonzero constant $\beta a_n$. If $\beta=0$, then $h(0)\ne0$, and comparison of the highest-degree coefficients in $m$ in
$\lambda^{-m}L_{m,0}f(x)=(x+mh(0))f(x+m)$
gives the nonzero constant $h(0)a_n$. In either case the submodule contains $1$, and repeated application of $L_{0}$ gives all of $\mathbb C[x]$. Hence the quotient is simple.

\rm(2) If $\beta=h(0)=0$, then $x\mathbb C[x]$ is a submodule of the quotient. Its inverse image under the canonical quotient map is $y\mathbb C[x,y]+x\mathbb C[x]$. Since
$\mathbb C[x,y]\big/\bigl(y\mathbb C[x,y]+x\mathbb C[x]\bigr)
 \cong \mathbb C[x]/x\mathbb C[x]\cong \C
$
is the one-dimensional trivial module, this inverse image is a maximal submodule.

\end{proof}

\begin{prop}\label{prop 5.3}
Suppose that ${\deg}_y (h)\ge2$, and let M be a codimension-one subspace of $\mathbb C[y]$ containing $D\mathbb C[y]$.
\begin{itemize}
    \item[(1)]  If $(\beta,s-1)\ne(0,0)$, then $\mathbb C[x]M$ is a maximal submodule.
    \item[(2)] If $\beta=0$ and $s=1$, then $\mathbb C[x]M+x\mathbb C[x,y]$ is a maximal submodule, and the corresponding quotient is the one-dimensional trivial module.
\end{itemize}
\end{prop}
\begin{proof}
Since ${\rm deg}_y (g)={\rm deg}_y (h)-1\ge1$, we have $g\ne 0$ and
${\rm deg}_y(D(q))={\rm deg}_y (q)+{\rm deg}_y (g)$
for every nonzero $q\in\mathbb C[y]$. 
Consequently,
$D\mathbb C[y]\cap
\operatorname{span}_{\mathbb C}\{1,y,\cdots,y^{{\rm deg}_y (g)-1}\}=0$.
Successive elimination of leading terms using the polynomials $D(y^j)$ gives
$$\mathbb C[y]
=
D\mathbb C[y]\oplus
\operatorname{span}_{\mathbb C}\{1,y,\cdots,y^{{\rm deg}_y (g)-1}\}.$$
Hence
$\dim\bigl(\mathbb C[y]/D\mathbb C[y]\bigr)={\rm deg}_y (g)\ge1$. The inverse image of any codimension-one subspace of this quotient is a subspace M of the required form.
For $q\in M$, we have
$D(q)\in D\mathbb C[y]\subseteq M,
yD(q)=D(yq)-q\in M$.
Therefore $\mathbb C[x]M$ is a submodule.
Choose $q\notin M$. Since $\mathbb C[y]=M\oplus\mathbb Cq$, every element of the quotient has a unique expression
$f(x)q+\mathbb C[x]M$.
Moreover,
$D(q)\equiv0, yD(q)=D(yq)-q\equiv-q\pmod M$.

 If $(\beta,s-1)\ne(0,0)$, the same coefficient comparison as in Proposition \ref{prop 5.2} proves that this quotient is simple.
If $\beta=0$ and $s=1$, the subspace represented by $x\mathbb C[x]q$ is a submodule of the quotient, with one-dimensional trivial quotient. Its inverse image is
$\mathbb C[x]M+x\mathbb C[x]q
=\mathbb C[x]M+x\mathbb C[x,y]$.
Thus this inverse image is maximal.
Finally, $M$ contains the nonzero subspace $D\mathbb C[y]$, while each quotient constructed above is nonzero. Hence the stated submodules are nonzero and proper.

\end{proof}

\begin{prop}
    If $h\in\mathbb C$ and $(\alpha,\gamma)\ne(0,0)$, then
$\Omega_{\mathcal{H}}(\lambda,\zeta,\alpha,\beta,\gamma,h)$
has no maximal submodules.
\end{prop}
\begin{proof}

Set
$\Omega=\Omega_{\mathcal H}(\lambda,\zeta,\alpha,\beta,\gamma,h),
\Omega_n=\bigoplus_{j=0}^{n}\mathbb C[x]y^j
,n\ge0$.
Then each $\Omega_n$ is an $\mathcal H$-submodule and
\begin{equation}\label{formula 5.4.1}
\aligned
\Omega_0\subsetneq\Omega_1\subsetneq\cdots,
\Omega=\bigcup_{n\ge0}\Omega_n.
\endaligned
\end{equation}

Let N be a proper submodule of $\Omega$. We claim that the y-degrees of the elements of N are bounded. Otherwise, N would contain elements of arbitrarily large y-degree. Since $(\alpha,\gamma)\neq(0,0)$, comparing coefficients in m or r in the above actions and repeatedly using the nonzero degree-lowering term
$-\alpha\frac{\partial}{\partial{y}}\ \text{or}\ \gamma\frac{\partial}{\partial{y}}$,
the same degree-reduction argument as above yields
$\Omega_n\subseteq N\ \text{for every }n\ge0$.
Thus (\ref{formula 5.4.1}) gives $N=\Omega$, a contradiction. Hence
$N\subseteq\Omega_n$
for some $n\ge0$.

Consequently,
$N\subseteq\Omega_n
\subsetneq\Omega_{n+1}
\subsetneq\Omega$,
so $N$ cannot be maximal. Since $N$ was arbitrary, $\Omega$ has no maximal submodules.

\end{proof}

\subsubsection{Isomorphism}
In this subsection, we will determine the isomorphism classes of the $\mathcal{H}$-modules
$\Omega_{\mathcal{H}}(\lambda,\zeta,\alpha,\beta,\gamma,h)$, including both irreducible and reducible types. 

\begin{lemm}\label{preth10}
    Let $\Omega=\Omega_{\mathcal{H}}(\lambda,\zeta,\alpha,\beta,\gamma,h),
\Omega'=\Omega_{\mathcal{H}}(\lambda',\zeta',\alpha',\beta',\gamma',h')$,
and set $g'(y)=\frac{h'(y)-h'(\alpha')}{y- \alpha'},  D'=\frac{\partial}{\partial{y}}+g'(y), s'=h'(\alpha'). $ If $\varphi:\Omega\longrightarrow\Omega'$ is an isomorphism of $\mathcal{H}$-modules, then the following statements hold:
\begin{itemize}
\item[{\rm (1)}] $(\lambda,\zeta,\beta)=(\lambda',\zeta',\beta')$ and $\alpha'=0\iff \alpha=0,\ \gamma'=0\iff\gamma=0$.
\item[{\rm (2)}]For any $f\in\C[x,y]$, $\varphi(xf)=x\varphi(f)$ and $\varphi(X^2D^2(f))=(X')^2(D')^2\varphi(f)$ for $X=\alpha,\gamma$.
\item[{\rm (3)}]If $(\alpha,\gamma)\neq(0,0)$, then $\varphi(\C[y])=\C[y]$, and there exists $t\in\C^*$ such that $(\alpha',\gamma')=t(\alpha,\gamma)$, $\varphi((s+yD)q)=(s'+yD')\varphi(q)$, $\varphi(D(q))=tD'(\varphi(q))$ for any $q\in\C[y]$.
\end{itemize}
\end{lemm}
\begin{proof}
Denote 
$f_m=f(x+m,y),\varphi_m(f)=(\varphi(f))(x+m,y)$ for any $f=f(x,y)\in\C[x,y],m\in\frac{1}{2}\Z$. For any $m,i\in \Z,f\in\C[x,y]$, from
$\varphi(L_{0,i}f)=L_{0,i}\varphi(f)$ and $\varphi(L_m1)=L_m\varphi(1)$, i.e. $$\lambda^m\left(x\varphi(1)+m\varphi(s+yg(y)) - m^2\alpha\varphi(g(y)) \right) =(\lambda')^m \left(x+m(s'+yD')- m^2\alpha'D'\right) \varphi_m(1),$$
 we have $\varphi(xf)=x\varphi(f)$, $\zeta'=\zeta$ and $\lambda' =\lambda$ 
 . For any $r\in\Z+\frac{1}{2}$, we obtain $\beta'=\beta$ by using Remark \ref{re 6.1}  and observing the constant term of $r$ in $$\varphi(I_r1)=\varphi\left(\lambda^{r-\frac{1}{2}}(\beta+r\gamma D)1\right)
=\lambda^{r-\frac{1}{2}}\beta\varphi(1)+r\gamma \lambda^{r-\frac{1}{2}}\varphi(g(y))
=I_{r}\varphi(1)
=\left(\lambda^{r-\frac{1}{2}}(\beta^\prime+r\gamma^\prime D')\right)\varphi_r(1).$$
For any $m\in\Z,r\in\Z+\frac{1}{2}$,$f=f(x,y)\in\C[x,y]$, 
we have
\begin{equation}\label{formula 5.4}
\aligned
\varphi&\left(\left(x^2+mx+m^2\left((s+yD)-(s+yD)^2-2\alpha xD\right)+m^4\alpha^2D^2\right)f\right)=\varphi(L_mL_{-m}f)=L_mL_{-m}\varphi(f)\\
=&\left(x^2+mx+m^2[(s'+yD')-(s'+yD')^2-2\alpha' xD']+m^4(\alpha')^2(D')^2\right)\varphi(f)\ \ \\
\endaligned
\end{equation}
and
\begin{equation}\label{formula 5.5}
\aligned
\varphi\Big((\beta^2-r^2\gamma^2D^2)f\Big)=\lambda\varphi(I_rI_{-r}f)=\lambda I_rI_{-r}\varphi(f)
=\Big(\beta^2-r^2(\gamma')^2(D')^2\Big)\varphi(f).
\endaligned
\end{equation}
Comparing the coefficients of $m^4$ and $r^2$ in (\ref{formula 5.4}) and (\ref{formula 5.5}), respectively, we can verify
 \begin{equation}\label{formula 5.6}
\aligned
\varphi(\alpha^2D^2(f))=(\alpha')^2(D')^2\varphi(f),
\endaligned
\end{equation}
 \begin{equation*}\label{formula 6.2.4.2}
\aligned
\varphi(\gamma^2D^2(f))=(\gamma')^2(D')^2\varphi(f).
\endaligned
\end{equation*}
Since neither $D^2$ nor $(D')^2$ is the zero operator, it follows that
\begin{equation}\label{formula 5.12}
\aligned
\alpha'=0\iff \alpha=0,\ \ \ \gamma'=0\iff\gamma=0.
\endaligned
\end{equation}
This proves (1) and (2).
  
We now prove (3). If $\alpha\ne 0$, then $Dy-yD=\operatorname{Id}_{\C[x,y]}$ gives $D^2y-yD^2=2D$. Thus, for any $f\in\C[x,y]$, $\alpha^2(D^2L_m-L_mD^2)f=2m\alpha^2\lambda^mD^2(f_m )$
for any $m\in\Z$. 
Applying $\varphi$ and using (\ref{formula 5.6}), 
we can obtain that $\varphi(\alpha^2D^2(f_m))=(\alpha')^2(D')^2\varphi_m(f).$
This identity is immediate for $m=0$. Using (\ref{formula 5.6}) again, we have
$(\alpha')^2(D')^2\varphi(f_m)=(\alpha')^2(D')^2\varphi_m(f).$
If $g'\ne0$, then $D'$ is injective, since the leading term of $g'p$ cannot be cancelled by $\frac{\partial}{\partial{y}}p$ for any nonzero polynomial $p$. As $\alpha'\ne0$, it follows that $\varphi(f_m)=\varphi_m(f)$, that is,
\begin{equation}\label{formula 5.10}
\aligned
 \varphi(u_m)=\varphi_m(u)  
\endaligned
\end{equation}
for any $u=u(x,y)\in\C[x,y]$.
If $g'=0$, then $(\alpha')^2(D')^2=(\alpha')^2\Big(\frac{\partial}{\partial{y}}\Big)^2$ is surjective. By (\ref{formula 5.6}), $\alpha^2D^2$ is also surjective. For any $u(x,y)\in\C[x,y]$, there exists $f(x,y)\in\C[x,y]$ such that $\alpha^2D^2(f)=u$. By Remark \ref{re 6.1}, we still obtain (\ref{formula 5.10}) from
\begin{equation*}
\aligned
 \varphi(u_m)=\varphi\big(\alpha^2D^2(f_m)\big)=(\alpha')^2(D')^2\varphi_m(f)=\Big(   (\alpha')^2(D')^2\varphi(f) \Big)_m=\varphi_m(u).
\endaligned
\end{equation*}
If $\alpha=0$ and $\gamma\neq 0$, replacing $\alpha^2D^2$ by $\gamma^2D^2$ in the preceding argument gives (\ref{formula 5.10}) as well. 
For any $q(y)\in\mathbb C[y]$, taking $m=1$ in (\ref{formula 5.10}) yields $\varphi(q)(x+1,y)=\varphi(q)(x,y)$. Hence $\varphi(q)$ is independent of $x$. Applying the same argument to $\varphi^{-1}$, we obtain $\varphi(\mathbb C[y])=\mathbb C[y]$.

For any $q(y)\in\mathbb C[y]$,
comparing coefficients in the identities
$\varphi(L_mq)=L_m(\varphi(q))$ and $\varphi(I_rq)=I_r(\varphi(q))$, we obtain
\begin{equation}\label{formula 5.11}
\aligned
 \varphi((s+yD)q)&=(s'+yD')\varphi(q),\alpha\varphi(D(q))=\alpha'D'(\varphi(q)), \gamma\varphi(D(q))=\gamma'D'(\varphi(q)).
\endaligned
\end{equation}
If $\alpha\neq0$, take $t=\frac{\alpha'}{\alpha}$; otherwise take $t=\frac{\gamma'}{\gamma}$.
By (\ref{formula 5.12}) and (\ref{formula 5.11}), and since $D'\varphi\neq0$, we have
\begin{equation}\label{formula 5.13}
\aligned
 (\alpha',\gamma')=t(\alpha,\gamma),\ \ \ \varphi(D(q))=tD'(\varphi(q)).
\endaligned
\end{equation}
The proof is complete.
\end{proof}

We recall the following definitions.
\begin{defi}
    An operator $T\in\operatorname{End}_{\mathbb C}(V)$ is locally finite on the complex vector space $V$ if for every $v\in V$,
$\operatorname{dim}\ \operatorname{span}
\{T^nv\mid n\in\N\}<+\infty$.
\end{defi}
\begin{defi}
    An operator $T\in\operatorname{End}_{\mathbb C}(V)$ is locally nilpotent on the complex vector space $V$ if for every $v\in V$, there exists $n
    \in\N$ such that
$
T^nv=0 $.
\end{defi}

\begin{theo}\label{theo 5.2}
Let
$\Omega=\Omega_{\mathcal{H}}(\lambda,\zeta,\alpha,\beta,\gamma,h)$ and
$\Omega'=\Omega_{\mathcal{H}}(\lambda',\zeta',\alpha',\beta',\gamma',h')$,
where $\lambda,\zeta,\lambda',\zeta'\in\mathbb C^*$, $\alpha,\beta,\gamma,\alpha',\beta',\gamma'\in\mathbb C$, and $h,h'\in\mathbb C[y]$. Then $\Omega\cong\Omega'$ as $\mathcal{H}$-modules if and only if
$\lambda=\lambda', \zeta=\zeta', \beta=\beta'$,
and one of the following mutually exclusive conditions holds:
\begin{itemize}
\item[{\rm (1)}] $\alpha=\gamma=\alpha'=\gamma'=0$, and $h=h'\in\mathbb C$.
\item[{\rm (2)}] $\alpha=\gamma=\alpha'=\gamma'=0$, and $h,h'$ are nonconstant polynomials with $\deg_y (h)=\deg_y (h')$.
\item[{\rm (3)}] $(\alpha,\gamma)\ne(0,0), h=h'\in\mathbb C$, and there exists $t\in\mathbb C^*$ such that $(\alpha',\gamma')=t(\alpha,\gamma)$.
\item[{\rm (4)}] $(\alpha,\gamma)\ne(0,0), h(y)=cy+d$ and $h'(y)=c'y+d'$ for some $c,c'\in\mathbb C^*$ and $d,d'\in\mathbb C$, and
$\alpha c=\alpha'c',\gamma c=\gamma'c'$.
\item[{\rm (5)}] $(\alpha,\gamma)\ne(0,0), \deg_y (h)\ge2$, and there exists $t\in\mathbb C^*$ such that
$(\alpha',\gamma')=t(\alpha,\gamma),h'(y)=h(y/t)$.
\end{itemize}
\end{theo}
\begin{proof}
We retain the notation of Lemma \ref{preth10}.
Part ($\mathrm{I}$) We first prove necessity. Suppose that $\varphi: \Omega\longrightarrow\Omega'$ is a module isomorphism. By Lemma \ref{preth10} (1) and (2), $(\lambda,\zeta,\beta)=(\lambda',\zeta',\beta'),$  $\alpha'=0\iff \alpha=0,\ \gamma'=0\iff\gamma=0$ and for any $f\in\C[x,y]$, $\varphi(xf)=x\varphi(f)$, $\varphi(X^2D^2(f))=(X')^2(D')^2\varphi(f)$ for $X=\alpha,\gamma$.

\rm(1) Suppose that $(\alpha,\gamma)=(0,0)$. Then $(\alpha',\gamma')=(0,0)$ and $s+yD=y\frac{\partial}{\partial{y}}+h(y)$. If $h=d\in\C$, we have $(s+yD)x^ky^l=(d+l)x^ky^l$ for any $k,l\in\N$, which means that $(s+yD)- (s+yD)^2$ is locally finite. If $h$ is nonconstant, then, for any $0\ne f\in\C[x,y]$, we have $${\rm deg}_y  \left(\left( (s+yD)- (s+yD)^2 \right)f\right)={\rm deg}_y f+2{\rm deg}_y h>{\rm deg}_y f,$$
which implies that $(s+yD)- (s+yD)^2$ is not locally finite. Comparing the coefficients of $m^2$ in (\ref{formula 5.4}), we see that these two cases cannot be isomorphic. 

Suppose now that $h=d,h'=d'$ are constant. Then
 \begin{equation}\label{formula 5.7}
\aligned
\Omega=\bigoplus_{k\ge0}\mathbb C[x]y^k,\ \Omega'=\bigoplus_{j\ge0}\mathbb C[x]y^j, 
\endaligned
\end{equation}
  where every summand is a submodule. A homomorphism from the $k$-th summand to the $j$-th summand is determined by $y^k\longmapsto p(x)y^j$, where 
  \begin{equation*}\label{formula 6.2.6}
\aligned
 (x+m(d+k))p(x) &=(x+m(d'+j))p(x+m),\ \beta p(x)=\beta p(x+r). 
\endaligned
\end{equation*}
Comparing degrees and coefficients in $m$, we see that a nonzero homomorphism is possible only when 
$d+k=d'+j, p(x)\in\mathbb C^*$, or when $\beta=0,d+k=1,d'+j=0,  p(x)\in\mathbb C^*x$. In the latter case, there is no nonzero homomorphism in the reverse direction.  
  Fix $k$, restrict $\varphi^{-1}\varphi=\operatorname{Id}$ to $\C[x]y^k$, and project back onto this summand. Since $\varphi(y^k)$ has only finitely many nonzero components, at least one composition through a summand $\C[x]y^j$ of $\Omega'$ is nonzero. The preceding classification implies $d+k=d'+j$. Applying the same argument to $\varphi\varphi^{-1}$, we obtain
 $\{d+k\mid k\in\mathbb N\} =\{d'+j\mid j\in\mathbb N\}$. In particular, $d=d'+j$ and $d'=d+k$ for some $j,k\in\N$, whence $j+k=0$ and $d'=d$.

\rm(2) Suppose that $h$ and $h'$ are nonconstant. For any $0\neq q(y)\in\mathbb C[y]$, write $\varphi(q(y))=\sum_{j=0}^{n}x^jq_j(y)$, where $q_n(y)\ne0$. In the identity $\varphi\left(\left(x+m(s+yD)\right)q\right) =\left(x+m(s'+yD')\right)\varphi_m(q),$ the left-hand side has degree at most $1$ in m. If $n\ge1$, the coefficient of $m^{n+1}$ on the right-hand side is $(s'+yD')q_n(y)\ne0$. Thus $n=0$. Applying the same argument to $\varphi^{-1}$, we obtain
 \begin{equation}\label{formula 5.8}
\aligned
 \varphi(\mathbb C[y])=\mathbb C[y].
\endaligned
\end{equation}
 Comparing the coefficients of $m$, we have 
 \begin{equation}\label{formula 5.9}
\aligned
 \varphi((s+yD)q)=(s'+yD')\varphi(q). 
\endaligned
\end{equation}
Since $\deg_y((s+yD)q)=\deg_y(q)+\deg_y (h)$, successive elimination of leading terms shows that the images of the monomials of degree less than ${\rm deg}_y (h)$ form a basis of $\mathbb C[y] / (s+yD)\mathbb C[y]$. Hence
   $dim\ \mathbb C[y] / (s+yD)\mathbb C[y]=\deg_y(h)$.
By (\ref{formula 5.8}) and (\ref{formula 5.9}), the corresponding quotient spaces are isomorphic.
Therefore $\deg_y (h')=\deg_y (h)$.

For the remainder of Part ($\mathrm{I}$), assume that $(\alpha,\gamma)\ne(0,0)$. By Lemma \ref{preth10} (3), (\ref{formula 5.8}) still holds. Moreover, (\ref{formula 5.11}) holds, and there exists $t\in\C^*$ satisfying (\ref{formula 5.13}).

\rm(3) 
If $h=d\in\mathbb C$, then $D=\frac{\partial}{\partial{y}}$ is locally nilpotent on $\mathbb C[y]$. By (\ref{formula 5.8}) and (\ref{formula 5.13}), $D'$ is also locally nilpotent. Thus $g'=0$, and $h'=d'\in\mathbb C$. Moreover, (\ref{formula 5.13}) gives $\varphi(\ker (D))=\ker (D')=\mathbb C1$. The operators $s+yD$ and $s'+yD'$ act on these kernels by the scalars $d$ and $d'$, respectively. Hence, (\ref{formula 5.11}) gives $h=d=d'=h'$.

\rm(4) If $h(y)=cy+d$ for $c\ne0$, then $D=\frac{\partial}{\partial{y}}+c$. Taking $q=1$ in (\ref{formula 5.13}), we obtain $tD'(\varphi(1))=c\varphi(1)$. 
Comparing degrees in $y$, we see that $g'$ must be a nonzero constant. Then write $g'=c', h'(y)=c'y+d'$.  
Comparison of the leading coefficients gives $c=tc'$. 
Consequently, $\alpha c=\alpha'c',\gamma c=\gamma'c'$.

\rm(5)  Suppose that $\deg_y (h)\ge2$. By the preceding two cases and the symmetry of isomorphism, we have $\deg_y (h')\ge2$. Thus $D,D'$ are injective, and both quotient spaces $\mathbb C[y]/D\mathbb C[y]$,  $\mathbb C[y]/D'\mathbb C[y]$ are nonzero. Since $yD=Dy-\operatorname{Id}_{\C[x,y]}$ (${\rm resp.}\ yD'=D'y-\operatorname{Id}_{\C[x,y]}$), the operator $s+yD$ (${\rm resp.}\ s'+yD'$) induces multiplication by $s-1$ (${\rm resp.}\ s'-1$) on $\mathbb C[y]/D\mathbb C[y]$ (${\rm resp.}\ \mathbb C[y]/D'\mathbb C[y]$). By (\ref{formula 5.8}), $(\ref{formula 5.11})$ and $(\ref{formula 5.13})$, $\varphi$  induces an isomorphism between these quotient spaces. Hence, $s=s'$ and consequently 
$\varphi(yD(q))=yD'(\varphi(q))$. Together with (\ref{formula 5.13}), this gives $D'\Big(\varphi(yq)\Big)=D'\Big( t^{-1}y\varphi(q)   \Big)$.
 Since $D'$ is injective, $\varphi(yq)=t^{-1}y\varphi(q)$. Therefore  $\varphi(q(y))=q(y/t)\varphi(1)$.
The surjectivity in (\ref{formula 5.8}) implies that $\varphi(1)$ is a nonzero constant. 
Taking $q=1$ in (\ref{formula 5.13}),
we obtain
$g\left(\frac{y}{t}\right)=tg'(y)$. Then, $h'(y)=s'+(y-\alpha')g'(y)=s+\left(\frac{y}{t}-\alpha\right)g\left(\frac{y}{t}\right)=h\left(\frac{y}{t}\right)$.

Part ($\mathrm{II}$) We next prove sufficiency. Throughout this part, assume that $(\lambda, \zeta,\beta)=(\lambda',\zeta',\beta')$.

\rm(1) If condition (1) holds, the actions on the two modules coincide, so the identity map is an isomorphism.

\rm(2) If condition (2) holds, set $n=\deg_y (h)=\deg_y (h')\ge1$. Since $\deg_y((s+yD)^ky^j)=kn+j$ for $k\in\mathbb N,\ 0\le   j<n$, with nonzero leading coefficient, the set 
$\{(s+yD)^ky^j\mid k\in\mathbb N,\ 0\le j<n\}$ is a basis of $\mathbb C[y]$. The corresponding primed expressions also form a basis. Hence, the assignment $\varphi\bigl(x^l(s+yD)^ky^j\bigr) =x^l(s'+yD')^ky^j $ for $(l,k\in\mathbb N,\ 0\le j<n)$ defines a linear bijection. This map satisfies (\ref{formula 5.9}). The action formulas therefore show that it is a module isomorphism.

\rm(3) and \rm(5) If condition (3) or (5) holds, define $\varphi(f(x,y))=f(x,y/t)$. Since $h'(y)=h(y/t)$ and $\alpha'=t\alpha$, we have $s'=s, g'(y)=t^{-1}g(y/t)$. Thus $\varphi(Df)=tD'\varphi(f),\varphi((s+yD)f)=(s'+yD')\varphi(f)$. Together with $\gamma'=t\gamma$, these identities show directly that $\varphi$ is a module isomorphism.

\rm(4) 
Suppose that condition (4) holds. Define 
\begin{equation}\label{formula 5.15}
\aligned
  \varphi(f(x,y)) = \sum_{j\ge0} \binom{s'-s}{j}(c')^{- j} \left(       \frac{\partial}{\partial{y}}\right)^j\left(f\left(x,\frac{c'}c\,y\right)\right).
\endaligned
\end{equation}
For every polynomial $f$, only finitely many terms in (\ref{formula 5.15}) are nonzero. The binomial identity shows that its inverse is $\varphi^{-1}(f(x,y)) = \sum_{j\ge0} \binom{s-s'}{j}c^{- j} \left(       \frac{\partial}{\partial{y}}\right)^j\left(f\left(x,\frac c{c'}\,y\right)\right).$
Using $\left(       \frac{\partial}{\partial{y}}\right)^j(yq) =y\left(       \frac{\partial}{\partial{y}}\right)^jq+j\left(       \frac{\partial}{\partial{y}}\right)^{j-1}q$ for $j\ge 1$, a direct calculation gives $\varphi(Df)=\frac c{c'}D'\big(\varphi(f)\big),\varphi((s+yD)f)=(s'+yD')\varphi(f)$.
Since $\alpha'=\frac c{c'}\alpha,\gamma'=\frac c{c'}\gamma$, we obtain $\varphi(L_{m,i}f)=L_{m,i}\varphi(f),  \varphi(I_{r,i}f)=I_{r,i}\varphi(f)$. The central elements act trivially on both modules. Hence $\varphi$ is an isomorphism of $\mathcal{H}$-modules, completing the proof.

\end{proof}

We have the following conclusion immediately from Theorem \ref{theo 5.2}.

\begin{coro}\label{coro 6.5}

Every $\Omega_{\mathcal{H}}(\lambda,\zeta,\alpha,\beta,\gamma,h)$ satisfying the irreducibility conditions of Theorem \ref{prop 6.1} is isomorphic to exactly one of the following modules:
\begin{equation*}\label{formula 6.5.1}
\aligned
\Omega_{\mathcal{H}}(\lambda,\zeta,a,\beta,b):=\Omega_{\mathcal{H}}(\lambda,\zeta,a,\beta,b,y),
\lambda,\zeta\in\mathbb C^*,
\beta\in\mathbb C,
(a,b)\in\mathbb C^2\setminus\{(0,0)\}.
\endaligned
\end{equation*}

More precisely, if $h(y)=cy+d$ with $c\ne0$, then
$$\Omega_{\mathcal{H}}(\lambda,\zeta,\alpha,\beta,\gamma,cy+d)
\cong
\Omega_{\mathcal{H}}(\lambda,\zeta,\alpha c,\beta,\gamma c,y)=\Omega_{\mathcal{H}}(\lambda,\zeta,\alpha c,\beta,\gamma c).$$
\end{coro}

\section{\texorpdfstring{ Tensor product $\mathcal{H}$-modules $\left(\bigotimes_{i=1}^{m}\Omega_{\mathcal{H}}(\lambda_i,\zeta_i, \alpha_i,\beta_i,\gamma_i)\right)\otimes V$}{Tensor product modules}}\label{section 6}


Fix any $m\in\Z_+$, $\{\lambda_i\}_{i=1}^m,\{\zeta_i\}_{i=1}^m\subset\C^*$, $\{\beta_i\}_{i=1}^m\subset\C$, $\{(\alpha_i, \gamma_i)\}_{i=1}^m\subset\C^2\backslash\{(0,0)\}$. Each $\Omega_{\mathcal{H}}(\lambda_i, \zeta_i,\alpha_i,\beta_i,\gamma_i)$ denotes the irreducible degenerate $\mathcal{H}$-module normalized as in Corollary \ref{coro 6.5} (so that $h_i(y_i)=y_i$) defined on $\C[x_i,y_i]$. By Corollary \ref{coro 6.5}, $\Omega_{\mathcal{H}}(\lambda_i, \zeta_i,\alpha_i,\beta_i,\gamma_i)$ is simple as an $\mathcal{H}$-module defined by
 \begin{equation}\label{act. def 7.1}
\aligned
&L_{k,l} f(x_i,y_i)=\lambda_i^k\zeta_i^l(x_i+ky_i-k(k-1)\alpha_i)f(x_i+k,y_i)+k\lambda_i^k\zeta_i^l(y_i-k\alpha_i)\frac{\partial }{\partial {y_i}}{f(x_i+k,y_i)},\\
&I_{r,l} f(x_i,y_i)=\lambda_i^{r-\frac{1}{2}}\zeta_i^l(\beta_i+r\gamma_i)f(x_i+r,y_i)+r\gamma_i\lambda_i^{r-\frac{1}{2}}\zeta_i^l\frac{\partial }{\partial {y_i}}{f(x_i+r,y_i)},\\ 
&C_{1,l}f(x_i,y_i)=C_{3,l}f(x_i,y_i)=0
\endaligned
\end{equation}
 for every $f(x_i,y_i)\in\C[x_i,y_i],k,l\in\Z,r\in\Z+\frac{1}{2},i=1,2,\cdots,m$. 

 On the other hand, we have the following definition of the restricted $\mathcal{H}$-modules.
 \begin{defi}
  An $\mathcal{H}$-module $V$ is called a restricted module if for every $v\in V$, there exists $K\in\N$ such that $\mathcal{H}_nv=0$ for all $n\in\frac{1}{2}Z$ with $n>K$, where $\mathcal{H}_n=\mathcal{H}\cap \mathcal{G}_n$ and $\mathcal{G}_n$ is defined in (\ref{grading 1.1}).
\end{defi}

 We consider a simple restricted $\mathcal{H}$-module $V$. 
Let us denote 
\[
    \mathrm{T}=\mathrm{T}(\boldsymbol{\lambda},\boldsymbol{\zeta},\boldsymbol{\alpha},\boldsymbol{\beta},\boldsymbol{\gamma})=\left(\bigotimes_{i=1}^{m}\Omega_{\mathcal{H}}(\lambda_i,\zeta_i, \alpha_i,\beta_i,\gamma_i)\right)\otimes V\cong\C[\boldsymbol{x},\boldsymbol{y}]\otimes V,
\] 
\[
    \boldsymbol{0}=(0,0,\cdots,0)\in\Z^m,\boldsymbol{e}_i=(0,\cdots,0,1,0,\cdots,0)\in\C^m\ \ \ \ \ (1\ {\rm{at\  the}}\      i{\rm{-th\  component}}), 
\]
where 
\begin{equation*}\label{formula 5.0.1}
\aligned
\boldsymbol{\lambda}=&(\lambda_1,\lambda_2,\cdots,\lambda_m),\boldsymbol{\zeta}=(\zeta_1,\zeta_2,\cdots,\zeta_m),
\boldsymbol{\alpha}=(\alpha_1,\alpha_2,\cdots,\alpha_m),
\boldsymbol{\beta}=(\beta_1,\beta_2,\cdots,\beta_m),\\
\boldsymbol{\gamma}=&(\gamma_1,\gamma_2,\cdots,\gamma_m),
\boldsymbol{x}=(x_1,\cdots,x_m),
\boldsymbol{y}=(y_1,\cdots,y_m),
\C[\boldsymbol{x},\boldsymbol{y}]=\C[x_1,\cdots,x_m;y_1,\cdots,y_m].
\endaligned
\end{equation*}
Define the directional derivative operator as $\boldsymbol{d}\cdot\nabla=\sum_{i=1}^md_i\frac{\partial}{\partial{x_i}}$
for any $\boldsymbol{d}=(d_1,d_2,\cdots,d_m)\in\C^m$ and define $(\boldsymbol{d}\cdot\nabla) (1^{\otimes m}\otimes v)=0$ for any $v\in V$. Then  $\boldsymbol{e}_i\cdot\nabla=\frac{\partial}{\partial{x_i}}$ and $\frac{\partial}{\partial{x_i}}(1^{\otimes m}\otimes v)=0$. Define a lexicographical order ``$\prec$" on $\mathbb{Z}^n$ as follows. For any 
$(p_1, \cdots, p_n), (q_1, \cdots, q_n) \in \mathbb{Z}^n$, we say 
$(p_1, \cdots, p_n) \prec (q_1, \cdots, q_n)$ if there exists $1\le k \le n$ such that 
$p_k < q_k$ and $p_i = q_i$ for all $1 \le i < k$.

Each nonzero element $f\in\mathrm{T}$ can be uniquely written in the form
\begin{equation}\label{formula 5.17}
\aligned
f=&f(\boldsymbol{x},\boldsymbol{y})=\sum_{(\boldsymbol{s},\boldsymbol{t})\in\N^{2m}}\boldsymbol{x}^{\boldsymbol{s}}\boldsymbol{y}^{\boldsymbol{t}}\otimes v_{\boldsymbol{s},\boldsymbol{t}},
\endaligned 
\end{equation}
where $v_{\boldsymbol{s},\boldsymbol{t}}\in V$, all but finitely many $v_{\boldsymbol{s},\boldsymbol{t}}=0$, $\boldsymbol{s}=(s_1,s_2,\cdots,s_m)$, $\boldsymbol{t}=(t_1,t_2,\cdots,t_m)$. Define the degree of $f$ in $\boldsymbol{x}$ (resp. $\boldsymbol{y}$) by ${\rm deg_{\boldsymbol{x}}}(f)=\boldsymbol{s}$ (resp. ${\rm deg_{\boldsymbol{y}}}(f)=\boldsymbol{t}$), where $\boldsymbol{x}^{\boldsymbol{s}}$ (resp. $\boldsymbol{y}^{\boldsymbol{t}}$) is the maximal element according to the partial order that appears in the sum. We also define 
$s^i={\rm deg}_{x_i}(f):={\rm max}\{s_i\mid v_{\boldsymbol{s},\boldsymbol{t}}\ne 0 \}$ and $t^i={\rm deg}_{y_i}(f):={\rm max}\{t_i\mid v_{\boldsymbol{s},\boldsymbol{t}}\ne 0 \}$. For convenience, we still refer to the elements in $\mathrm{T}$ as polynomials, i.e. \[
    f\boldsymbol{x}^{\boldsymbol{s^\prime}}\boldsymbol{y}^{\boldsymbol{t^\prime}}=\boldsymbol{x}^{\boldsymbol{s^\prime}}\boldsymbol{y}^{\boldsymbol{t^\prime}}f=\sum_{(\boldsymbol{s},\boldsymbol{t})\in\N^{2m}}\boldsymbol{x}^{\boldsymbol{s}+\boldsymbol{s^\prime}}\boldsymbol{y}^{\boldsymbol{t}+\boldsymbol{t^\prime}}\otimes v_{\boldsymbol{s},\boldsymbol{t}},\ \ \ \ \ \Big((\boldsymbol{s^{\prime}},\boldsymbol{t^{\prime}})\in\N^{2m}\Big).
\]
We also denote $\boldsymbol{x}_i=\boldsymbol{x}^{\boldsymbol{e}_i}=1\otimes\cdots \otimes x_i\otimes\cdots\otimes 1$ and $\boldsymbol{y}_i=\boldsymbol{y}^{\boldsymbol{e}_i}=1\otimes\cdots \otimes y_i\otimes\cdots \otimes1$ for every $i=1,\cdots,m$. Then for any $p(x_i)\in\C[x_i],q(y_i)\in\C[y_i],r(x_i,y_i)\in\C[x_i,y_i]$, we have $p(\boldsymbol{x}_i)=1\otimes\dots \otimes p(x_i)\otimes\cdots\otimes 1,q(\boldsymbol{y}_i)=1\otimes\cdots \otimes q(y_i)\otimes\cdots \otimes1,r(\boldsymbol{x}_i,\boldsymbol{y}_i)=1\otimes\cdots \otimes r(x_i,y_i)\otimes\cdots \otimes1$.

\begin{rema}\label{rema 7.1}
    For any polynomial $f(\boldsymbol{x},\boldsymbol{y})$ in $\mathrm{T}$ of the form (\ref{formula 5.17}) and any $\boldsymbol{d}=(d_1,d_2,\cdots,d_m)\in\C^m$,  we have 
\begin{equation*}\label{formula 7.0.1}
\aligned
f(\boldsymbol{x}+\boldsymbol{d},\boldsymbol{y})=&\sum_{j\ge0}\frac{1}{j!}\left(\boldsymbol{d}\cdot \nabla\right)^jf(\boldsymbol{x},\boldsymbol{y})=e^{\boldsymbol{d}\cdot \nabla}f(\boldsymbol{x},\boldsymbol{y}).
\endaligned 
\end{equation*}
    Specifically, for any $k\in \C$, we have $f(\boldsymbol{x}+k\boldsymbol{e}_i,\boldsymbol{y})=\sum_{j=0}^{s^i}\frac{k^j}{j!}\left(\frac{\partial}{\partial{x_i}}\right)^jf(\boldsymbol{x},\boldsymbol{y})
=e^{k\frac{\partial}{\partial{x_i}}}f(\boldsymbol{x},\boldsymbol{y}).$
    
\end{rema}

 
\subsection{Irreducibility}
We first consider the irreducibility of $\mathrm{T}$.
 
\begin{theo}\label{th5.1}
 Suppose that $0\neq v\in V$ and that the pairs $(\lambda_i,\zeta_i)$, $1\le i \le m$, are pairwise distinct. Then,
 \begin{itemize}
\item[\rm(1)]$1^{\otimes m}\otimes v$ generates the module $\mathrm{T}$;
\item[\rm(2)] $\mathrm{T}$ is an irreducible $\mathcal{H}$-module.
\end{itemize}
\end{theo}
\begin{proof}
    \rm(1) Let $\mathrm{M}$ be the nonzero submodule of $\mathrm{T}$ generated by $1^{\otimes m}\otimes v$. Then, for all positive integers $k>K$ and any $l\in\Z$, we have
\begin{equation}\label{formula 5.18}
\aligned
L_{k,l}(1^{\otimes m}\otimes v)
=&\sum_{i=1}^m\lambda_i^k\zeta_i^l1\otimes\dots \otimes\Big(x_i+ky_i-k(k-1)\alpha_i\Big)\otimes\dots \otimes 1\otimes v\\
=&\sum_{i=1}^m\lambda_i^k\zeta_i^l\Big(\boldsymbol{x}_i+k\Big(\boldsymbol{y}_i+\alpha_i1^{\otimes m}\Big)-k^2\alpha_i1^{\otimes m}
\Big)\otimes v.
\endaligned 
\end{equation}
From Lemma \ref{le2.3}, by observing the coefficients of each power of $k$ in (\ref{formula 5.18}), we can deduce that $\boldsymbol{x}_i\otimes v,\boldsymbol{y}_i\otimes v\in \mathrm{M}$ for $i=1,\cdots, m$, which means that $\C[\boldsymbol{x},\boldsymbol{y}]\otimes v\subset \mathrm{M}$.

Let $V_0 := \{u \in V \mid \C[\boldsymbol{x},\boldsymbol{y}] \otimes u \subset \mathrm{M}\}$. The previous argument implies
that $V_0 \neq \{0\}$. Moreover, $X(f(\boldsymbol{x},\boldsymbol{y})\otimes u)=X(f(\boldsymbol{x},\boldsymbol{y}))\otimes u+f(\boldsymbol{x},\boldsymbol{y})\otimes X(u)\in \mathrm{M}$ for all $X\in\mathcal{H},f(\boldsymbol{x},\boldsymbol{y})\in\C[\boldsymbol{x},\boldsymbol{y}]$ and $u\in V_0$, which implies that $V_0$ is an $\mathcal{H}$-submodule of $V$ . Since $V$ is irreducible, we get $V_0=V$ and $\mathrm{M}=\mathrm{T}$.

\rm(2) 
Again, let $\mathrm{M}$ be an arbitrary nonzero submodule of $\mathrm{T}$. Take any nonzero polynomial $f(\boldsymbol{x},\boldsymbol{y})$ in $\mathrm{M}$. Let $\boldsymbol{x}^{\boldsymbol{s}}\boldsymbol{y}^{\boldsymbol{t}}$ be the nonzero term with maximal degree in $\boldsymbol{x}$ and assume $\boldsymbol{0}\prec\boldsymbol{s}$, which means there exists $s^i\ge 1$. For all $k>K$ and $l\in\Z$, we have
\begin{equation*}\label{formula 7.1.3.0}
\aligned
L_{k,l}f(\boldsymbol{x},\boldsymbol{y})
=&\sum_{i=1}^m\lambda_i^k\zeta_i^l\Big(1\otimes\cdots 1\otimes\Bigg(x_i+ky_i-k(k-1)\alpha_i+k(y_i-k\alpha_i)\frac{\partial}{\partial{y_i}}\Bigg)\otimes 1\cdots\otimes 1\Big)    \\            
 &\ \ \ \ \ \Big(\sum_{(\boldsymbol{s},\boldsymbol{t})\in\N^{2m}}x_1^{s_1}y_1^{t_1}x_2^{s_2}y_2^{t_2}\cdots (x_i+k)^{s_i}y_i^{t_i}\cdots x_m^{s_m}y_m^{t_m}\otimes v_{\boldsymbol{s},\boldsymbol{t}}   \Big)  \\
=&\sum_{i=1}^m\lambda_i^k\zeta_i^l\Bigg{(}\boldsymbol{x}_i+k\Big(\boldsymbol{y}_i+\alpha_i1^{\otimes m}\Big)-k^2\alpha_i1^{\otimes m}+k \boldsymbol{y}_i\frac{\partial}{\partial{y_i}}-k^2\alpha_i\frac{\partial}{\partial{y_i}}\Bigg{)}f(\boldsymbol{x}+k\boldsymbol{e}_i,\boldsymbol{y}).
\endaligned 
\end{equation*}
By Remark \ref{rema 7.1}, we can obtain
\begin{equation}\label{formula 5.19}
\aligned
L_{k,l}f(\boldsymbol{x},\boldsymbol{y})
=&\sum_{i=1}^m\sum_{j= 0}^{s^i }\frac{k^j}{j!}\lambda_i^k\zeta_i^l\Bigg(\Bigg(\boldsymbol{x}_i
+k\Bigg(\boldsymbol{y}_i+\alpha_i1^{\otimes m}+\boldsymbol{y}_i\frac{\partial}{\partial{y_i}}  \Bigg)-k^2\alpha_i\Bigg(1^{\otimes m}     
  +\frac{\partial}{\partial{y_i}}\Bigg)\Bigg)\left(\frac{\partial}{\partial{x_i}}     \right)^{j}
\Bigg) f(\boldsymbol{x},\boldsymbol{y}) 
.
\endaligned 
\end{equation}
If $\alpha_i\neq 0$, then extracting the coefficient of $k^{s^i+2}$, we can obtain that
\begin{equation}\label{formula 5.20}
\aligned
\left(\left(\frac{\partial}{\partial{x_i}}     \right)^{s^i} +\frac{\partial}{\partial{y_i}}\left(\frac{\partial}{\partial{x_i}}     \right)^{s^i}\right)f(\boldsymbol{x},\boldsymbol{y})\in \C[x_1,\cdots,x_{i-1},\widehat{x_i},x_{i+1},\cdots,x_m;\boldsymbol{y}]\otimes V,
\endaligned 
\end{equation}
where 
     $\C[x_1,\cdots,x_{i-1},\widehat{x_i},x_{i+1},\cdots,x_m;\boldsymbol{y}]$
 represents that $x_i$ does not appear. If $\alpha_i=0$ and $\gamma_i\ne 0$, similarly, we have
\begin{equation}\label{formula 5.21}
\aligned
I_{r,l}f(\boldsymbol{x},\boldsymbol{y})
=&\sum_{i=1}^m\sum_{j=0}^{s^i}\frac{r^{j}}{j!}\lambda_i^{r-\frac{1}{2}}\zeta_i^l\Bigg{(} \Bigg(  \beta_i1^{\otimes m}+r\gamma_i\Bigg(1^{\otimes m}+\frac{\partial}{\partial{y_i}}\Bigg)\Bigg)\left(\frac{\partial}{\partial{x_i}}     \right)^{j}\Bigg{)}f(\boldsymbol{x},\boldsymbol{y}).
\endaligned 
\end{equation}
for all $r>K$ and $l\in\Z$. Therefore, extracting the coefficient of $r^{s^i+1}$, we can still obtain  (\ref{formula 5.20}). Then $\mathrm{M}\cap (\C[y]\otimes V)\neq \{0\}$. 

Still let $f\in \mathrm{M}\cap (\C[y]\otimes V)$ be a nonzero polynomial with ${\rm deg_{\boldsymbol{x}}}(f)=\boldsymbol{0}$ and minimal ${\rm deg_{\boldsymbol{y}}}(f)$. Let $\boldsymbol{y}^{\boldsymbol{t}}$ be the nonzero term with maximal degree in $\boldsymbol{y}$ and assume $\boldsymbol{0}\prec\boldsymbol{t}$. This implies there exists $t^i\ge 1$. Then (\ref{formula 5.19})  and (\ref{formula 5.21}) can be simplified to
\begin{equation}\label{formula 5.22}
\aligned
L_{k,l}f(\boldsymbol{y})
=&\sum_{i=1}^m\lambda_i^k\zeta_i^l\Bigg(\boldsymbol{x}_i
+k\left(\boldsymbol{y}_i+\alpha_i1^{\otimes m}+\boldsymbol{y}_i\frac{\partial}{\partial{y_i}}  \right)
-k^2\alpha_i\left(1^{\otimes m}+  \frac{\partial}{\partial{y_i}}  \right)
\Bigg)f(\boldsymbol{y})  
\endaligned 
\end{equation}
and
\begin{equation}\label{formula 5.23}
\aligned
I_{r,l}f(\boldsymbol{y})
=&\sum_{i=1}^m\lambda_i^{r-\frac{1}{2}}\zeta_i^l\Bigg{(}\beta_i1^{\otimes m}+r\gamma_i\left(1^{\otimes m}+\frac{\partial}{\partial{y_i}}\right) \Bigg{)} f(\boldsymbol{y}).
\endaligned 
\end{equation}
Observing the quadratic term of $k$ and the linear term of $r$ in (\ref{formula 5.22}) and (\ref{formula 5.23}) respectively, we have 
\[
    \alpha_i\Big(f(\boldsymbol{y})+  \frac{\partial}{\partial{y_i}}f(\boldsymbol{y})  \Big)
 \  ,\gamma_i\Big(f(\boldsymbol{y})+  \frac{\partial}{\partial{y_i}}f(\boldsymbol{y})  \Big) 
\in\mathrm{M}.
\]
Since $(\alpha_i,\gamma_i)\neq (0,0)$, we can obtain $\frac{\partial}{\partial{y_i}}f(\boldsymbol{y})\in\mathrm{M}$. This contradicts the choice of $f$, so $1^{\otimes m}\otimes v\in\mathrm{M}$. 
  Hence,  
  we can get $\mathrm{M}=\mathrm{T}$ and $\mathrm{T}$ is simple.
\end{proof}

From Theorem \ref{th5.1}, we have the following corollary.
\begin{coro}
If  the pairs $(\lambda_i,\zeta_i)$, $1\le i \le m$, are pairwise distinct, then $\bigotimes_{i=1}^{m}\Omega_{\mathcal{H}}(\lambda_i,\zeta_i, \alpha_i,\beta_i,\linebreak\gamma_i)$ is an irreducible $\mathcal{H}$-module.
\end{coro}

\begin{prop}\label{theo5.2}
    Let $V_n$ be the subspace of $\Omega_{\mathcal{H}}(\lambda, \zeta,\alpha_1,\beta_1,\gamma_1)\otimes\Omega_{\mathcal{H}}(\lambda,\zeta, \alpha_2,\beta_2,\gamma_2)$ spanned by $\{x_1^p(x_1+x_2)^qf_1f_2\mid 0\le p\le n, q\in\N,f_1\in\C[y_1],f_2\in\C[y_2]\}$. Then $V_n$ is a proper submodule of $\Omega_{\mathcal{H}}(\lambda, \zeta,\alpha_1,\beta_1,\gamma_1)\otimes\Omega_{\mathcal{H}}(\lambda, \zeta,\alpha_2,\beta_2,\gamma_2)$ for any $n\in\N$.
\end{prop}
\begin{proof}
Set $Z=x_1+x_2$ and, for $i=1,2$, write $f_i'=\frac{d}{dy_i}f_i$. By
(\ref{act. def 7.1}), $L_{k,l}$ and $I_{r,l}$ act on the two tensor factors with
$\lambda_i=\lambda,\ \zeta_i=\zeta$ and the same $(\lambda,\zeta)$. Hence, for
every $0\le p\le n$, $q\in\N$, and $f_1\in\C[y_1], f_2\in\C[y_2]$, a direct
calculation with $h_{i,k}(y_i)=ky_i-k(k-1)\alpha_i$ gives
\begin{equation}\label{theo5.2L}
\begin{aligned}
L_{k,l}\Bigl(x_1^pZ^qf_1f_2\Bigr)
=&\lambda^k\zeta^l(x_1+x_2)x_1^p(x_1+x_2+k)^qf_1f_2+\varepsilon_L,\\
I_{r,l}\Bigl(x_1^pZ^qf_1f_2\Bigr)
=&\lambda^{r-\frac12}\zeta^l\Big(
   (x_1+r)^p(\beta_1+r\gamma_1)+x_1^p(\beta_2+r\gamma_2)\Big)(x_1+x_2+r)^qf_1f_2\\
&+r\gamma_1\lambda^{r-\frac12}\zeta^l(x_1+r)^p(x_1+x_2+r)^qf_1'f_2\\
&+r\gamma_2\lambda^{r-\frac12}\zeta^lx_1^p(x_1+x_2+r)^qf_1f_2'+\varepsilon_I,
\end{aligned}
\end{equation}
where $\varepsilon_L,\varepsilon_I\in V_n$. Indeed, (\ref{theo5.2L}) is a
direct binomial expansion: writing $(x_1+x_2+k)^q
=\sum_{t=0}^q\binom{q}{t}k^{q-t}Z^t$, $(x_1+k)^p
=\sum_{t=0}^p\binom{p}{t}k^{p-t}x_1^t$ and $(x_1+r)^p
=\sum_{t=0}^p\binom{p}{t}r^{p-t}x_1^t$, the terms of maximal $x_1$-degree
combine into the displayed leading parts, and the remainder
$x_1^{p'}Z^{q'}f_1^{(a)}f_2^{(b)}$ satisfies $p'\le n$ for every monomial
in $\varepsilon_L,\varepsilon_I$. (Here $q'\in\N$ and $a,b\in\N$ arise from
the $\frac{\partial}{\partial y_i}$-parts, which never increase the
$x_1$-degree.) Each displayed leading part has $x_1$-degree $p\le n$, so
the right-hand side of (\ref{theo5.2L}) lies in $V_n$. Hence
$\mathcal{H}V_n\subseteq V_n$, and $V_n$ is an $\mathcal{H}$-submodule.
It is proper: by the flat-change-of-variables $Z=x_1+x_2$, every element of
$V_n$ has $x_1$-degree $\le n$, whereas
$T\cong\C[x_1,x_2,y_1,y_2]$ plainly contains, e.g., $x_1^{n+1}\notin V_n$.
\end{proof}

We have the following conclusion immediately from Theorem \ref{th5.1} and Proposition \ref{theo5.2}.

\begin{coro}\label{coro6.6}
  $\mathrm{T}$ is irreducible if and only if the pairs $(\lambda_i,\zeta_i)$, $1\le i \le m$, are pairwise distinct.
\end{coro}


\subsection{Isomorphism} In this subsection, we will determine the necessary and sufficient conditions for two simple tensor product modules of the form $\left(\bigotimes_{i=1}^{m}\Omega_{\mathcal{H}}(\lambda_i, \zeta_i,\alpha_i,\beta_i,\gamma_i)\right)\otimes V$ to be isomorphic.

\begin{theo}\label{theo5.3}
Fix any $m,n\in\Z_+,\{\lambda_i\}_{i=1}^m,\{\zeta_i\}_{i=1}^m,\{\theta_j\}_{j=1}^n,\{\eta_j\}_{j=1}^n\subset\C^*$, $\{\beta_i\}_{i=1}^m,\{\nu_j\}_{j=1}^n\subset\C$ and $\{(\alpha_i, \gamma_i)\}_{i=1}^m,\{(\mu_j, \omega_j)\}_{j=1}^n\subset\C^2\backslash\{(0,0)\}$. The pairs $(\lambda_i,\zeta_i)$, $1\le i \le m$, are pairwise distinct as well as the pairs $(\theta_j,\eta_j)$ $1\le j \le n$, are pairwise distinct. Let $V$ and $V^\prime$ be irreducible restricted $\mathcal{H}$-modules. Then 
\[
    \left(\bigotimes_{i=1}^{m}\Omega_{\mathcal{H}}(\lambda_i, \zeta_i,\alpha_i,\beta_i,\gamma_i)\right)\otimes V\cong 
    \left(\bigotimes_{j=1}^{n}\Omega_{\mathcal{H}}(\theta_j, \eta_j,\mu_j,\nu_j,\omega_j)\right)\otimes V^\prime
\]
if and only if $m=n$, $V\cong V^\prime$ and $(\lambda_i,\zeta_i,\alpha_i,\beta_i, \gamma_i)=(\theta_i,\eta_i,\mu_i,\nu_i ,\omega_i)$ for $1\le i\le m$, after renumbering $(\theta_i,\eta_i,\mu_i, \nu_i,\omega_i)$ if necessary.
\end{theo}
\begin{proof}
The sufficiency is clear. We only need to prove the necessity. 
Set 
\[    \mathrm{T}^\prime=\left(\bigotimes_{j=1}^{n}\Omega_{\mathcal{H}}(\theta_j,\eta_j, \mu_j,\nu_j,\omega_j)\right)\otimes V^\prime\cong\C[\boldsymbol{x}^\prime,\boldsymbol{y}^\prime]\otimes V^\prime
\]
and let $(1^\prime)^{\otimes n}$ denote the multiplicative identity of $\C[\boldsymbol{x^\prime},\boldsymbol{y^\prime}]$. 
Let $\psi:\mathrm{T} \to \mathrm{T}^\prime$ be an $\mathcal{H}$-module isomorphism.
 For any $0\ne v\in V$, write 
 $$0\ne F_{v}=F_v(\boldsymbol{x^\prime},\boldsymbol{y^\prime})=\psi(1^{\otimes m}\otimes v) =\sum_{(\mathbf s,\mathbf t)\in\mathbb N^{2n}} \boldsymbol{x^\prime}^{\mathbf s}\boldsymbol{y^\prime}^{\mathbf t}\otimes v'_{\mathbf s,\mathbf t},$$ 
 where only finitely many $v'_{\mathbf s,\mathbf t}\in V'$ are nonzero. And denote $s^i={\rm deg}_{x'_i}(F_v)={\rm max}\{s_i\mid v'_{\mathbf s,\mathbf t}\ne0\}$. Since $V'$ is restricted, 
for any $v'\in V'$, there exists $K'\in\N$ such that $L_{k,l}v'=I_{r,l}v'=0$ for all $k,r>K'$.
By (\ref{act. def 7.1}), $\psi(L_{k,l}(1^{\otimes m}\otimes v))=L_{k,l}(\psi(1^{\otimes m}\otimes v))$  for all $ k>{\rm max}\{ K,K'\},l\in\mathbb Z$,
we have 
\begin{equation}\label{formula 5.24}
\begin{aligned}
&\sum_{i=1}^m\lambda_i^k\zeta_i^l\Big(\psi(\boldsymbol{x}_i\otimes v)
+k\Big(\psi(\boldsymbol{y}_i\otimes v)+\alpha_iF_v  \Big)-k^2\alpha_iF_v     
  \Big)
\\
=&\sum_{j=1}^n\sum_{p= 0}^{s^j}\frac{k^p}{p!}\theta_j^k\eta_j^l\Bigg(\Bigg(\boldsymbol{x^\prime}_j+k\Bigg(\boldsymbol{y^\prime}_j+\mu_j(1')^{\otimes n}+\boldsymbol{y^\prime}_j\frac{\partial}{\partial{y'_j}}  \Bigg)-k^2\mu_j\Bigg((1')^{\otimes n} 
  +\frac{\partial}{\partial{y'_j}}\Bigg)\Bigg)\left(\frac{\partial}{\partial{x'_j}}     \right)^{p}
\Bigg) F_v.
\end{aligned}
\end{equation}
Applying Lemma \ref{le2.3},
if some $(\theta_j,\eta_j)$ does not occur among $(\lambda_i,\zeta_i)$, then the polynomial corresponding to $\theta_j^k\eta_j^l$ on the right side must vanish. Its constant term in $k$ is $\boldsymbol{x^\prime}_jF_v\ne0$, which is a contradiction. Hence $\{(\theta_j,\eta_j)\mid 1\le j\le n\} \subseteq \{(\lambda_i,\zeta_i)\mid 1\le i\le m\}$. Applying the same argument to $\psi^{-1}$, we obtain the reverse inclusion. Since the pairs in each family are pairwise distinct, it follows that $m=n$, and, after renumbering if necessary, $\lambda_i=\theta_i, \zeta_i=\eta_i,  1\le i\le m$. Therefore, we can directly verify that
\begin{equation}\label{formula 5.25}
\aligned
\psi(\boldsymbol{x}_i\otimes v)=\boldsymbol{x^\prime}_iF_v=\boldsymbol{x^\prime}_i\psi(1^{\otimes m}\otimes v).
\endaligned 
\end{equation}
For all $r>{\rm max}\{ K,K'\},l\in\mathbb Z$
and $\psi(I_{r,l}(1^{\otimes m}\otimes v))=I_{r,l}(\psi(1^{\otimes m}\otimes v))$, we have
\begin{equation}\label{formula 5.26}
\aligned
&\sum_{i=1}^m\lambda_i^{r-\frac{1}{2}}\zeta_i^l (\beta_i+r\gamma_i)F_v
=\sum_{i=1}^m\sum_{p=0}^{s^i}\frac{r^{p}}{p!}\lambda_i^{r-\frac{1}{2}}\zeta_i^l\Bigg{(} \Bigg(  \nu_i1^{\otimes m}+r\omega_i\Bigg(1^{\otimes m}+\frac{\partial}{\partial{y'_i}}\Bigg)\Bigg)\left(\frac{\partial}{\partial{x'_i}}     \right)^{p}\Bigg{)}F_v.
\endaligned 
\end{equation}
Similarly, we can obtain that $\nu_i=\beta_i$.

Suppose, for a contradiction, that $s^i>0$ for some $i$.
Write $$F_v(\boldsymbol{x^\prime},\boldsymbol{y^\prime})=\sum_{p=0}^{s^i}\boldsymbol{x'}_i^pf_p(\boldsymbol{x'},\boldsymbol{y'})\in\C[\boldsymbol{x'},\boldsymbol{y'}]\otimes V' 
$$
 where $f_{p}(\boldsymbol{x^\prime},\boldsymbol{y^\prime})\in\C[x'_1,\cdots,x'_{i-1},\widehat{x'_i},x'_{i+1},\cdots,x'_m;\boldsymbol{y'}]\otimes V', f_{s^i}\neq0$ 
 . If $\mu_i\ne 0$, by comparing the coefficients of $k^{s^i+2}$ in (\ref{formula 5.24}), we have
\begin{equation}\label{formula 5.27}
\begin{aligned}
0=-\mu_i\Bigg( f_{s^i}+\frac{\partial}{\partial{y'_i}}f_{s^i}    \Bigg)\neq0.
\end{aligned}
\end{equation}
If $\mu_i=0$ and $\omega_i\ne0$, 
by comparing the coefficients of $r^{s^i+1}$ in (\ref{formula 5.26}), we have
\begin{equation}\label{formula 5.28}
\begin{aligned}
0=\omega_i\Bigg( f_{s^i}+\frac{\partial}{\partial{y'_i}}f_{s^i}    \Bigg)\neq0.
\end{aligned}
\end{equation}
Therefore, we always derive contradiction (\ref{formula 5.27}) or (\ref{formula 5.28}), which implies that $s^i=0$ for all $i=1,\cdots,m$.  Then $F_v(\boldsymbol{x^\prime},\boldsymbol{y^\prime})=F_v(\boldsymbol{y^\prime})\in\C[\boldsymbol{y^\prime}]\otimes V'$ and comparing the coefficients of $k,k^2$ and $r$ in (\ref{formula 5.24}) and (\ref{formula 5.26}), respectively, we obtain that

\begin{equation}\label{formula 5.29}
\aligned
\psi(\boldsymbol{y}_i\otimes v)=\boldsymbol{y^\prime}_iF_v+(\mu_i-\alpha_i)F_v+\boldsymbol{y^\prime}_i\frac{\partial}{\partial{y'_i}}F_v,
\endaligned 
\end{equation}
\begin{equation}\label{formula 5.30}
\aligned
\alpha_iF_v=\mu_i\Bigg(1^{\otimes m}+\frac{\partial}{\partial{y'_i}}\Bigg)F_v,
\endaligned 
\end{equation}
\begin{equation}\label{formula 5.31}
\aligned
\gamma_iF_v=\omega_i\Bigg(1^{\otimes m}+\frac{\partial}{\partial{y'_i}}\Bigg)F_v
\endaligned 
\end{equation}
for any $i=1,\cdots,m$.
 Due to $(\alpha_i,\gamma_i)\neq(0,0)$, observing  the coefficients of the maximal degree terms on both sides of (\ref{formula 5.30}) and (\ref{formula 5.31}) in $y'_i$, we can deduce that $\mu_i=\alpha_i$,$\omega_i=\gamma_i$ and 
 ${\rm deg}_{\boldsymbol{y^\prime}}(F_v(\boldsymbol{y^\prime}))=0$.
Hence, 
 we have  
\begin{equation}\label{formula 5.32}
\aligned
\psi(\boldsymbol{y}_i\otimes v)=\boldsymbol{y^\prime}_iF_v=\boldsymbol{y^\prime}_i\psi(1^{\otimes m}\otimes v)
\endaligned 
\end{equation}
through (\ref{formula 5.29}). ${\rm deg}_{\boldsymbol{y^\prime}}(F_v(\boldsymbol{y^\prime}))=0$ also means that
 there exists a unique $\Gamma(v)\in V^\prime$ such that
  $F_v=(1^\prime)^{\otimes m}\otimes \Gamma(v)$, where $\Gamma$ is a vector space homomorphism from $V$ to $V^\prime$. Combining (\ref{formula 5.25}) and (\ref{formula 5.32}), for  any $X\in\{L_{k,l}, I_{r,l}, C_{j,l} \mid k,l\in \mathbb{Z}, r\in\mathbb{Z}+\tfrac{1}{2}, j=1,3  \}$ and $v\in V$, we can obtain that $\psi\Big(X\Big(1^{\otimes m}\Big)\otimes v  \Big)=X\Big( (1^\prime)^{\otimes m} \Big)\otimes \Gamma(v)$
and
\begin{equation*}\label{formula 5.4.2.3}
\aligned
(1^\prime)^{\otimes m}\otimes \Gamma(X(v))=&\psi\Big( 1^{\otimes m}\otimes  X(v) \Big)
=\psi\Big(X\Big( 1^{\otimes m}\otimes  v \Big)-  X\Big(1^{\otimes m} \Big) \otimes v\Big)\\
=&X\Big( (1^\prime)^{\otimes m}\otimes  \Gamma(v) \Big)-  X\Big((1^\prime)^{\otimes m} \Big) \otimes \Gamma(v)
=(1^\prime)^{\otimes m}\otimes X(\Gamma(v)).
\endaligned 
\end{equation*}
Then $\Gamma(X(v))=X(\Gamma(v))$ for any $v\in V$. Since $\psi$ is injective and $\Gamma$ is nonzero,  $\Gamma$ is a $\mathcal{H}$-module  isomorphism, i.e.  $V\cong V^\prime$.

This completes the proof of the theorem.
\end{proof}

Based on the previous discussion, we conclude with the following Remark.
\begin{rema}

Let  $\rho_{\mathrm{T}}:\mathcal H\to\operatorname{End}_{\mathbb C}(\C[\boldsymbol{x},\boldsymbol{y}]\otimes V)$ and $\rho_{V}:\mathcal H\to\operatorname{End}_{\mathbb C}(V)$ be the corresponding representations of $\mathcal H$-module $\mathrm{T}$ and $V$, respectively. By Remark \ref{rema 7.1}, the $\mathcal H$-action on $\mathrm{T}$ admits the equivalent realization
\begin{equation*}\label{formula 5.4.2.4}
\aligned
\rho_{\mathrm{T}}(L_{k,l})
&=\sum_{i=1}^{m}\lambda_i^k\zeta_i^l
\left(
\boldsymbol{x}_i+k\boldsymbol{y}_i-k(k-1)\alpha_i1^{\otimes m}
+k(\boldsymbol{y}_i-k\alpha_i1^{\otimes m})\frac{\partial}{\partial{y_i}}
\right)
e^{k\frac{\partial}{\partial{x_i}}}\otimes\operatorname{Id}_V\\
&\quad+\operatorname{Id}_{\C[\boldsymbol{x},\boldsymbol{y}]}\otimes\rho_V(L_{k,l}),\\[1mm]
\rho_{\mathrm{T}}(I_{r,l})
&=\sum_{i=1}^{m}\lambda_i^{r-\frac12}\zeta_i^l
\left(
\beta_i1^{\otimes m}+r\gamma_i1^{\otimes m}+r\gamma_i\frac{\partial}{\partial{y_i}}
\right)
e^{r\frac{\partial}{\partial{x_i}}}\otimes\operatorname{Id}_V\\
&\quad+\operatorname{Id}_{\C[\boldsymbol{x},\boldsymbol{y}]}\otimes\rho_V(I_{r,l}),\\
\rho_{\mathrm{T}}(C_{j,l})
&=\operatorname{Id}_{\C[\boldsymbol{x},\boldsymbol{y}]}\otimes \rho_V(C_{j,l}),
\endaligned 
\end{equation*}
where $k,l\in\Z,r\in\frac{1}{2}+\Z,j=1,3$.

\end{rema}

\appendix

\section{The proof of Proposition \ref{prop 3.1}}\label{ Proposition 3.1}
\begin{proof}

    For any $m,n,i,j\in\mathbb{Z}$ and $f(x,y)\in\mathbb{C}[x,y]$, we check
\begin{equation*}
\aligned
L_{m,i}W_{n,j}f(x,y)
=&\lambda^{m+n}\zeta^{i+j}
\bigl(x+h_m(y)+i\sigma(y-m\alpha)\bigr)
(y-n\alpha)f(x+m+n,y)\\
&+m\lambda^{m+n}\zeta^{i+j}
(y-m\alpha)f(x+m+n,y)\\
&+m\lambda^{m+n}\zeta^{i+j}
(y-m\alpha)(y-n\alpha)
\frac{\partial }{\partial y}f(x+m+n,y),
\endaligned
\end{equation*}
\begin{equation*}
\aligned
W_{n,j}L_{m,i}f(x,y)
=&\lambda^{m+n}\zeta^{i+j}(y-n\alpha)
\bigl(x+n+h_m(y)+i\sigma(y-m\alpha)\bigr)
f(x+m+n,y)\\
&+m\lambda^{m+n}\zeta^{i+j}
(y-n\alpha)(y-m\alpha)
\frac{\partial }{\partial y}f(x+m+n,y),
\endaligned
\end{equation*}

Then, we have
\begin{equation*}
\aligned
\relax[L_{m,i},W_{n,j}]f(x,y)
=&\lambda^{m+n}\zeta^{i+j}
\bigl(m(y-m\alpha)-n(y-n\alpha)\bigr)f(x+m+n,y)\\
=&(m-n)\lambda^{m+n}\zeta^{i+j}
(y-(m+n)\alpha)f(x+m+n,y)\\
=&(m-n)W_{m+n,i+j}f(x,y)+\frac{m^3-m}{12}\delta_{m+n,0}C_{2,i+j}f(x,y).
\endaligned
\end{equation*}

For any $m,i,j\in\mathbb{Z}, r\in\Z+\frac{1}{2}, f(x,y)\in\mathbb{C}[x,y]$, direct calculation gives
\begin{equation*}
\aligned
L_{m,i}I_{r,j}f(x,y)
&=\lambda^{m+r-\frac12}\zeta^{i+j}
\bigl(x+h_m(y)+i\sigma(y-m\alpha)\bigr)
\bigl(\beta+\gamma g_r(y)+j\sigma\gamma\bigr)f(x+m+r,y)\\
&\quad+m\gamma\lambda^{m+r-\frac12}\zeta^{i+j}
(y-m\alpha)\frac{d}{dy}g_r(y)f(x+m+r,y)\\
&\quad+m\lambda^{m+r-\frac12}\zeta^{i+j}
(y-m\alpha)\bigl(\beta+\gamma g_r(y)+j\sigma\gamma\bigr)
\frac{\partial }{\partial y}f(x+m+r,y)\\
&\quad+r\gamma\lambda^{m+r-\frac12}\zeta^{i+j}
\bigl(x+h_m(y)+i\sigma(y-m\alpha)\bigr)
\frac{\partial }{\partial y}f(x+m+r,y)\\
&\quad+mr\gamma\lambda^{m+r-\frac12}\zeta^{i+j}
(y-m\alpha)\frac{\partial^2}{\partial y^2}f(x+m+r,y),
\endaligned
\end{equation*}
\begin{equation*}
\aligned
I_{r,j}L_{m,i}f(x,y)
&=\lambda^{m+r-\frac12}\zeta^{i+j}
\bigl(x+r+h_m(y)+i\sigma(y-m\alpha)\bigr)
\bigl(\beta+\gamma g_r(y)+j\sigma\gamma\bigr)f(x+m+r,y)\\
&\quad+r\gamma\lambda^{m+r-\frac12}\zeta^{i+j}
\left(\frac{d}{dy}h_m(y)+i\sigma\right)f(x+m+r,y)\\
&\quad+r\gamma\lambda^{m+r-\frac12}\zeta^{i+j}
\bigl(x+r+h_m(y)+i\sigma(y-m\alpha)\bigr)
\frac{\partial }{\partial y}f(x+m+r,y)\\
&\quad+m\lambda^{m+r-\frac12}\zeta^{i+j}
(y-m\alpha)\bigl(\beta+\gamma g_r(y)+j\sigma\gamma\bigr)
\frac{\partial }{\partial y}f(x+m+r,y)\\
&\quad+mr\gamma\lambda^{m+r-\frac12}\zeta^{i+j}
\frac{\partial }{\partial y}f(x+m+r,y)\\
&\quad+mr\gamma\lambda^{m+r-\frac12}\zeta^{i+j}
(y-m\alpha)\frac{\partial^2}{\partial y^2}f(x+m+r,y).
\endaligned
\end{equation*}
Subtracting these expressions and using (\ref{formula 2.4}), we obtain
\begin{equation*}
\aligned
\relax[L_{m,i},I_{r,j}]f(x,y)
=&-r\lambda^{m+r-\frac12}\zeta^{i+j}
\bigl(\beta+\gamma g_{m+r}(y)+(i+j)\sigma\gamma\bigr)
f(x+m+r,y)\\
&-r(m+r)\gamma\lambda^{m+r-\frac12}\zeta^{i+j}
\frac{\partial }{\partial y}f(x+m+r,y)
=-rI_{m+r,i+j}f(x,y).
\endaligned
\end{equation*}
Similarly, we have
\begin{equation*}
\aligned
L_{m,i}B_{r,j}f(x,y)
=&\gamma\lambda^{m+r-\frac12}\zeta^{i+j}
\bigl(x+h_m(y)+i\sigma(y-m\alpha)\bigr)f(x+m+r,y)\\
&\quad+m\gamma\lambda^{m+r-\frac12}\zeta^{i+j}
(y-m\alpha)\frac{\partial }{\partial y}f(x+m+r,y),\\
B_{r,j}L_{m,i}f(x,y)
=&\gamma\lambda^{m+r-\frac12}\zeta^{i+j}
\bigl(x+r+h_m(y)+i\sigma(y-m\alpha)\bigr)f(x+m+r,y)\\
&\quad+m\gamma\lambda^{m+r-\frac12}\zeta^{i+j}
(y-m\alpha)\frac{\partial }{\partial y}f(x+m+r,y)
\endaligned
\end{equation*}
and
\begin{equation*}
\aligned
\relax[L_{m,i},B_{r,j}]f(x,y)
=-r\gamma\lambda^{m+r-\frac12}\zeta^{i+j}f(x+m+r,y)
=-rB_{m+r,i+j}f(x,y)
\endaligned
\end{equation*}
for all $m,i,j\in\mathbb{Z}, r\in\Z+\frac{1}{2}, f(x,y)\in\mathbb{C}[x,y]$.

For any $m,i,j\in\mathbb{Z}, r\in\Z+\frac{1}{2}, f(x,y)\in\mathbb{C}[x,y]$, we can obtain
\begin{equation*}
\aligned
W_{m,i}I_{r,j}f(x,y)
=&\lambda^{m+r-\frac12}\zeta^{i+j}
(y-m\alpha)\bigl(\beta+\gamma g_r(y)+j\sigma\gamma\bigr)
f(x+m+r,y)\\
&+r\gamma\lambda^{m+r-\frac12}\zeta^{i+j}
(y-m\alpha)\frac{\partial }{\partial y}f(x+m+r,y),\\
I_{r,j}W_{m,i}f(x,y)
=&\lambda^{m+r-\frac12}\zeta^{i+j}
(y-m\alpha)\bigl(\beta+\gamma g_r(y)+j\sigma\gamma\bigr)
f(x+m+r,y)\\
&+r\gamma\lambda^{m+r-\frac12}\zeta^{i+j}f(x+m+r,y)\\
&+r\gamma\lambda^{m+r-\frac12}\zeta^{i+j}
(y-m\alpha)\frac{\partial }{\partial y}f(x+m+r,y),
\endaligned
\end{equation*}
\begin{equation*}
\aligned
\relax[W_{m,i},I_{r,j}]f(x,y)
=&-r\gamma\lambda^{m+r-\frac12}\zeta^{i+j}f(x+m+r,y)=-rB_{m+r,i+j}f(x,y).
\endaligned
\end{equation*}

For any $i,j\in\Z,r,s\in\Z+\frac{1}{2}, f(x,y)\in\mathbb{C}[x,y]$, we verify
\begin{equation*}
\aligned
I_{r,i}I_{s,j}f(x,y)
=&\lambda^{r+s-1}\zeta^{i+j}
\bigl(\beta+\gamma g_r(y)+i\sigma\gamma\bigr)
\bigl(\beta+\gamma g_s(y)+j\sigma\gamma\bigr)f(x+r+s,y)\\
&+r\gamma^2\lambda^{r+s-1}\zeta^{i+j}
\frac{d}{dy}g_s(y)f(x+r+s,y)\\
&+s\gamma\lambda^{r+s-1}\zeta^{i+j}
\bigl(\beta+\gamma g_r(y)+i\sigma\gamma\bigr)
\frac{\partial }{\partial y}f(x+r+s,y)\\
&+r\gamma\lambda^{r+s-1}\zeta^{i+j}
\bigl(\beta+\gamma g_s(y)+j\sigma\gamma\bigr)
\frac{\partial }{\partial y}f(x+r+s,y)\\
&+rs\gamma^2\lambda^{r+s-1}\zeta^{i+j}
\frac{\partial^2}{\partial y^2}f(x+r+s,y),
\endaligned
\end{equation*}
\begin{equation*}
\aligned
I_{r,i}B_{s,j}f(x,y)
=&\gamma\lambda^{r+s-1}\zeta^{i+j}
\bigl(\beta+\gamma g_r(y)+i\sigma\gamma\bigr)f(x+r+s,y)\\
&+r\gamma^2\lambda^{r+s-1}\zeta^{i+j}
\frac{\partial }{\partial y}f(x+r+s,y)
=B_{s,j}I_{r,i}f(x,y).
\endaligned
\end{equation*}
Therefore, we have
\begin{equation*}
\aligned
\relax[I_{r,i},I_{s,j}]f(x,y)
=\gamma^2\lambda^{r+s-1}\zeta^{i+j}
\left(r\frac{d}{dy}g_s(y)-s\frac{d}{dy}g_r(y)\right)
f(x+r+s,y)
=r\delta_{r+s,0}C_{3,i+j}f(x,y),
\endaligned
\end{equation*}
\[
[I_{r,i},B_{s,j}]f(x,y)
=0=r\delta_{r+s,0}C_{4,i+j}f(x,y).
\]
The remaining relations  are easy to obtain,  so we omit
 the details.
This completes the verification of Proposition \ref{prop 3.1}.
\end{proof}

\section*{Acknowledgements}
This work was supported by the National Natural Science Foundation of China (No. 12171132). The authors would like to express sincere gratitude to Professor Naihuan Jing for the invaluable discussions and constructive suggestions throughout this work.

\section*{Data Availability Statement}
This manuscript has no associated data.

\section*{Conflicts of Interest} The authors declare that they have no conflicts of interest in this work.

\end{document}